\documentclass{amsart}

\usepackage{amssymb,amsfonts}
\usepackage[dvips]{graphicx}

\theoremstyle{plain}
\newtheorem{thm}{Theorem}[section]
\newtheorem{cor}[thm]{Corollary}
\newtheorem{lem}[thm]{Lemma}
\newtheorem{prop}[thm]{Proposition}
\newtheorem{conj}[thm]{Conjecture}
\newtheorem{main}[thm]{Main~Theorem}

\newtheorem{thmsub}{Theorem}[subsection]
\newtheorem{corsub}[thmsub]{Corollary}
\newtheorem{lemsub}[thmsub]{Lemma}
\newtheorem{propsub}[thmsub]{Proposition}

\theoremstyle{definition}
\newtheorem{defn}{Definition}
\theoremstyle{remark}

\newtheorem{rem}[thm]{Remark}
\newtheorem{remsub}[thmsub]{Remark}

\numberwithin{equation}{section}

 \newcommand{\Stw}{\textbf{S}$^2$}
 \newcommand{\Rtw}{\textbf{R}$^2$}

 \newcommand{\Rth}{\textbf{R}$^3$}

 \newcommand{\sixty}{\frac{\pi}{3}}
 \newcommand{\thirty}{\frac{\pi}{6}}
 \newcommand{\onetwenty}{120^{\circ}}
 
 \newcommand{\half}{\frac{1}{2}}

 \newcommand{\bl}{band lens}

 \newcommand{\sdb}{standard double bubble}

\newcommand{\Path}{.}

\newcommand{\mbf}{\mathbf}

\begin{document}

%ARTICLE INFO
\title[The Double Bubble Problem on the Flat Two-Torus]
      {The Double Bubble Problem on the Flat Two-Torus}
\date{\today}
\author[J. Corneli]{Joseph Corneli}
\author[P. Holt]{Paul Holt}
\author[G. Lee]{George Lee}
\author[N. Leger]{Nicholas Leger}
\author[E. Schoenfeld]{Eric Schoenfeld}
\author[B. Steinhurst]{Benjamin Steinhurst}

\address{Mailing address: C/O Frank Morgan\\
         Department of Mathematics and Statistics\\
         \indent Williams College\\
         Williamstown, MA 01267}
\email{Frank.Morgan@williams.edu}

\address{Joseph Corneli and Nicholas Leger\\
          Department of Mathematics \\
          University of Texas\\
          \indent Austin, TX 78712}
\email{jcorneli@math.utexas.edu}
\email{nickleger@mail.utexas.edu}

\address{Paul Holt, Eric Schoenfeld, and Benjamin Steinhurst\\
         Department of Mathematics \indent and Statistics \\ Williams College\\
         Williamstown, MA 01267}
\email{pholt@wso.williams.edu} \email{eschoenf@wso.williams.edu}
\email{Benjamin.A.Steinhurst@williams.edu}

\address {George Lee \\
          Department of Mathematics \\
          Harvard University\\
          Cambridge, MA 02138}
\email{lee43@fas.harvard.edu}

\begin{abstract}
We characterize the perimeter-minimizing double bubbles on all flat two-tori and, as corollaries, on the flat infinite cylinder and the flat infinite strip with free boundary.  Specifically, we show that there are five distinct types of minimizers on flat two-tori, depending on the areas to be enclosed. 
\end{abstract}

\maketitle
\tableofcontents

\footnotetext[1]{2000 Mathematics Subject Classification 53A10(49Q10)}

\section{Introduction} \label{S:intro}
\begin{figure}[ht]
\begin{center}
\begin{tabular} {@{\extracolsep{1.0 cm}} cc}
\includegraphics*{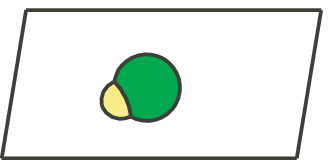}%
&
\includegraphics*{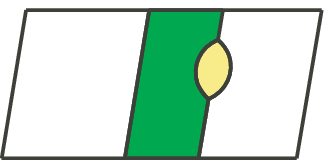}%
\\ \medskip
\textsc{Standard Double Bubble} & \textsc{Band Lens} \\
\includegraphics*{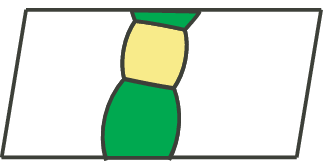}%
&
\includegraphics*{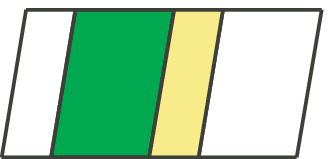}%
\\  \medskip
\textsc{Standard Chain} & \textsc{Double Band} \\
\end{tabular}\\

\includegraphics*[height=.7in]{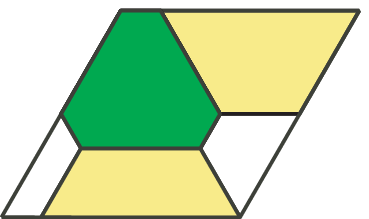}%
\\
\textsc{Standard Hexagon Tiling}

\end{center}

\caption{\label{fig:introduce_minimizers}There are five types of
perimeter-minimizing double bubbles on a flat two-torus. The
standard hexagon tiling occurs only on the hexagonal torus.}
\end{figure}

\begin{figure}[htbp]
\begin{center}
\includegraphics*[height=3.25in]{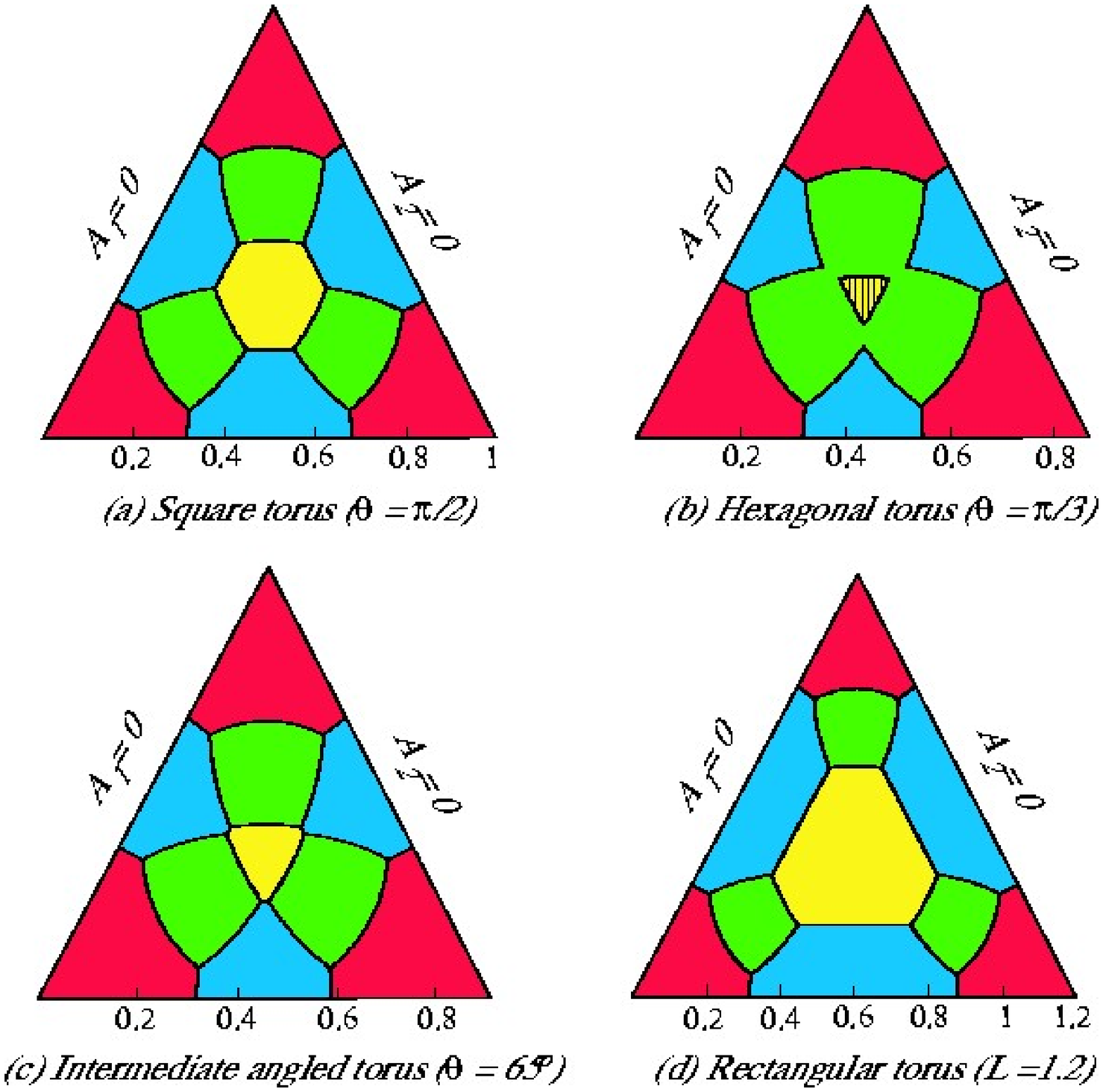} \\
\bigskip{}
\includegraphics*[height=.35in]{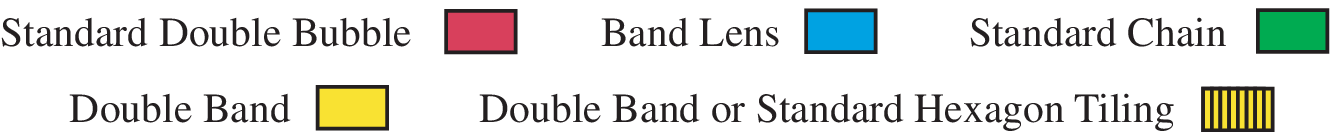}
\end{center}
\caption{\label{fig:param-diags} These phase portraits show which
type of double bubble is minimizing given prescribed area pairs
$(A_1, A_2)$, on four specific flat tori. Each torus for (a)-(c) has a rhombic fundamental domain with interior angle $\theta$; the torus for (d) has a rectangular fundamental domain with side lengths $1$ and $L = 1.2$.}
\end{figure}

Our Main Theorem \ref{thm:main theorem} shows that on any flat
two-torus, the least-perimeter way to enclose and separate two
prescribed areas is a double bubble of one of the five types shown
in Figure \ref{fig:introduce_minimizers}: the standard double
bubble, the band lens, a standard chain, the double band, or the
standard hexagon tiling. Figure \ref{fig:param-diags} gives
``phase portraits,'' computational plots of the minimizing type
for given areas, for four different tori. (Section \ref{S:
numerical comparisons} describes the creation of these phase
portraits, using formulas for the perimeter and area of the five
minimizers derived in Section \ref{S:perimeter-and-area}.)

%%%%%%%%%%%%%%%%%%%%%%%%%%%%%%%%%%%%%%%%%%%%%%%%%%%%%%%%%%%%
%%%%%%%%%%%%%%%%%%%%%%%%%%%%%%%%%%%%%%%%%%%%%%%%%%%%%%%%%%%%

%%%%%%%%%%%%%%%%%%%%%%%%%%%%%%%%%%%%%%%%%%%%%%%%%%%%%%%%%%%%
%%%%%%%%%%%%%%%%%%%%%%%%%%%%%%%%%%%%%%%%%%%%%%%%%%%%%%%%%%%%

The strategy of the proof is to separate all candidates into five
classes, as described in Proposition \ref{prop:topo
classification}. Proposition \ref{prop:topo classification} says
that a minimizer must either be the double band or (possibly after
relabelling the two interior regions and the exterior) belong to
one of four other topological classes: Section
\ref{SS:contract_bcs} deals with contractible double bubbles;
Section \ref{SS:swaths} deals with double bubbles whose components
taken together wrap around the torus, and whose complement is not
contractible; Section \ref{SS:double_bubbles_with_single_band}
deals with double bubbles with one band adjacent to a contractible
set of components; and Section \ref{SS:tilings} deals with double
bubbles for which both enclosed regions and the exterior region are contractible.

The last of these topological classes (with both regions and the
exterior contractible) we refer to as ``tilings" because the three
regions taken together lift to a tiling of the plane.  Despite
escaping our initial conjectures, this class turned out to be one
of the most interesting. Whereas the non-tiling types tie only for
prescribed areas in the transitions between phases (Figure
\ref{fig:param-diags}), the hexagon tiling and the double band
provide distinct minimizers over an open set of areas (Figure
\ref{fig:param-diags}(b)) -- the first example of such a
phenomenon in double bubble history.  Another candidate tiling,
the octagon square, with one region consisting of two
four-sided components (Figure \ref{fig:oct-square}), was the
last candidate to be eliminated (Proposition \ref{prop:oct-square
bad}).

We also consider the double bubble problem on the flat infinite
cylinder, which is simpler because there are no candidate tilings.
Corollary \ref{cor:cylinder corollary} provides a complete
characterization of the minimizers in this space as a trivial
consequence of the torus result, and the cylinder result is used in Corollary \ref{cor:free-boundary} to characterize all minimizers on the flat infinite strip with free boundary.

The torus and cylinder lack the useful symmetries used in other
spaces to establish connectedness of the exterior and bounds on
the number of components.  Luckily, the recent work of Wichiramala
[MW] on the triple bubble problem in the plane gave us a very
useful component
bound (Proposition \ref{prop:component_bound}). \\

\noindent \paragraph{\textbf{History and recent developments}}

The double bubble problem is a generalization of the isoperimetric
problem. Given a geometric space, the latter seeks the least-area
way to enclose a single volume.  (A solution to the isoperimetric problem on the flat two-torus can be found in [HHM].)  The double bubble problem seeks
the least-area way to enclose and separate two volumes. In recent
years, this problem has been solved in a number of
spaces, including \Rtw\ [F], \Stw\ [Ma], and still more recently \Rth\ [HMRR] and $\mathbf{R}^4$ [RHLS].
In these cases the minimizers are all standard double bubbles, in contrast
to the multiplicity of types on the torus. Some of the major contributions were due to
undergraduates affiliated with Williams College and the Williams
``SMALL" REU program.  Underlying all of these results is geometric
measure theory as developed by Federer, Fleming, Almgren, and
others.  For a more complete discussion of the history of the double bubble
problem as well as many related problems on films, foams, and
other efficient surfaces, see [M].

The triple bubble problem in the plane has recently been solved by Wichiramala [W]. On
the other hand, even the single bubble problem remains unsolved on a torus of revolution. We have contributed to a paper that gives experimental
evidence on the solutions to the double bubble problem in a cubic flat three-torus
[CCWB], another space where the single bubble problem is still
open.

\bigskip

\paragraph{\textbf{Acknowledgements}}

This is work of the Williams College National Science Foundation
``SMALL" undergraduate research Geometry Group, advised by Frank
Morgan. The authors would like to thank first of all Professor
Morgan, whose patient advising and perseverance have managed to
always keep us moving in the right direction.  We also owe our
thanks to the National Science Foundation and to Williams College
for sponsoring the SMALL REU.  In addition, we would like very
much to thank Joel Hass and David Hoffman, organizers of the
Clay/MSRI Summer School on the Global Theory of Minimal Surfaces,
as well as MSRI and the Clay Institute for bravely supporting four
of the authors during what happily turned out to be several very
productive weeks in Berkeley in July, 2001.  We are very grateful
to Gary Lawlor, who  came up with the idea for showing that
octagon-square tilings are not minimizers.  We also would like to
thank John Sullivan and Michael Hutchings for helpful
conversations about tilings on the torus and Wacharin Wichiramala
for helping us understand area-preserving variations in hexagon
tilings. Joseph Masters [Ma1] made the first stab at this problem,
and we were grateful for his notes.

\section{Definitions}

This section contains definitions for, among other
things, the classes of potential minimizers mentioned above.

\subsection{Tori and the Cylinder}

\begin{defn} \label{defn:torus}
A \emph{flat two-dimensional torus} can be represented by a planar
parallelogram with opposite sides identified. Throughout the paper
the torus is normalized so that each shortest closed geodesic has
unit length. \end{defn}

\begin{remsub} \label{rem:torus-angle-and-area} The parallelogram can be chosen such that one shortest side
is a shortest closed geodesic and no interior angle is less than
sixty degrees. Then the area of the parallelogram, and hence the
area of the torus, is at least $\sqrt{3}/2$ with equality
precisely for the hexagonal torus (see Definition
\ref{defn:hexagonal-torus}).
\end{remsub}

\begin{defn} \label{defn:hexagonal-torus}
The \emph{hexagonal torus} can be represented either by a regular
hexagon with opposite sides identified with the same orientation,
or by a 60-degree rhombus with opposite sides identified with the
same orientation.
\end{defn}

\begin{defn}
We say that each shortest closed geodesic lies along a
\textit{short direction} of the torus. On certain nonhexagonal
tori (those that can be represented by a parallelogram with four
sides of unit length), there are two short directions; on the
hexagonal torus, there are three.
\end{defn}

\begin{defn}
The \emph{flat infinite cylinder} can be represented as the
surface contained between two identified infinite parallel lines
in the plane, where each segment perpendicular to and connecting
the lines is a closed geodesic. The \emph{flat infinite strip with free boundary} is the surface between the same two lines, except that the lines are not identified; portions of the lines may be included in the boundary of a double bubble without contributing to the perimeter of the double bubble.
\end{defn}

\subsection{Topological Definitions}

\begin{defn}
A \emph{double bubble} on a smooth Riemannian surface consists of two disjoint regions (i.e.,
nonempty open sets) bounded by piecewise smooth curves, such that
the exterior of these regions has nonzero area.  The boundary of
the double bubble divides the surface into three regions: the two
regions belonging to the double bubble, which we call the
\textit{enclosed} regions, and the complement of the closure of
the double bubble, which we call the \textit{exterior} region. Any
of the three regions may have multiple components.
\end{defn}

\begin{remsub} On the torus, by relabelling which of the three regions in a double bubble is
the exterior, we can obtain two other double bubbles with the same
boundary. On the infinite cylinder, the exterior region has infinite area.
\end{remsub}

\begin{defn}
A \emph{tiling} is a double bubble on the torus for which each
component, including each component of the exterior, is
contractible, such as the standard hexagon tiling of Figure
\ref{fig:introduce_minimizers}.  It is called a tiling because the
two enclosed regions plus the exterior region lift to a tiling of
the plane.
\end{defn}

\begin{defn} \label{defn:octagon-square}
An \emph{octagon-square tiling} is a tiling in which one of the
three regions consists of two curvilinear quadrilaterals, and each
of the other two regions consists of one curvilinear octagon (see
Figure \ref{fig:oct-square}).
\end{defn}

\begin{defn}
A \emph{band} is a non-contractible annulus.
\end{defn}

\begin{defn} \label{defn:definition swath}
A \emph{swath} is a set of contractible components  with
non-contractible complement and whose union has non-contractible
and connected closure, such as the standard chain of Figure
\ref{fig:introduce_minimizers}.
\end{defn}

\begin{figure}[htbp]
\begin{center}
\includegraphics*[height=1.2in]{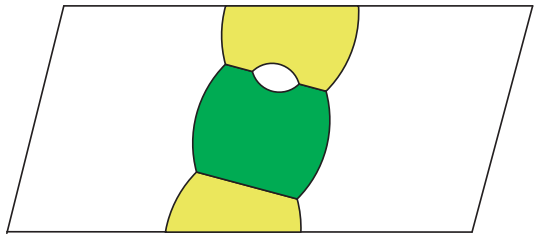}%
\end{center}
\caption{\label{fig:funky_chain} The two shaded six-sided
components form a \textit{chain}, even though there is a curvilinear digon
embedded between them. The two shaded components plus the
curvilinear digon do not form a chain, because the swath they form
is not minimal.}
\end{figure}

\begin{defn}
A \emph{chain}, such as the standard chain of Figure
\ref{fig:introduce_minimizers} or the chain of Figure
\ref{fig:funky_chain}, is a minimal swath --- that is, a swath
without a proper subset of components that themselves form a
swath. Note that every swath contains at least one chain and that
every chain is itself a swath.

\end{defn}

%%%%%%%%%%%%%%%%%%%%%%%%%%%%%%%%%%%%%%%%%%%%%%%%%%%%%%%%
%%%%%%%%%%%%%%%%%%%%%%%%%%%%%%%%%%%%%%%%%%%%%%%%%%%%%%%%

\begin{defn} Two closed curves lie in the same \emph{homology
class}, or have the same homology, if they can be smoothly
deformed into each other along the surface of the torus. Given two
different oriented directions of the torus (i.e., directions of
closed geodesics on the torus), the homology class of a closed
curve is given by an ordered pair $(p, q)$, where the curve wraps
around the torus $p$ times in the first direction and $q$ times in
the other direction. (The ordered pair is well-defined up to
negating both $p$ and $q$.) We will refer to the homology of a
band or chain $K$, by which we mean the homology of any
non-contractible component of the boundary of $\overline{K}$.
\end{defn}

\subsection{Geometric Definitions}

\begin{defn} \label{defn:definition sdb}
The \emph{\sdb}\  is comprised of three circular arcs
meeting in threes at $\onetwenty$, such that the curvature of the
arc separating the contractible components is the difference of the curvatures of the outer
caps. For prescribed areas, there is a standard double bubble in the Euclidean plane, unique up to congruence [M1, Proposition 14.1]; it may or may not fit on the torus.
\end{defn}

\begin{defn}
The \emph{double band} consists of two adjacent bands, bounded by three shortest closed geodesics.
\end{defn}

\begin{defn}
A \emph{symmetric chain} is a chain that is symmetric about a
closed geodesic with the same homology as the chain.
\end{defn}

\begin{defn} \label{defn:definition standard chain}
A \emph{standard chain} is a symmetric chain enclosing two
four-sided components, whose boundary is comprised of six circular
arcs that satisfy the regularity conditions in Proposition \ref{prop:regularity}: that is, the arcs meet in threes at four common vertices at $\onetwenty$
such that the curvature of each arc separating the two components
in the chain is the difference of the curvatures of the outer arcs
(see Figure \ref{fig:introduce_minimizers}). The
\textit{axis-length} of a standard chain is the length of the
geodesic with the same homology as the chain.
 \end{defn}

\begin{figure}[htbp]
\begin{center}
\includegraphics*[height=3.25in]{\Path/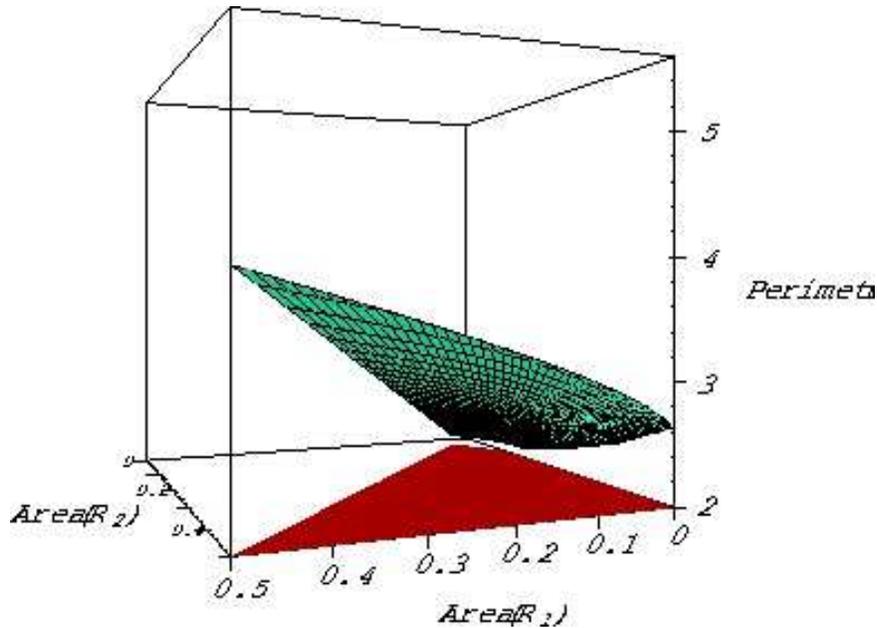}%
\end{center}
\caption{\label{fig:showsnofold} This parametric plot from
\emph{Mathematica}  of perimeter as a function of areas for standard chains with axis-length 1. The plot has no folds in it and contains only one sheet,
suggesting uniqueness for prescribed area pairs $A_1 \le A_2$.}
\end{figure}

\begin{remsub} \label{rem:symm_chain is unique}
Our parameterization of the standard chain (Section
\ref{SS:Param_2_Comp_Sym_Chains}, see Remark \ref{rem:standard
chain exists}) proves existence of some family of standard chains.
Furthermore, there exists graphical evidence that the only chains
that are minimizing are those with the same homology as a shortest closed geodesic (Remark
\ref{rem:standard chain short direction}).

We believe the standard chain is also unique up to isometries of
the torus for given areas and homology class. Indeed, a parametric
plot in \emph{Mathematica} (see Figure \ref{fig:showsnofold})
suggests that different standard chains with the same
homology enclose different area pairs.  \end{remsub}

%%%%%%%%%%%%%%%%%%%%%%%%%%%%%%%%%%%%%%%%%%%%%%%%%%%%%%%%
%%%%%%%%%%%%%%%%%%%%%%%%%%%%%%%%%%%%%%%%%%%%%%%%%%%%%%%%

\begin{defn} \label{defn:lens}
A \emph{lens} is a pair of congruent circular arcs meeting at two
vertices at $\onetwenty$.
\end{defn}

\begin{remsub} \label{rem:lens-fits}
It is easily shown that if the distance between the two vertices
is less than or equal to 1, then the lens can be embedded on any
torus with its axis along any short direction.
\end{remsub}

\begin{defn} \label{defn:definition band lens}
A \emph{band lens} consists of a band and a lens, such that the
boundary curves meet in threes at two vertices at $\onetwenty$, as
in Figure \ref{fig:introduce_minimizers}. \end{defn}

\begin{remsub}  By Lemma \ref{lem:bl_short_direction}, a
minimizing band lens on the torus has the same homology as a
shortest closed geodesic (that is, the closed geodesic bordering
the band has length one). It is easy to see that when a minimizing band lens enclosing two given areas exists, it is essentially unique for those areas (up
to isometries of the torus).
\end{remsub}

\begin{defn}
A \emph{hexagon tiling} is a tiling in which every
component of the double bubble and exterior is a hexagon with interior angles of $\onetwenty$.
\end{defn}

\begin{defn}
The \emph{standard hexagon tiling} is a hexagon tiling on the
hexagonal torus, which partitions the torus into three hexagonal
regions, whose boundary edges lie parallel to the three short
directions of the torus. By Corollary
\ref{cor:uniqueness_of_standard_hex_tiling}
and Lemma
\ref{lem:existence-of-some-hex-tilings}, the standard hexagon
tiling exists and is unique for a particular range of area pairs.
\end{defn}

\section{Existence, Regularity, and Basic Properties}\label{S:basic_propositions}
This section provides existence and regularity of
perimeter-minimizing double bubbles (Proposition
\ref{prop:regularity}), a bound on the number of components for
minimizers (Proposition \ref{prop:component_bound}), and an easy
but important perimeter bound (Proposition
\ref{prop:perimeter_bound}).

\begin{prop}\label{prop:regularity}
 {\rm Existence and Regularity Theorem [M2, 2.3
and 2.4].} In a smooth Riemannian surface $S$ with compact quotient
by its isometry group, for any two areas $A_1$
and $A_2$ (whose sum is less than the area of $S$), there
exists a least-perimeter enclosure of the two areas.  This
enclosure consists of finitely many smooth constant-curvature curves meeting in
threes at $\onetwenty$ angles. No boundary curve separates components of the same region, and all curves separating a specific
pair of the three regions have the same (signed) curvature.  Moreover, the enclosure satisfies
the cocycle condition: the sum of the
signed curvatures around any closed path is zero.
\end{prop}

\begin{rem} \label{rem:finite} In particular, given any minimizing
double bubble on the torus, each of the three regions contains finitely many
components. \end{rem}

\begin{rem} \label{rem:connected}
It follows that in a minimizing double bubble on the torus, the
union of the closures of any two regions is connected. Otherwise,
two disjoint pieces of this union can be translated to be made tangent, forming
a new minimizer that violates regularity, a contradiction. In
particular, each component (in either of the two enclosed regions
or the exterior region) is adjacent to components of both other
regions.
\end{rem}

\begin{rem} \label{rem:alternating_regions}
It also follows that alternating sides of a component must border different regions. Hence, any contractible component in a minimizing double bubble or its exterior has an even number of sides.
\end{rem}

\begin{rem} \label{rem:equilibrium and variation}
A perimeter-minimizing double bubble is in \textit{equilibrium}: under any
smooth variation that preserves areas, the derivative of perimeter with
respect to time is initially zero.  (Here by ``variation,'' we mean a
perturbation of the double bubble over time.) For a piecewise smooth double bubble, being in equilibrium is equivalent to having boundary curves that satisfy the given curvature conditions, satisfy the cocycle condition, and meet in threes at $\onetwenty$ (see, for instance, [HMRR], Lemma 3.1).
\end{rem}

\begin{rem} \label{rem:pressure}
The curvature and cocycle conditions imply the existence of a well-defined
\emph{pressure} for each region with the following property: the difference between the
pressures of any two regions in an equilibrium double bubble is
the curvature of the separating arcs between those regions, signed
so that a curve separating two regions bows into the region of
lower pressure. (For instance, we could set the pressure of the exterior region to equal zero and set the pressure of each enclosed region to equal the curvature of the interface separating that region from the exterior.) A least-pressure nonpolygonal contractible
component must have more than six sides; and a highest-pressure
nonpolygonal contractible component must have fewer than six
sides.
\end{rem}

\begin{prop}\label{prop:component_bound} {\rm (Wichiramala Component Bound) [MW, Proposition
3.1]}  Unless all three regions have equal pressure, each region
of highest pressure in a perimeter-minimizing double bubble in the
flat two-torus has at most two contractible components.
\end{prop}

\begin{proof}
The proof in [MW], which we sketch here, uses variations that shrink or expand a high-pressure nonpolygonal component at a constant rate along its boundary. If there are at least three such components in some region, then some non-trivial linear combination of the corresponding ``shrink and expand'' variations will preserve the areas of each region to first order. Each of the three variations will contribute a negative term to the second variation of the combined variation, so that the new variation has negative second variation (i.e., negative second derivative of perimeter with respect to time). In this case, we call the double bubble ``unstable,'' and no such double bubble is perimeter minimizing. The topology of the torus does not affect this proof.

\end{proof}

\begin{prop}\label{prop:perimeter_bound}
 {\rm (Perimeter Bound)} The perimeter of a perimeter-minimizing double bubble in the flat two-torus is at most 3.
\end{prop}

\begin{proof}
A double band of perimeter three can enclose any pair of prescribed areas, so no minimizer
can have greater perimeter.
\end{proof}

\section{Topological Classification of Minimizers}\label{S:topological classification}

This section proves two important propositions. Proposition
\ref{prop:three_equal_pressures} says that if the three regions in
a perimeter-minimizing double bubble all have equal pressure, then
it is a hexagon tiling or the double band. Proposition \ref{prop:topo
classification} is the backbone for the proof of our Main Theorem
\ref{thm:main theorem}. It shows that any minimizing double bubble
must be a tiling, a contractible double bubble, a swath, a single
band adjacent to a contractible set of components, or the double
band.  Section \ref{S:beyond} will treat in turn each of these
topological classes.

\begin{lem}\label{lem:non-contractible is band}
A non-contractible component $C$ in a minimizing double bubble
must either be a band or have contractible
complement.
\end{lem}

\begin{proof}
If every component of the boundary of $C$ is
contractible, then the complement of $C$ is contractible.

Otherwise, some component $\sigma$ of the boundary of $C$ is non-contractible, and any
other non-contractible component of the boundary of $C$ has the same
homology as $\sigma$. Thus, $C$ must be a band with that homology, except
possibly with some collection of contractible components in its interior.
But such contractible components never exist in a minimizing double bubble
(Remark \ref{rem:connected}), so $C$ must indeed be a band.
\end{proof}

\begin{lem} \label{lem:wrap only once}
Consider a fundamental domain of the torus, a parallelogram with
sides $1, L$ and interior angle $\theta \in [\pi/3, \pi/2]$.
Suppose that a chain or band is part of a minimizing double bubble
and has homology $(p, q)$, defined with respect to two directions
along the sides of the parallelogram. Then $|p|, |q| \le 1$.
\end{lem}

\begin{proof}
Each band or chain has least two disjoint, non-contractible boundary
curves. If $|p| \ge 2$ or $|q| \ge 2$, then the length of each of
these boundary curves is at least $\sqrt{3}$, and the total
perimeter is at least $2\sqrt{3}$. By the perimeter bound (Proposition
\ref{prop:perimeter_bound}), the double bubble cannot be a
minimizer.
\end{proof}

\begin{lem} \label{lem:less than two bands or chains}
A perimeter-minimizing double bubble on the torus, together with its exterior, has at
most three disjoint bands or chains. If there are three,
then the double bubble must be the double band.
\end{lem}

\begin{proof}
Each band or chain has at least two
non-contractible boundary curves, and each such boundary curve
borders at most two bands or chains. Hence, if there are $n$ bands and
chains, then there will be at least $n$ such non-contractible
boundary curves, each with length at least one. Then by
Proposition \ref{prop:perimeter_bound}, $n \leq 3$. If $n = 3$,
then the boundary must consist entirely of three parallel closed
geodesics, each of length one. In this case, the double bubble
must be the double band.

\end{proof}

\begin{lem} \label{lem:equal-pressures-components}
In a perimeter-minimizing double bubble with three equal
pressures, every component in the double bubble plus exterior is
either a vertex-free band or a hexagon.
\end{lem}

\begin{proof} Consider any component in the double bubble plus exterior, and let $\mathcal{C}$ be the boundary of this component.

Because the double bubble has three equal pressures, any vertex-free component of $\mathcal{C}$ is linear and thus must be a closed geodesic.

Suppose that a component of $\mathcal{C}$ has at least one vertex. Then it has positive total turning angle, whereas any non-contractible closed curve has total turning angle 0. Hence, if a component of $\mathcal{C}$ has at least one vertex, then it is contractible. Also, it consists of straight edges meeting at interior angles measuring $2\pi/3$ --- precisely the configuration of a hexagon.

Therefore, every component of $\mathcal{C}$ is a closed geodesic or a contractible hexagon, and the desired result follows easily. \end{proof}

\begin{prop} \label{prop:three_equal_pressures}
A perimeter-minimizing double bubble with three equal pressures is
either a hexagon tiling or the double band.
\end{prop}

\begin{proof} By Lemma \ref{lem:equal-pressures-components}, each
component in the double bubble plus exterior is a vertex-free band or a
hexagon. If every component is a hexagon, then the double bubble
is a hexagon tiling. If every component is a band, then because
there are at least three components, by Lemma \ref{lem:less than
two bands or chains} the double bubble is a double band. Finally,
the double bubble cannot contain both bands and hexagons, because
no band can be adjacent to a hexagon without violating regularity
(Proposition \ref{prop:regularity}). \end{proof}

\begin{prop}\label{prop:topo classification}
A minimizing double bubble on the torus must be a contractible
double bubble, a tiling, a swath, a band adjacent to a set of
components such that the closure of their union is contractible,
or the double band (after perhaps relabelling the two regions and
the exterior).
\end{prop}

\begin{proof}
We will categorize the minimizing double bubble by the number and
type of non-contractible components in the double bubble plus its
exterior.  By Lemma \ref{lem:non-contractible is band} any
non-contractible component must either be a band or have
contractible complement. If the double bubble plus its exterior
has in any region (say, in the exterior) a component with
contractible complement, then the double bubble is contractible.
Thus, for our remaining cases we can assume that all
non-contractible components are bands.

If the double bubble plus its exterior has no bands, then every
component (of the double bubble plus its exterior) is contractible, and by definition the double bubble is a tiling.

If the double bubble plus its exterior has exactly one band, say in the exterior, then the double bubble consists of contractible components.  But the closure of the union of this set of components is non-contractible, because it shares a non-contractible boundary curve with the band. Furthermore, any minimizing double bubble is connected (Remark \ref{rem:connected}). Then by definition the double bubble is a swath.

If the double bubble plus its exterior has exactly two bands (say, with at least one
band in the exterior), then consider the set $C$ of the contractible components --- that is,
the set of components different from the bands. Suppose some
subset of these components has non-contractible and connected
closure.  Then $C$ contains a swath, and therefore a chain.  Now
the double bubble plus its exterior has two bands and at least one
chain, so by Lemma \ref{lem:less than two bands or chains} it is
not perimeter minimizing.  Therefore, $C$ must be a contractible
set of components.  The two bands must then be adjacent to each
other, and hence lie in different regions. Then the double bubble contains exactly one band, together with a subset of $C$ --- that is, the double bubble is a band adjacent to contractible
set of components.

If the double bubble plus its exterior has three or more bands,
then by Lemma \ref{lem:less than two bands or chains}, the double
bubble must be the double band.

\end{proof}

\section{The Four Classes with Contractible Components} \label{S:beyond}
By Proposition \ref{prop:topo classification}, upon relabelling,
every minimizing double bubble on the torus must be a contractible
double bubble, a tiling, a swath, a band adjacent to a set of
components whose union is contractible, or the double band.  While
the double band describes a specific geometric configuration, each
of the remaining possibilities is merely a topological class of
configurations.  This section identifies the potential minimizer
in each of these four remaining topological classes: the standard
double bubble among contractible double bubbles (Proposition
\ref{prop:sdb_is_minimizer}), a standard chain among swaths
(Proposition \ref{prop:chain_is_minimizer}), the band lens among
double bubbles with a single band adjacent to a contractible set
of components (Proposition \ref{prop:band lens_is_minimizer}), and
the standard hexagon tiling among tilings (Proposition
\ref{prop:only tiling is 3 hex}).

\subsection{Contractible Double Bubbles}\label{SS:contract_bcs}
Proposition \ref{prop:sdb_is_minimizer} uses the result from
$\textbf{R}^2$ that perimeter-minimizing double bubbles are
standard [F], to show that the only potential minimizer among
contractible double bubbles is the standard double bubble.

\begin{propsub} \label{prop:perimeter is 3 times diameter}
The perimeter of any planar standard double bubble is greater than $\pi$
times its diameter: $P > \pi D$.
\end{propsub}

\begin{proof}
We obtain a formula for $P/D$ as a function $f(\theta)$, where
$\theta \in [0, \pi/3)$ is the angle between the interior arc of
the standard double bubble and the chord that connects the
endpoints of the arc (see Figure \ref{fig:sdb_param}). Then we
show $f(\theta) > \pi$ for the entire domain $0 < \theta < \pi/3$.
(A separate computation easily verifies the inequality $P > \pi D$
when $\theta = 0$.)

%%%%%%%%%%%%%%%%%%%%%%%%%%%%%%%%%%%%%%%%%%%%%%%%
\begin{figure}[ht]
\begin{center}
\includegraphics*[height=1.75in]{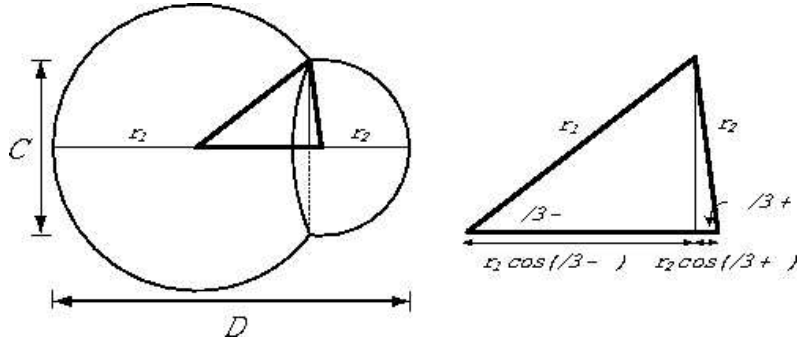}%

\end{center}
\caption{\label{fig:sdb_param} The perimeter of the standard
double bubble is greater than $\pi$ times its diameter.}
\end{figure}
%%%%%%%%%%%%%%%%%%%%%%%%%%%%%%%%%%%%%%%%%%%%%%%%

Let $r_1, r_2, c_1, c_2$ be the radii and centers of the longer
and shorter exterior arcs of the standard double bubble,
respectively.  Then the diameter is the sum of $r_1$, $r_2$, and
the distance between $c_1$ and $c_2$:

$$D = r_1 + r_2 + r_1\cos{(\pi/3-\theta)} + r_2\cos{(\pi/3+\theta)}.$$
Applying the formula for radii (Equation \ref{eqn:R}), we obtain
an expression that can be solved for $C$, the length of the chord
connecting the two vertices:
$$
D = C \left(\frac{1+\cos{(\pi/3-\theta)}} {2\sin{(2\pi/3+\theta)}} +
\frac{1+\cos{(\pi/3+\theta)}} {2\sin{(2\pi/3-\theta)}}\right).
$$
Substituting for $C$ in the formula for perimeter (Equation
\ref{eqn:L}) yields after some simplification the desired function
in $\theta$:
$$
f(\theta) = P/D = \frac{8\pi\sin{\theta}\cos{\theta} +
3\sqrt{3}\theta}{6\sin{\theta}\cos{\theta} + 3\sin{\theta}}.
$$
Now, showing $f(\theta) > \pi$ is equivalent to showing
$$
8\pi\sin\theta\cos\theta + 3\sqrt{3}\theta > \pi(6\sin\theta\cos\theta
+ 3\sin\theta).
$$
Moving the terms to the left hand side and using
the identity $\sin(2\theta) = 2\sin\theta\cos\theta$ gives the
equivalent inequality
$$g(\theta) = \pi\sin(2\theta) + 3\sqrt{3}\theta - 3\pi\sin\theta > 0,$$
with
$$ g''(\theta) = -4\pi\sin(2\theta) + 3\pi\sin\theta =
\pi\sin\theta(-8\cos\theta + 3).$$
For $\theta$ in $(0, \pi/3)$, $\sin\theta > 0$ while $\cos\theta >
\cos(\pi/3) > 3/8$. Therefore, $g''(\theta)$ is strictly negative.
Because $g(\theta)$ equals 0 at both endpoints of this interval, it must
be positive along the interior, as desired.

\end{proof}

\begin{propsub}[Suggested by Masters \cite{[Ma1]}.] \label{prop:sdb_is_minimizer} The standard double bubble is the only contractible double bubble that may be a minimizer.
\end{propsub}

\begin{proof}

Assume that there is a perimeter-minimizing contractible double
bubble $\Sigma$ enclosing areas $A_1$ and $A_2$ that is not the
standard double bubble.  By [F], the unique perimeter-minimizing
solution in \Rtw\ for the same prescribed areas is a standard
double bubble, $\Theta$. By the perimeter bound (Proposition
\ref{prop:perimeter_bound}), $\Sigma$ has perimeter at most three.
Since $\Theta$ has less perimeter than $\Sigma$, the perimeter $P$
of $\Theta$ satisfies $P < 3$.

%%%%%%%%%%%%%%%%%%%%%%%%%%%%%%%%%%%%%%%%%%%%%%%%
\begin{figure}[ht]
\begin{center}
\includegraphics*{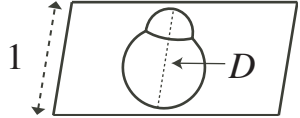}%

\end{center}
\caption{\label{fig:when sdb will not fit} A standard double
bubble with minimizing perimeter has diameter $D \le 1$ and hence
fits on the torus.}
\end{figure}
%%%%%%%%%%%%%%%%%%%%%%%%%%%%%%%%%%%%%%%%%%%%%%%%

By Proposition \ref{prop:perimeter is 3 times diameter}, the
diameter $D$ of $\Theta$ satisfies $D < P/\pi < 1$. Hence, it fits
on the torus with its axis of symmetry in any direction (Figure
\ref{fig:when sdb will not fit}). Therefore, $\Theta$ is a double
bubble on the torus enclosing areas $A_1$ and $A_2$ with less
perimeter than $\Sigma$, a contradiction. \end{proof}

\begin{figure}[ht]
\begin{center}
\includegraphics*[height=1.0in]{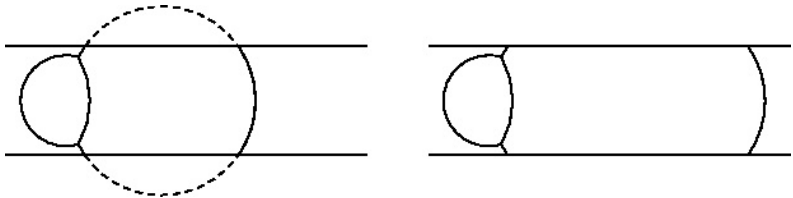}%
\end{center}
\caption{\label{fig:standard doesnt fit}If a standard
double bubble does not fit on the cylinder, we can easily
construct a double bubble with even less perimeter.}
\end{figure}

\begin{remsub}
There is an easier argument on the infinite cylinder: when
$\Theta$ does not fit with its diameter parallel to the sides of
the cylinder, we can modify its boundary to create a new double
bubble with less perimeter than $\Theta$ and enclosing the same
areas, contradicting the assumption that there is a non-standard
contractible minimizer (see Figure \ref{fig:standard doesnt fit}).
\end{remsub}

%%%%%%%%%%%%%%%%%%%%%%%%%%%%%%%%%%%%%%%%%%%%%%%%
%%%%%%%%%%%%%%%%%%%%%%%%%%%%%%%%%%%%%%%%%%%%%%%%

\subsection{Swaths of Contractible Components}\label{SS:swaths}
Recall that by Definition \ref{defn:definition swath} a swath is a
set of contractible components with non-contractible and connected
closure, and with non-contractible complement, such as the
standard chain of Figure \ref{fig:introduce_minimizers}. This
section shows that minimizing chains are standard chains
(Definition \ref{defn:definition standard chain}), exactly as in
Figure \ref{fig:introduce_minimizers}. Lemma \ref{lem:chains}
shows that minimizing swaths must be chains consisting of two or
four curvilinear quadrilateral components. Lemma
\ref{lem:asymmetric_chains} shows that asymmetric chains
satisfying regularity cannot be minimizers since they have
perimeter greater than three. Finally, Proposition
\ref{prop:chain_is_minimizer} shows that minimizing symmetric
chains are standard.

\begin{remsub} In this section, we often speak about contractible four-sided components as if they sat in the plane rather than a torus, because it is easier to talk about lines, reflections, and circles in the plane. (Lines can wrap around a torus infinitely many times, reflections across lines on a torus do not always exist, and circles on a torus may overlap themselves.) \end{remsub}

\begin{lemsub} \label{lem:ww-cyc-quad} {\rm [W, Lemma 5.31].} Given a contractible four-sided component in a regularity-satisfying double bubble or its exterior, the vertices of the component lie on a circle or a line. They lie clockwise on the circle or line in the same or opposite order as they appear clockwise on the boundary of the four-sided component. \end{lemsub}

\begin{proof} {\rm [W, Lemma 5.31]} actually considers four-sided components
in regularity-satisfying triple bubbles in the plane. However, the proof
looks at each contractible four-sided component in isolation from the rest
of the triple bubble, using only the facts that each side has constant
curvature, that the sides do not intersect each other, and that the sides
meet at angles of $2\pi/3$. Thus, the proof also applies to contractible
four-sided components of double bubbles in the torus.

[W, Lemma 5.31] is also slightly more limited than our result, not proving
that if the vertices of the component lie on a line, then they do so in
the same order that they appear on the boundary of the component.
However, it does offer a proof of the analogous fact about the vertices
when they lie on a circle, and this proof can be easily modified to give
our complete result.
\end{proof}

\begin{lemsub}\label{lem:isos-or-same-circle} Given a contractible four-sided component with vertices $A,B,C,D$ in a regularity-satisfying double bubble or its exterior, $BC = DA$ unless sides $BC$ and $DA$ of the component lie on the same circle or line. \end{lemsub}

\begin{proof} By Lemma \ref{lem:ww-cyc-quad}, vertices $A,B,C,D$ lie on a circle or line $\omega$. Suppose that sides $AB$ and $CD$ do not lie on the same circle or line. By regularity in double bubbles, sides $AB$ and $CD$ have constant curvature and lie on some circles or lines $\omega_1$ and $\omega_2$, respectively. Also by regularity, sides $AB$ and $CD$ have the same curvature, and they meet side $BC$ at equal angles (namely, $2\pi/3$). Hence, $\omega_1$ and $\omega_2$ are reflections of each other across the perpendicular bisector of $\overline{BC}$. Thus, $\omega_1 \cap \omega = \{A,B\}$ and $\omega_2 \cap \omega = \{C, D\}$ are reflections of each other across the perpendicular bisector of $\overline{BC}$ as well, implying that $AB = CD$. \end{proof}

\begin{lemsub}\label{lem:cyclic_quadrilateral} Given a contractible four-sided component in a regularity-satisfying double bubble or its exterior, the vertices of the component do not lie on a line. \end{lemsub}

\begin{proof} Let the vertices of the component be $A, B, C, D$ in that order along the boundary of the component. Suppose, for sake of contradiction, that the vertices lie on a line $\ell$. By Lemma \ref{lem:ww-cyc-quad}, $A, B, C, D$ lie in that order along $\ell$ --- without loss of generality, $\ell$ is horizontal with $A$ at the left and $D$ at the right. Because $BC \ne DA$, Lemma \ref{lem:isos-or-same-circle} implies that sides $DA$ and $BC$ of the component lie on the same circle or line. The only circle or line that contains all four vertices $A, B, C, D$ is $\ell$, so that sides $DA$ and $BC$ of the component are linear and lie along $\ell$. However, this is impossible because then sides $DA$ and $BC$ would overlap, a contradiction. \end{proof}

%%%%%%%%%%%%%%%%%%%%%%%%%%%%%%%%%%%%%%%%%%%%%%%%%%%%%%%%%%%%
\begin{figure}[ht]
\begin{center}
\begin{tabular} {@{\extracolsep{1.5 cm}} cc}
\includegraphics*[height=1in]{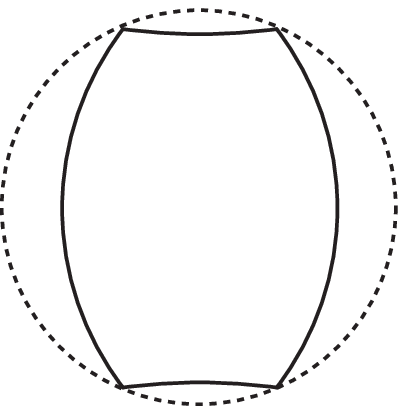}%
&
\includegraphics*[height=1in]{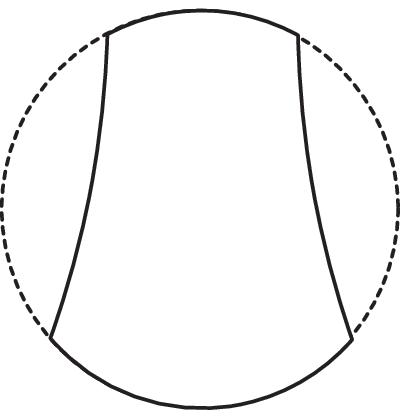}%
\\
(a) & (b)
\end{tabular}
\end{center}
\caption{\label{fig:possible-quads}
The vertices of any four-sided component lie on a circle. Either (a) the vertices form a rectangle, or (b) two opposite sides of the component lie along the same circle.}
\end{figure}
%%%%%%%%%%%%%%%%%%%%%%%%%%%%%%%%%%%%%%%%%%%%%%%%%%%%%%%%%%%%

\begin{lemsub}\label{lem:isosceles trapezoid} {\rm(Compare [W,
proof of Lemma 5.3])} Given a contractible four-sided component in
a regularity-satisfying double bubble or its exterior, the vertices of the component form a non-degenerate isosceles trapezoid. The vertices appear clockwise on this trapezoid in the same order as they
appear clockwise on the boundary of the component. If the
trapezoid is not a rectangle, then the two boundary curves
subtended by the parallel sides of the trapezoid lie on a circle
circumscribing the polygonal quadrilateral.
\end{lemsub}

\begin{proof} By Lemmas \ref{lem:ww-cyc-quad} and \ref{lem:cyclic_quadrilateral}, the polygonal quadrilateral $ABCD$ formed by the vertices of the component is inscribed
in some circle $\omega$; the vertices appear on the quadrilateral
in the same order that they appear on the circle (although perhaps
with different orientation).

By regularity, curvilinear sides $AB, BC, CD, DA$ cannot all lie along $\omega$. Without loss of generality, assume that sides $AB$ and $CD$ do not lie on the same circle. By Lemma \ref{lem:isos-or-same-circle}, we have $AB = CD$. Hence, because quadrilateral $ABCD$ is inscribed in a circle, it must be an isosceles trapezoid with $\overline{BC} \parallel \overline{DA}$.

If sides $BC$ and $DA$ of the curvilinear component do not lie on the same circle, then by the same analysis we have $BC = DA$, so that quadrilateral $ABCD$ is a rectangle (Figure \ref{fig:possible-quads}(a)). And if trapezoid $ABCD$ is not a rectangle, say because $BC \ne DA$, then sides $BC$ and $DA$ of the curvilinear component \textit{do} lie on the same circle (Figure \ref{fig:possible-quads}(b)).

Finally, it remains to show that the vertices appear clockwise on
the boundary of the component in the same order as they appear
clockwise on the trapezoid. Because the vertices form a trapezoid (say, with $\overline{BC} \parallel \overline{DA}$), the four-sided component is symmetric about the perpendicular bisector $m$ of sides $BC$, $DA$ of the component. It follows easily that neither side $AB$ nor side $CD$ intersects $m$, implying the desired fact about the order of the vertices.
\end{proof}

\begin{lemsub}\label{lem:High-Pressure Quadrilateral's Low-Pressure Arcs Are
Parallel} In a contractible four-sided component of a
highest-pressure region in a regularity satisfying double bubble or its exterior, the arcs bordering a
least-pressure region are subtended by parallel chords, and the
arcs bordering a medium-pressure region are subtended by congruent
chords.
\end{lemsub}

\begin{proof} By Lemma \ref{lem:isosceles trapezoid}, the chords subtending the
boundary curves of the component form an isosceles trapezoid
inscribed in a circle $\omega$. If the trapezoid is a rectangle,
we are done. Otherwise, the two boundary curves $\sigma, \sigma'$
subtended by parallel chords lie on $\omega$, and the other
two boundary curves $\mu, \mu'$ are subtended by congruent chords. Because $\sigma, \sigma'$ meet the remainder of $\omega$ at angles \textit{greater}
than $2\pi/3$ (namely, $\pi$), the other boundary
curves $\mu, \mu'$ of the component bow within $\omega$. In other
words, the curvature of $\mu$ is less than the curvature of
$\omega$, which equals the curvature of $\sigma$. It follows that
$\sigma, \sigma'$ bound a region of less pressure than $\mu, \mu'$
do. \end{proof}

\begin{lemsub}\label{lem:six_edges}
For a minimizing tiling on a torus, the average number of edges
per component is exactly six.
\end{lemsub}

\begin{proof}
The Euler characteristic on the torus is $v-e+c=0$, where $v$ is
the number of vertices, $e$ is the number of edges, and $c$ is the
number of components.  By regularity,
(Proposition \ref{prop:regularity} and Remark \ref{rem:finite}), $v$, $e$ and $c$ are all finite.
Now, let $a$ be the average number of edges per component.  Then
$c = \frac{2}{a} e$, since each edge is adjacent to two
components. Similarly, $v = \frac{2}{3} e$, since each edge is
adjacent to two vertices and (by regularity, Proposition
\ref{prop:regularity}) each vertex is adjacent to three edges.
Substituting into the Euler formula, we obtain $(\frac{2}{3} - 1 +
\frac{2}{a} ) e = 0$, from which we obtain $a = 6$.
\end{proof}

\begin{propsub}\label{prop:Lowest_Pressure_Region_Non-Contractible} A minimizing double bubble
with a contractible component in a region of least pressure must
be either an octagon-square tiling or a hexagon tiling.
\end{propsub}

\begin{proof}
If the pressures of the three regions are equal, then by
Proposition \ref{prop:three_equal_pressures} either the double
bubble is a hexagon tiling or the double band; the latter case is
impossible because the double bubble must have a contractible
component by assumption. Hence, we may assume that the three
regions do not all have the same pressure.

If any component in the double bubble plus exterior is a vertex-free band, then its boundary curves are closed geodesics. This
either implies that both regions and the exterior have equal
pressure, contradicting our previous assumption; or we can slide the
band until it collides with some other boundary curve, for a
violation of regularity (Remark \ref{rem:connected}). Thus, no component is a vertex-free band.

If any region has contractible complement,
then by Proposition \ref{prop:sdb_is_minimizer}, the double bubble
is the standard double bubble. Again, the least-pressure region
consists solely of non-contractible components, a contradiction.
Hence, we may assume that every region has non-contractible complement.

Let $R_1$ be a region of highest pressure. From the previous
paragraph, any non-contractible component of $R_1$ is a band.
Because $R_1$ is of highest pressure, such a band cannot have any
vertices; if it did, then each boundary edge of the band would
have positive turning angle, a contradiction. However, we have
already shown that no region contains a vertex-free band.
Therefore, any component in $R_1$ is contractible. By Proposition
\ref{prop:component_bound}, $R_1$ has at most two contractible
components. Also, each component of $R_1$ contains at most four
vertices (Remark \ref{rem:pressure}). Hence, $R_1$ contains at
most eight vertices.

Let $R_2$ be a region of lowest pressure, containing a
contractible component with $n$ sides (and perhaps other
contractible, or non-contractible, components). This component has
least pressure but is not polygonal (because the three pressures
are not all equal), so $n > 6$ (Remark \ref{rem:pressure}). By
Remark \ref{rem:alternating_regions}, $n$ is even, implying that
$n \ge 8$.

Because every vertex of $R_1$ is a vertex of $R_2$ (and vice
versa), and because neither $R_1$ nor $R_2$ contains a vertex-free
band, it follows that $R_1$ consists of two contractible
four-sided components and that $R_2$ consists of a single
contractible curvilinear octagon.

Because no region has contractible complement, $R_2$ and some
four-sided component $C$ in $R_1$ must form a chain that wraps around the torus. We now consider two cases.

%%%%%%%%%%%%%%%%%%%%%%%%%%%%%%%%%%%%%%%%%%%%%%%%%%%%%%%%%%%%
\begin{figure}
\begin{center}
\includegraphics*{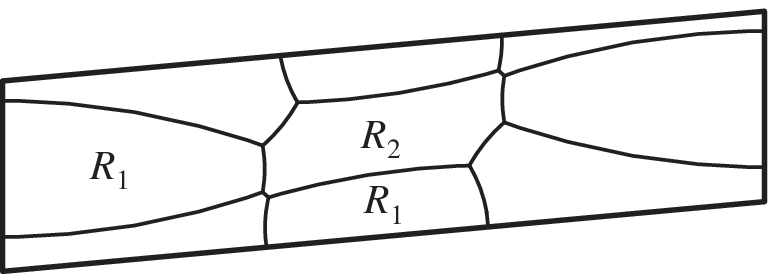}%
\\ (a) \\ \bigskip
\includegraphics*{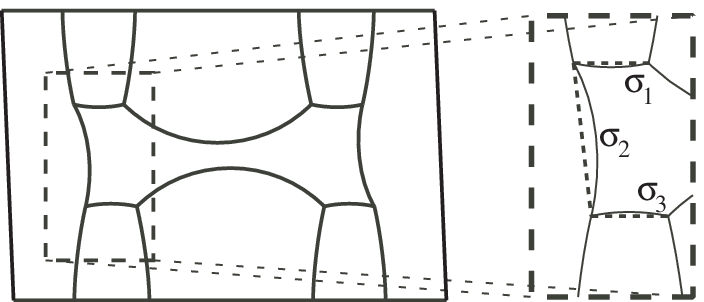}%
\\ (b) \\
\end{center}
\caption{\label{fig:low-pressure_cases} Two double
bubbles with a purported contractible low-pressure component.}
\end{figure}
%%%%%%%%%%%%%%%%%%%%%%%%%%%%%%%%%%%%%%%%%%%%%%%%%%%%%%%%%%%%

\textit{Case 1}: $C$ borders two opposite sides of $R_2$. Then just as $R_2$ and $C$ form a chain wrapping around the torus, so too do $R_2$ and the \textit{other} four-sided component of $R_1$ --- wrapping around the torus in a different direction. In this case, the double bubble is a
tiling. Because $R_1$ and $R_2$ each have eight vertices, $R_3$
must have eight vertices as well, for a total count of 24
vertices. But by Lemma \ref{lem:six_edges}, the average number of
vertices per component must be exactly six. Therefore, there are
exactly four components, implying that $R_3$ consists of one
curvilinear octagon. Therefore, as shown in Figure
\ref{fig:low-pressure_cases}(a), the configuration is an
octagon-square tiling.

\textit{Case 2}: $C$ borders two sides $\sigma_1, \sigma_3$ of
$R_2$ as depicted in Figure \ref{fig:low-pressure_cases}(b). Let
$\sigma_2$ be the side of the eight-sided component between them.
By Lemma \ref{lem:High-Pressure Quadrilateral's Low-Pressure Arcs
Are Parallel}, the chords subtending $\sigma_1$ and $\sigma_3$ are
parallel. Therefore, one of these chords must meet the chord
subtending $\sigma_2$ at a non-obtuse angle. But then the
corresponding arcs meet at an angle less than $2\pi/3$, a
contradiction.

Therefore, if a least-pressure region contains a contractible component, then the double bubble is a tiling with only straight edges or an octagon-square tiling. \end{proof}

\begin{corsub} \label{cor:hex-or-oct} A perimeter-minimizing
tiling is either a hexagon tiling or an octagon-square tiling.
\end{corsub}

\begin{lemsub} \label{lem:circles_propagate} Let $\gamma$ be a
boundary edge separating two components $X_1$ and $X_2$ in a minimizing double bubble,
where both $X_1$ and $X_2$ have different pressures from the third region of the double bubble. Let the two boundary arcs of $X_1$ adjacent to $\gamma$ be $\sigma_1$ and
$\sigma_2$, and let the two boundary arcs of $X_2$ adjacent to
$\gamma$ be $\mu_1$ and $\mu_2$. (The $\sigma_i$ coincide if $X_1$
has two sides, and the $\mu_i$ coincide if $X_2$ has two sides.)
Suppose that $\sigma_1$ and $\sigma_2$ have the following
property: they lie on a common arc, and this arc forms an immersed
curvilinear digon with $\gamma$. Then $\mu_1$ and $\mu_2$ have
this property as well.
\end{lemsub}

\begin{proof} Consider the curvilinear digon bordered by $\gamma$ and the arc containing
$\sigma_1$ and $\sigma_2$. We may add a third arc $\mu$ to the
curvilinear digon to form an immersed standard double bubble.
Because the cocycle condition holds for both the given double bubble and for the immersed standard double bubble,
$\mu_1$, $\mu_2$, and $\mu$ have the same (signed) curvatures.
Therefore, $\mu_1$ and $\mu_2$ lie on the same circular arc $\mu$,
which forms a curvilinear digon with $\gamma$. \end{proof}

\begin{lemsub}\label{lem:chains}
A minimizing swath must be a chain of four-sided components.  The exterior is a band with less pressure than either enclosed region of the swath.
\end{lemsub}

\begin{proof}
By Proposition \ref{prop:Lowest_Pressure_Region_Non-Contractible},
every component of a region of least pressure is non-contractible.
This implies two facts: first, because the swath contains contractible components from two regions, the exterior region is the only region of least pressure, as claimed. Second, because the exterior is a region of least pressure, each component of the exterior region is non-contractible. Also, the \textit{complement} of each component in the exterior is non-contractible, because this complement contains the swath formed by the two enclosed regions. Thus, by Lemma \ref{lem:non-contractible is band}, each component of the exterior is a band. Because the swath formed by the two enclosed regions is connected (Remark \ref{rem:connected}), the exterior region consists of exactly one band, as claimed. Also, the swath cannot contain a contractible cycle of components, because otherwise the cycle would encircle a contractible component of the exterior.

Consider a chain $K$ in the swath.  The components in $K$ must
alternate between two regions $R_1$ and $R_2$ with higher pressure than the exterior. We
will show that the swath contains no additional component of $R_1$
or $R_2$ attached to the chain, as in Figure
\ref{fig:no_lens_in_2-component_chain}.

%%%%%%%%%%%%%%%%%%%%%%%%%%%%%%%%%%%%%%%%%%%%%%%%%%%%%%%%
\begin{figure}[htbp]
\begin{center}
\begin{tabular}{@{\extracolsep{.25 cm}}cc}
\includegraphics*[height=1.5in]{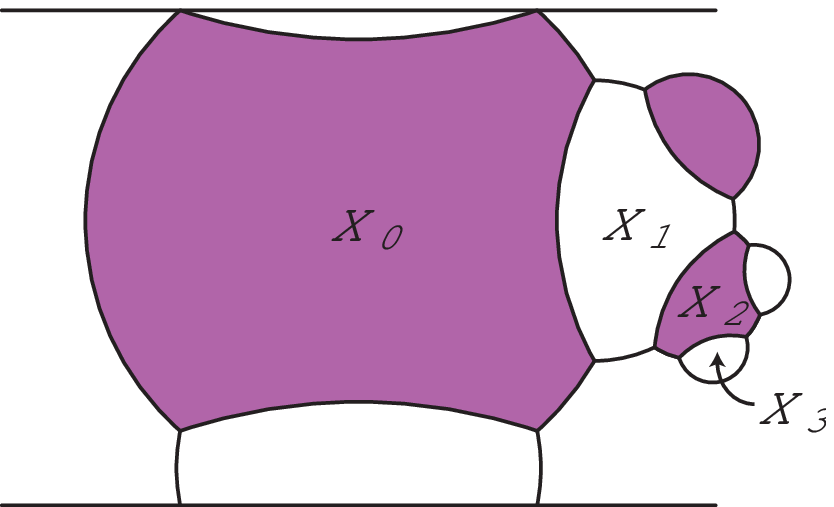}%
&
\includegraphics*[height=1.5in]{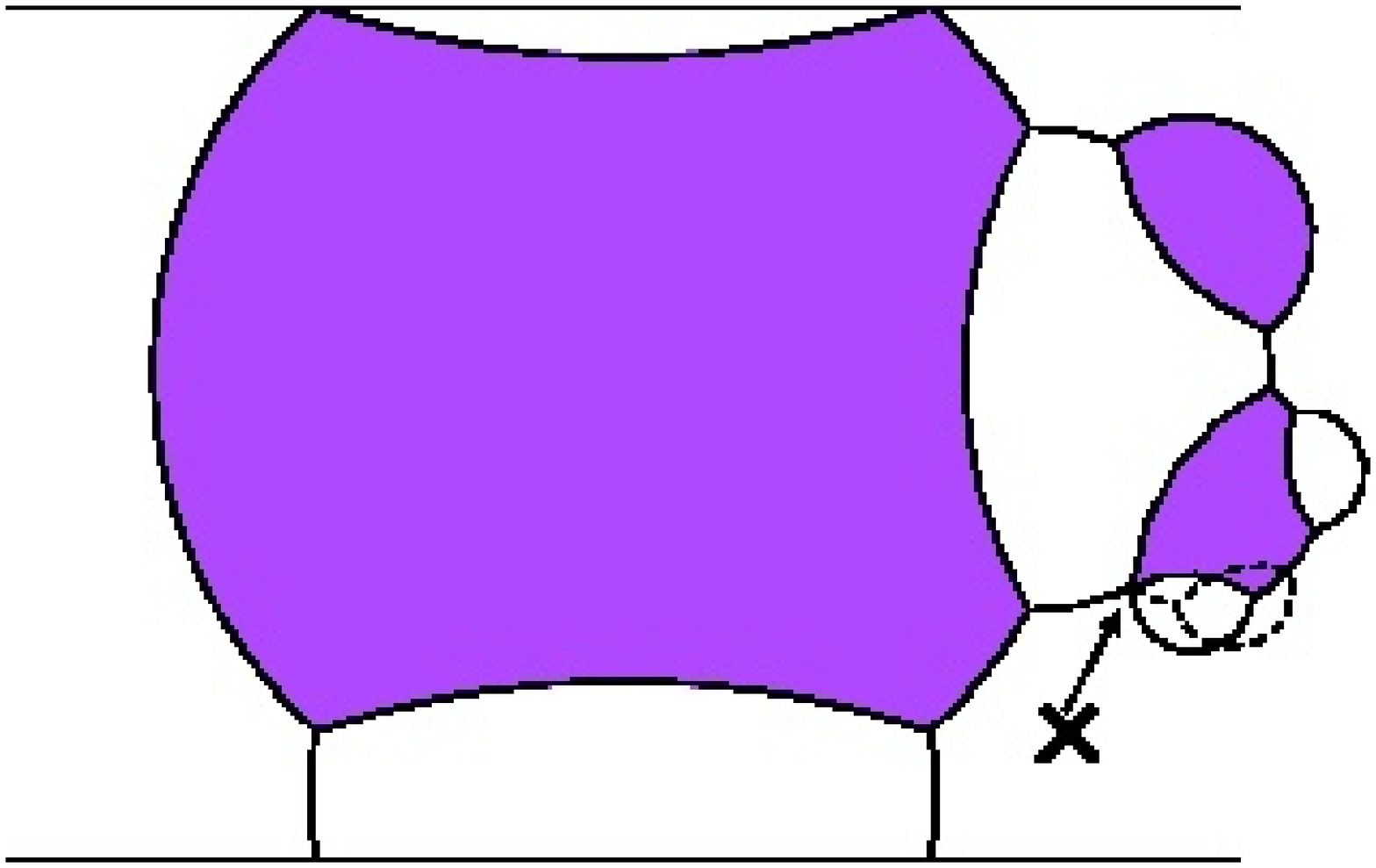}%
\end{tabular}
\caption{\label{fig:no_lens_in_2-component_chain} Extra components
attached to a chain could be slid to create a double bubble
violating regularity.}
\end{center}
\end{figure}
%%%%%%%%%%%%%%%%%%%%%%%%%%%%%%%%%%%%%%%%%%%%%%%%%%%%%%%%

Suppose, for sake of contradiction, that the swath contained an
additional component $X_1$ attached to some component $X_0$ in the
chain. Either $X_1$ has two sides, or another component $X_2$
attaches to it. In this case, either $X_2$ has two sides, or
another component $X_3$ attaches to it. Because the swath contains
no contractible cycle, this process terminates after yielding some
sequence of components $X_0, X_1, \dots, X_k$, where $X_k$ has two
sides. By Lemma \ref{lem:circles_propagate}, the boundary edges of
$X_{k-1}$ meeting the curvilinear digon $X_k$ lie on the same
circle. Therefore, $X_k$ can be slid along the boundary of this
circle until it bumps into the boundary of another component,
creating a minimizer with an illegal singularity (Proposition
\ref{prop:regularity}), a contradiction. Therefore, our original
assumption was false, and no additional components of $R_1$ or
$R_2$ attach to the chain.

Thus, the double bubble is a chain. Because the exterior consists of a single band, each component in the chain is enclosed by four boundary curves: two separating it from the other components of the chain, and two separating it from the exterior. Therefore, the double bubble is a chain of four-sided components.
\end{proof}

\begin{lemsub}\label{lem:circle_chains}
Consider a chain of four-sided components in a perimeter-minimizing double bubble.  If the arcs separating each component from the exterior of the chain lie on the same circle,
then the perimeter of the chain is greater than 3. \end{lemsub}

%%%%%%%%%%%%%%%%%%%%%%%%%%%%%%%%%%%%%%%%%%%%%%%%%%%%%%%%
\begin{figure}[ht]
\begin{center}
\includegraphics[height=2in]{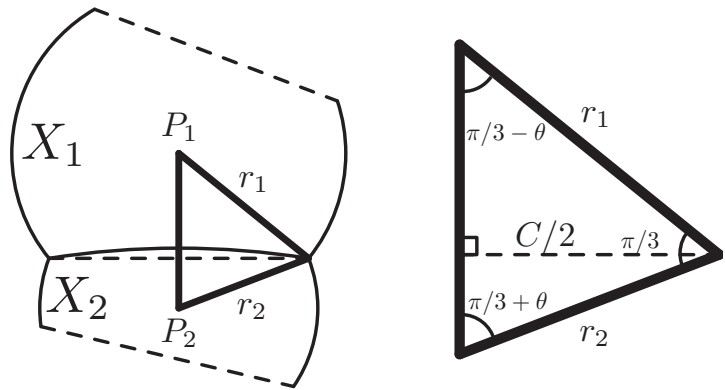}%
\end{center}
\caption{\label{fig:adjacent circle chain components}
Triangulating chains in which opposite arcs are cocircular leads
to a perimeter estimate.}
\end{figure}
%%%%%%%%%%%%%%%%%%%%%%%%%%%%%%%%%%%%%%%%%%%%%%%%%%%%%%%%

\begin{proof} Consider two adjacent components $X_1$ and $X_2$ of the chain, as
drawn in Figure \ref{fig:adjacent circle chain components} so that
the chord $C$ of the arc separating them is horizontal and lies below $X_1$.
Consider the four sides of the polygonal quadrilateral inscribed
on the vertices of $X_1$. By Lemma \ref{lem:chains}, the chain
neighbors a least-pressure region; hence, by Lemma
\ref{lem:High-Pressure Quadrilateral's Low-Pressure Arcs Are
Parallel}, $C$ and the side opposite it in the quadrilateral are
congruent. It follows that $C$ lies below $P_1$, the center of the
circle passing through the vertices of $X_1$. Similarly, $C$ lies
above the center $P_2$ of the circle passing through the vertices
of $X_2$.

We draw the triangle whose vertices are $P_1$, $P_2$, and an endpoint $P_3$ of the chord $C$. The triangle has angles $\pi/3 - \theta, \pi/3, \pi/3 + \theta$ and side lengths $r_1, r_2, P_1P_2$ for
some $\theta, r_1, r_2$. Because $P_1$ and $P_2$ lie on opposite
sides of $C$, angles $\angle P_1P_2P_3$ and $\angle P_2P_1P_3$ are acute, implying that $\theta \in (-\pi/6,
\pi/6)$. Using the Law of Sines in this triangle, we have $r_1 = P_1P_2 \frac{\sin(\pi/3+\theta)}{\sin(\pi/3)}$. Also, $|C|/2 = r_1 \sin(\pi/3 - \theta)$. Hence,
\[
|C| = \frac{4}{\sqrt{3}}P_1P_2 \sin(\pi/3+\theta)\sin(\pi/3-\theta).
\]
Applying the standard transformation $2 \sin x \sin y = \cos(x-y) -
\cos(x+y)$ shows that
\[
|C| = \frac{2}{\sqrt{3}}P_1P_2 (\cos(2\theta) - \cos(2\pi/3)).
\]
With the restriction $\theta \in (-\pi/6, \pi/6)$, we see that $|C|
\ge \frac{2}{\sqrt{3}} P_1P_2$.

As we vary over pairs of adjacent components, the segments
analogous to $\overline{P_1P_2}$ form a closed non-contractible
curve with total length greater than or equal to 1. Hence, the
total length of the chords analogous to $C$ is at least
$\frac{2}{\sqrt{3}} > 1$, implying that the total length of the
corresponding arcs is also greater than one. Furthermore, the arcs
separating the exterior and the chain form two closed
non-contractible curves with total length greater than or equal to
two. Therefore, the total perimeter of the chain is greater than
three. \end{proof}

\begin{lemsub}\label{lem:asymmetric_chains}
No asymmetric chain is perimeter minimizing.
\end{lemsub}

\begin{proof} Suppose, for sake of contradiction, that an
asymmetric chain is minimizing. By Lemma \ref{lem:chains}, it
consists of four-sided components, and the exterior is a
least-pressure region consisting of a single band. By Lemma
\ref{lem:isosceles trapezoid}, each component of the chain has
vertices that form an isosceles trapezoid. If the trapezoid
corresponding to each component is a rectangle, then the
rectangles (and hence the chain) are all symmetric about a closed
geodesic (with the same homology as the chain), a contradiction.
Therefore, one of the components has vertices that form a
non-rectangular isosceles trapezoid. By Lemma \ref{lem:isosceles
trapezoid} and Lemma \ref{lem:High-Pressure Quadrilateral's
Low-Pressure Arcs Are Parallel}, the arcs separating this
component from the least-pressure exterior lie on a single circle.
By Lemma \ref{lem:circles_propagate}, the same is true for all
other components in the chain. Applying Lemma
\ref{lem:circle_chains}, the perimeter of the chain is greater
than 3. By Proposition \ref{prop:perimeter_bound}, the asymmetric
chain cannot be minimizing, a contradiction.
\end{proof}

\begin{lemsub} \label{lem:two_beats_four}
A minimizing symmetric chain must be a standard chain.
\end{lemsub}

\begin{proof}
By Lemma \ref{lem:chains}, such a chain consists of four-sided
components, where alternate components lie in each of two regions.
Thus, the chain contains an even number of components. By Lemma
\ref{lem:isosceles trapezoid}, each four-sided component is
circumscribed about an isosceles trapezoid. Because the chain is
symmetric, each four-sided component must be circumscribed about a
rectangle. It easily follows that any pair of adjacent components
in the chain is congruent to any other pair.

%%%%%%%%%%%%%%%%%%%%%%%%%%%%%%%%%%%%%%%%%%%%%%%%%%%%%%%%
\begin{figure}[htbp]
\begin{center}
\includegraphics*[height=1.5in]{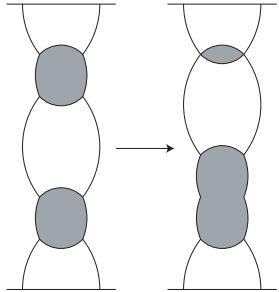}
\end{center}
\caption{\label{fig:rearrange_the_components} Rearranging the
components of a symmetric chain with more than two components
brings about a violation of regularity, proving that a minimizing
chain has two components.}
\end{figure}
%%%%%%%%%%%%%%%%%%%%%%%%%%%%%%%%%%%%%%%%%%%%%%%%%%%%%%%%

Now if the chain contains more than two components, then
rearranging four adjacent components as in Figure
\ref{fig:rearrange_the_components} maintains area and perimeter
while violating regularity, so that the original chain cannot be
perimeter minimizing. Therefore, the chain contains exactly two
components. No curvilinear quadrilateral lies in a least-pressure
region, so the two components of the chain belong to the highest-
and medium-pressure regions. Applying the cocycle condition
(Proposition \ref{prop:regularity}) shows that the chain is
standard.
\end{proof}

\begin{propsub} \label{prop:chain_is_minimizer}
Every perimeter-minimizing swath is a standard chain.
\end{propsub}

\begin{proof}
By Lemma \ref{lem:chains}, a minimizing swath must be a chain of
four-sided components. By Lemma \ref{lem:asymmetric_chains} a
minimizing chain must be symmetric.  By Lemma
\ref{lem:two_beats_four}, a minimizing symmetric chain must be a
standard chain.
\end{proof}

\begin{remsub} \label{rem:standard chain short direction}
While we will show that a minimizing band lens must have the same
homology as a shortest closed geodesic (Lemma
\ref{lem:bl_short_direction}), we do not have a proof of the
analogous fact for minimizing standard chains.  However,
\emph{Mathematica} plots for perimeter suggest that this property
holds at least for the particular tori described by Figure
\ref{fig:param-diags}.

\end{remsub}

%%%%%%%%%%%%%%%%%%%%BBLOB%%%%%%%%%%%%%%%%%%
\subsection{A Single Band Adjacent to a Contractible Set of Components}
\label{SS:double_bubbles_with_single_band}

Proposition \ref{prop:band lens_is_minimizer} uses regularity
(Proposition \ref{prop:regularity}) and Proposition
\ref{prop:Lowest_Pressure_Region_Non-Contractible} to show that
the only potential minimizer in the class of double bubbles with a
single band adjacent to a contractible set of components is the
band lens.

Throughout this section, we will consider a proposed minimizer
with regions $R_1$, $R_2$ and the exterior, such that $R_1$ is
contractible and $R_2$ contains a single band (and possibly some contractible components).

\begin{lemsub} \label{lem:p_band_equals_p_ext_in bblob}
In a minimizing double bubble with some region $R_1$ contractible
and another region $R_2$ containing a single band, the pressure
difference between $R_2$ and the third (exterior) region must be
zero; that is, all curves separating $R_2$ from the exterior must
be straight line segments.
\end{lemsub}

\begin{proof}
Consider a curve separating the band in $R_2$ from the exterior.
If it is a straight line segment, the proof is complete, since all
curves separating $R_2$ from the exterior have the same curvature
by regularity (Proposition \ref{prop:regularity}).

If the curve has nonzero curvature, continuously
straighten the curve to reduce perimeter, while sliding one component of the boundary of the band in $R_2$
in order to maintain area. Either we will
be able to straighten all such pieces while reducing perimeter, or
we will first cause a violation of regularity.
\end{proof}

\begin{lemsub}\label{lem:r1_highest_pressure}
In a minimizing double bubble with some region $R_1$ contractible
and another region $R_2$ containing a single band, $R_1$ must have
higher pressure than each of $R_2$ and the exterior region.
\end{lemsub}

\begin{proof}
By Lemma \ref{lem:p_band_equals_p_ext_in bblob}, $R_2$ and the
exterior have the same pressure.  Since the double bubble is
neither a tiling nor the double band, $R_1$ cannot also have the
same pressure (Proposition \ref{prop:three_equal_pressures}).
Moreover, by Proposition
\ref{prop:Lowest_Pressure_Region_Non-Contractible}, since $R_1$ is
contractible it cannot be a region of least pressure. Thus, $R_1$
must have a higher pressure than both $R_2$ and the exterior.
\end{proof}

\begin{corsub}\label{cor:only high pressure in blob}
All contractible components of the double bubble plus exterior
belong to $R_1$, and $R_2$ consists of a single band.
\end{corsub}

\begin{proof}
By Lemma \ref{lem:p_band_equals_p_ext_in bblob} and Lemma
\ref{lem:r1_highest_pressure}, $R_2$ and the exterior have equal
pressure, less than that of $R_1$.  Thus, by Proposition
\ref{prop:Lowest_Pressure_Region_Non-Contractible}, their
components must be non-contractible.  Hence only components of
$R_1$ can be contractible, and $R_2$ consists of the only
non-contractible component of the double bubble --- a band.
\end{proof}

\begin{lemsub}
\label{lem:bblob_contractible_are_lenses} In a minimizing double
bubble with some region $R_1$ contractible and another region
$R_2$ containing a single band, all components of $R_1$ are
lenses.
\end{lemsub}

\begin{proof}
By Lemma \ref{lem:r1_highest_pressure}, $R_1$ has higher pressure
than both $R_2$ and the exterior. By Remarks \ref{rem:pressure}
and \ref{rem:alternating_regions}, all components of $R_1$ have
either two or four sides. We will assume some component of $R_1$
is a curvilinear quadrilateral $Z$, and show that this leads to a
contradiction. Let the four sides of $Z$ be labelled clockwise $x_1$, $y_1$, $x_2$, $y_2$, where $x_1, x_2$ border $R_2$ (a single band, by Corollary \ref{cor:only high pressure in blob}) and $y_1, y_2$ border $Y_1, Y_2$ in the exterior.

If $x_1$ and $x_2$ lie on the same side of $R_2$, then either $y_1$ or $y_2$ lies along that same side of the band; the corresponding component $Y_1$ or $Y_2$ is contractible, contradicting Corollary \ref{cor:only high pressure in blob}. Otherwise, $x_1$ and $x_2$ lie on opposite sides of $R_2$. In this case, the closure of $R_2 \cup Z$ has contractible complement, so that both $Y_1$ and $Y_2$ are contractible --- again contradicting Corollary \ref{cor:only high pressure in blob}.

Thus, all components of $R_1$ are curvilinear digons. Now, by Lemma \ref{lem:p_band_equals_p_ext_in bblob}, $R_2$ and the exterior have equal pressure.  Therefore, all curves bounding components of $R_1$ have equal, nonzero curvature, and so all curvilinear digons in $R_1$ must be lenses. \end{proof}

\begin{propsub} \label{prop:band lens_is_minimizer}
A perimeter-minimizing double bubble with a single band adjacent to a contractible set of components must be the band lens.
\end{propsub}

\begin{proof}
%%%%%%%%%%%%%%%%%%%%%%%%%%%%%%%%%%%%%%%%%%%%%%%%%%%%%%%%%%%%%%%%
\begin{figure}[htbp]
\includegraphics*{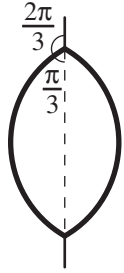}%
\caption{\label{fig:lens diameter} Any lens from $R_1$ must be
embedded in a closed geodesic.}
\end{figure}
%%%%%%%%%%%%%%%%%%%%%%%%%%%%%%%%%%%%%%%%%%%%%%%%%

By Corollary \ref{cor:only high pressure in blob} and Lemma
\ref{lem:bblob_contractible_are_lenses}, one region $R_1$ consists
only of lenses and the other region $R_2$ consists of a single
band. Corollary \ref{cor:only high pressure in blob} also implies
that the exterior contains no contractible components. It follows
from Lemma \ref{lem:non-contractible is band} that the exterior
consists solely of bands. Because the double bubble is connected (Remark \ref{rem:connected}),
the exterior consists of a single band.

Each lens in $R_1$ has two vertices, each the endpoint of some
curve separating $R_2$ and the exterior. By Lemma
\ref{lem:p_band_equals_p_ext_in bblob}, these curves are straight
line segments. Furthermore, both segments must lie on the line
passing through the diameter of the lens, as depicted in Figure
\ref{fig:lens diameter}. Hence, the diameters of the lenses,
together with the segments separating $R_2$ from the exterior,
form closed geodesics. Since $R_2$ and the exterior each consist
of exactly one band, the boundary of the double bubble consists of
exactly two closed geodesics, with lenses from $R_1$ embedded in
at least one of them. Each of these two closed geodesics has
length $L' \ge 1$.

Let the lenses in $R_1$ have total perimeter $p$ and let their
diameters have total length $d$. If $d \ge 1$, then the perimeter
of the double bubble is at least $(2L' - d) + p > (2L' - d) + 2d
\ge 3$, so by Proposition \ref{prop:perimeter_bound} the double
bubble is not perimeter minimizing. Otherwise, we can move all the
lenses so that they are adjacent to each other and embedded in the
same closed geodesic, creating a new double bubble enclosing and
separating the same areas. If there is more than one lens, then
the new double bubble violates regularity (Figure \ref{fig:merging
blobs}), a contradiction. Therefore, the original minimizing
double bubble contains exactly one lens, and so must be the band
lens.

%%%%%%%%%%%%%%%%%%%%%%%%%%%%%%%%%%%%%%%%%%%%%%%%%%%%%%%%%%%%%%%%%%%%%%%%%%%%%%%%%
\begin{figure}[htbp]
\includegraphics*{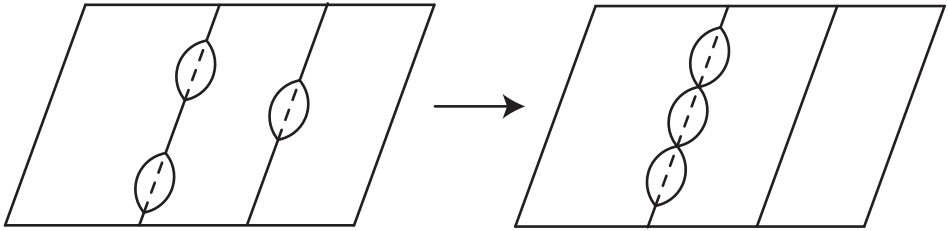}%
\caption{\label{fig:merging blobs} Multiple lenses can be slid
together to violate regularity. Therefore, there can only be one
lens.}
\end{figure}
%%%%%%%%%%%%%%%%%%%%%%%%%%%%%%%%%%%%%%%%%%%%%%%%%%%%%%%%%%%%%%%%%
\end{proof}

\begin{lemsub} \label{lem:bl_short_direction}
In a minimizing band lens, the band must have the same homology as
a shortest closed geodesic. \end{lemsub}

\begin{proof}
We show that any pair of areas that can be enclosed by a band lens
wrapping around a longer direction of the torus can be enclosed by
an embedded band lens wrapping around a short direction. Suppose
not. By Proposition \ref{prop:perimeter_bound}, the perimeter of
the original band lens is less than or equal to three, so the
perimeter of the lens is less than two. Therefore, the diameter of
the lens is less than one, and the lens can be embedded on the
torus with its axis along a short direction (Remark
\ref{rem:lens-fits}). Thus, if the band lens wrapping around the
short direction does not exist, the closed geodesic not meeting
the lens cannot get close enough to make the band thin enough. But
if this happens in the short direction, then it will certainly
happen for the band lens wrapping around the longer direction, a
contradiction.

Hence, if there is a minimizing band lens wrapping around a longer
direction of the torus --- say, parallel to a closed geodesic of
length $L_0 > 1$
--- there is a band lens in a short direction of the torus
enclosing the same areas. However, the new band lens has perimeter
$2(L_0-1)$ less than the initial band lens, a contradiction.
Therefore, any minimizing band lens has the same homology as a
shortest closed geodesic.
\end{proof}

\subsection{Tilings}\label{SS:tilings}

This section will show that every perimeter-minimizing tiling is a
standard hexagon tiling (Proposition \ref{prop:only tiling is 3
hex}). Corollary \ref{cor:hex-or-oct} showed that if a tiling is
minimizing, then it is a hexagon tiling or an octagon-square
tiling. Proposition \ref{prop:only hex is standard hex} shows that
every minimizing hexagon tiling is standard, and Proposition
\ref{prop:oct-square bad} shows that no octagon-square
tiling is perimeter minimizing. \\

\paragraph{\textbf{Hexagon Tilings}}

Lemmas \ref{lem:small hexagons lose}, \ref{lem:perimeter outside
big hexagon}, and \ref{lem:6 or more hexes is bad} show that a
minimizing hexagon tiling divides the torus into exactly three
hexagonal components. Lemma \ref{lem:trans preserves perimeter} shows
that such a tiling can be transformed without changing perimeter
to make the three hexagons congruent translations of each other.
By Lemma \ref{lem:small hexagons lose}, the perimeter of the new
hexagon tiling is greater than 3 except on the hexagonal torus,
implying that the original minimizing hexagon tiling lies on the
hexagonal torus. Finally, Proposition \ref{prop:only
hex is standard hex} shows that any such hexagon tiling is the standard
hexagon tiling.

\begin{lemsub} \label{lem:multiple of three}
The number of components in a perimeter-minimizing hexagon tiling plus exterior
is a multiple of three.
\end{lemsub}

\begin{proof} By Proposition \ref{prop:regularity}, there are finitely
many vertices in the tiling --- say, $v$ --- each bordering three
distinct regions. Consider any of the three regions $R_i$. Each of
its components borders six distinct vertices in the tiling, and
each vertex in the tiling borders exactly one component of $R_i$.
Hence, $R_i$ contains $v/6$ components. Therefore, $v/6$ is an
integer, and there are $3(v/6)$ components in all. \end{proof}

\begin{lemsub} \label{lem:anti-smoothing}
Given $n \ge 3$ real numbers $a_1, \dots, a_n \in [0, A/3]$ such
that $\sum_{i=1}^n a_i = A$, we have $\sum_{i=1}^n \sqrt{a_i} \ge
3\sqrt{A/3}$.
\end{lemsub}

\begin{proof}
Because $\sqrt{A}$ is strictly concave, $\sum_{i=1}^n \sqrt{a_i}$ attains
its minimum at some vertex of the given domain. At each vertex, three of the $a_i$ equal $A/3$,
and the other $a_i$ equal 0. The result follows immediately. \end{proof}

\begin{lemsub} \label{lem:small hexagons lose}
In a hexagon tiling on a torus with area $A$
such that each component has area at most $A/3$, the perimeter of
the tiling is greater than or equal to 3. Equality can hold
only for a tiling on a hexagonal torus with three components in the tiling plus exterior.
\end{lemsub}

\begin{proof}
By Lemma \ref{lem:multiple of three}, if there are $n$ components
in the tiling plus exterior, then $n \ge 3$. Let $a_1, \dots, a_n$
be the areas of the $n$ components $C_1, \dots, C_n$, and let
$p_1, \dots, p_n$ be their perimeters. By the isoperimetric
inequality for hexagons, $p_i \ge \sqrt{8\sqrt{3}}\sqrt{a_i}$
for each $i$, with equality if and only if $C_i$ is a
regular hexagon.

The total perimeter of the tiling, $\frac{1}{2} \sum_{i=1}^{n}
p_i$, is thus greater than or equal to $\sqrt{2\sqrt{3}}
\sum_{i=1}^{n} \sqrt{a_i}$. By Lemma \ref{lem:anti-smoothing},
$\sqrt{2\sqrt{3}} \sum_{i=1}^{n} \sqrt{a_i}$ is in turn
greater than or equal to $\sqrt{2\sqrt{3}} \cdot 3\sqrt{A/3} =
\sqrt{6\sqrt{3}A}$.

Because $A \ge \sqrt{3}/2$ with equality only for the hexagonal
torus (Remark \ref{rem:torus-angle-and-area}), the perimeter of
the tiling is at least $\sqrt{6\sqrt{3}A} \ge \sqrt{6\sqrt{3}
\frac{\sqrt{3}}{2}} = 3$. The perimeter is exactly 3 only if
equality holds in all the intermediate inequalities.  In other
words, the torus must have area $\sqrt{3}/2$, which only holds for
the hexagonal torus; each component must be a regular hexagon and
there must be 3 components each with area $A/3$.
\end{proof}

\begin{remsub}
This method of comparing hexagon tilings to tilings of three equal
areas was inspired by the result of Hales that the regular hexagon
is the most efficient way to tile the plane --- not necessarily with polygons --- into unit areas [H]. Earlier, Fejes-T\'{o}th proved the
result for polygonal tilings ([FT1, Chapter III, Section 9, p. 84]
or [FT3, Section 26, Corollary, p. 183] after [FT4]), a result more
comparable than Hales' to ours.
\end{remsub}

\begin{lemsub}\label{lem:perimeter outside big hexagon} Consider a hexagon tiling with $n \ge 6$ components in the tiling plus exterior. For each component $H$ in any of the three regions, the perimeter of the tiling exceeds 1 plus the perimeter of $H$. \end{lemsub}

%%%%%%%%%%%%%%%%%%%%%%%%%%%%%%%%%%%%%%%%%%%%%%%%%%%%%%%%%%%%
\begin{figure}[ht]
\begin{center}
\includegraphics*{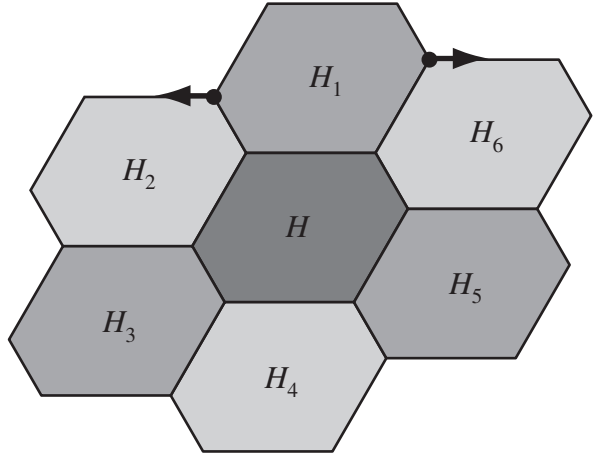}%
\end{center}
\caption{\label{fig:seven_hexagons} A hexagon $H$ on the torus and
its six neighbors, shaded according to which regions they lie in.
(Some of the six hexagons neighboring $H$ in this picture may
actually represent the same component.) In order to find a
non-contractible curve not intersecting $H$, we avoid travelling
along a downward edge from either of the two marked vertices.}
\end{figure}
%%%%%%%%%%%%%%%%%%%%%%%%%%%%%%%%%%%%%%%%%%%%%%%%%%%%%%%%%%%%

\begin{proof}
We orient the torus so that the edges of the tiling form angles of
0 and $\pm \pi/3$ with the horizontal. Consider $H$ along with the
hexagons neighboring it, labelled $H_1, \dots, H_6$ as in Figure
\ref{fig:seven_hexagons}.

Although some of $H_1, \dots, H_6$ may coincide (if $H$ and some component share more than one edge), we claim that
$H_2$, $H_4$, and $H_6$ cannot \textit{all} coincide. Assume not.
$H_3$ borders $H_2$ from below, and $H_5$ borders $H_6$ from
below; because $H_2 = H_6$, we have $H_3 = H_5$. We can likewise
show that $H_1$, $H_3$, and $H_5$ are all the same component. The
top edge, bottom-right, and bottom-left edges of this component
border $H_2 = H_4 = H_6$ (since, as shown in Figure
\ref{fig:seven_hexagons}, the top edge of $H_3$ borders $H_2$, the
bottom-right edge of $H_1$ borders $H_6$, and the bottom-left edge
of $H_5$ borders $H_4$); similarly, the other edges of this
component border $H$. Thus, $H_1$ only neighbors two distinct
components: $H$ and $H_2$. Likewise, $H$ and $H_2$ must only
neighbor components in $\{H, H_1, H_2\}$. However, this is
impossible: at least one of the $n-3$ hexagons not in $\{H, H_1,
H_2\}$ must neighbor some hexagon in $\{H, H_1, H_2\}$. Therefore,
our original assumption was false, and we do not have $H_2 = H_4 =
H_6$.

Without loss of generality, assume that $H_2 \ne H_6$. Notice also
that $H_2 \ne H_5$ and $H_6 \ne H_3$, since alternate sides of $H$
must border different regions (Remark
\ref{rem:alternating_regions}).

%%%%%%%%%%%%%%%%%%%%%%%%%%%%%%%%%%%%%%%%%%%%%%%%%%%%%%%%%%%%
\begin{figure}[ht]
\begin{center}
\begin{tabular} {@{\extracolsep{.5 cm}} cc}
\includegraphics*[height=1in]{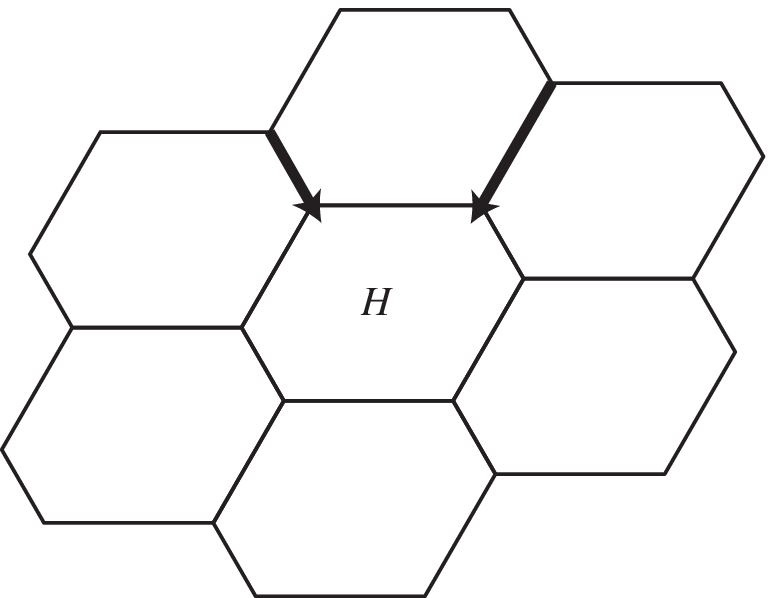}%
&
\includegraphics*[height=1in]{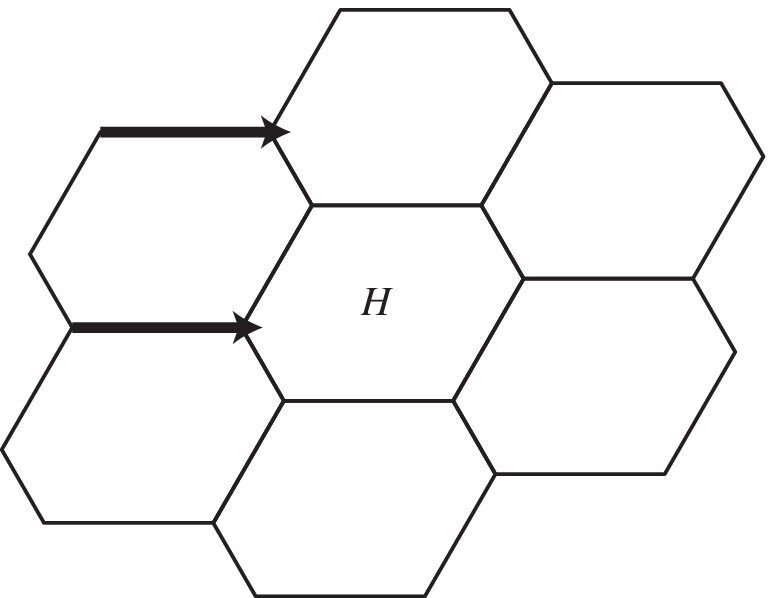}%
\\
(a) & (b)
\end{tabular}
\\
\includegraphics*[height=1in]{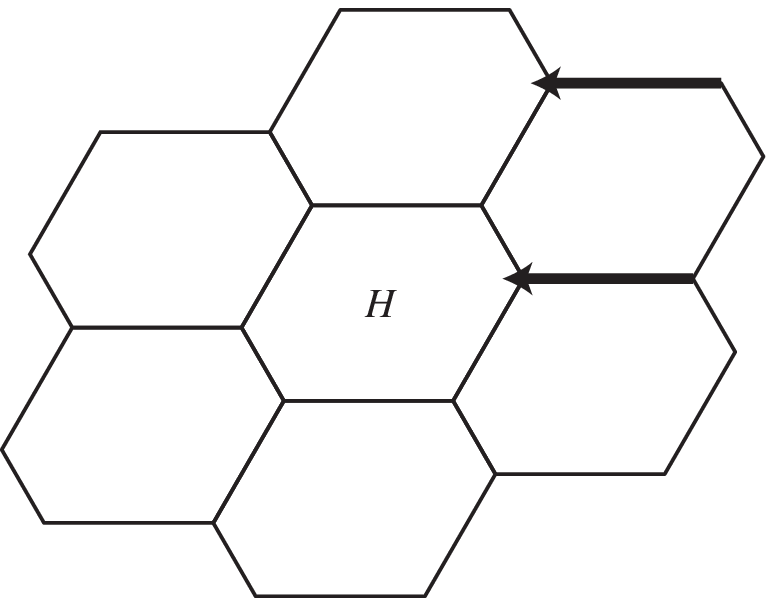}%
\\
(c)
\end{center}

\caption{\label{fig:impossible traps along hexagon journey} In
order to bound the perimeter outside any component, we travel
along the boundary in an attempt to form a closed non-contractible
path. These diagrams show various movements that we do not
perform during our journey.}
\end{figure}
%%%%%%%%%%%%%%%%%%%%%%%%%%%%%%%%%%%%%%%%%%%%%%%%%%%%%%%%%%%%

The tiling contains $6n/3$ vertices, which is at least twelve, so
we may pick a vertex $v$ not on the boundary of $H$. Starting at
$v$, we travel as follows: if we are at one of the vertices
marked in Figure \ref{fig:seven_hexagons}, travel horizontally as shown in Figure \ref{fig:seven_hexagons}; otherwise, travel downward along an edge. It follows that we never travel downward along the edges depicted in Figure
\ref{fig:impossible traps along hexagon journey}(a). Furthermore,
because $H_6 \not\in \{H_2, H_3\}$, we never travel to the right
along the edges depicted in Figure \ref{fig:impossible traps along
hexagon journey}(b). Similarly, we never travel to the left along
the edges depicted in Figure \ref{fig:impossible traps along
hexagon journey}(c). It easily follows that we never travel along
a horizontal edge only to immediately backtrack along it, and that
we never travel along the boundary of $H$ or along any of the
edges stemming from $H$.

%%%%%%%%%%%%%%%%%%%%%%%%%%%%%%%%%%%%%%%%%%%%%%%%%%%%%%%%%%%%
\begin{figure}[ht]
\begin{center}
\includegraphics*{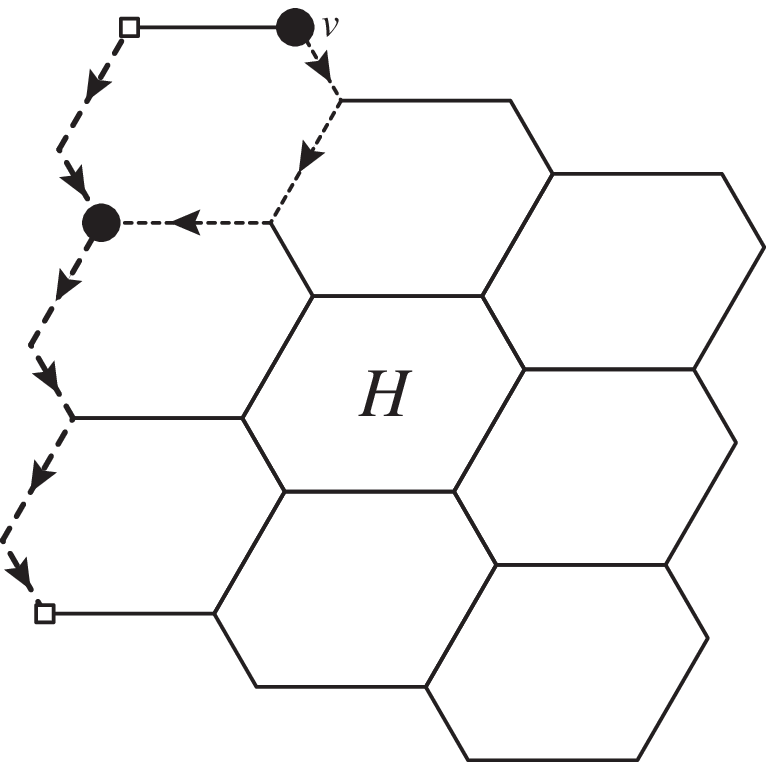}%
\end{center}
\caption{\label{fig:sample hexagon journey} Starting at some vertex $v$ and travelling along downward edges --- except in those cases where doing so would lead
us to a vertex of $H$ --- yields a closed, non-contractible curve
not passing through the boundary of $H$, for a total perimeter
greater than 3.}
\end{figure}
%%%%%%%%%%%%%%%%%%%%%%%%%%%%%%%%%%%%%%%%%%%%%%%%%%%%%%%%%%%%

There are finitely many vertices in the tiling (Proposition
\ref{prop:regularity}), implying that some portion of our journey
begins and ends at the same vertex. Because we never backtrack
along a horizontal edge, and because we never travel upward along
an edge, this portion of our journey occurs along a closed,
non-contractible curve (Figure \ref{fig:sample hexagon journey}).
This curve has length at least 1, and it does not overlap the
boundary of $H$. Nor does it overlap any of the edges stemming
from $H$, which have positive length. Hence, the boundary curves
not enclosing $H$ have total length greater than $1$, as desired.
\end{proof}

\begin{lemsub} \label{lem:6 or more hexes is bad}
No hexagon tiling with six or more components in the double bubble plus exterior is
perimeter minimizing.
\end{lemsub}

\begin{proof}
We consider a hexagon tiling with six or more components. We
first consider the case when the area $\alpha$ of some component
is at least one-third the total area $A$ of the torus, where $A
\ge \sqrt{3}/2$ (Remark \ref{rem:torus-angle-and-area}). By the
isoperimetric inequality for hexagons, the perimeter of $H$ is at
least $\sqrt{8\sqrt{3}\alpha} \ge \sqrt{8\sqrt{3}(A/3)} \ge 2$.
Therefore, by Lemma \ref{lem:perimeter outside big hexagon}, the
perimeter of the tiling is greater than $2 + 1 = 3$.

Next consider the case in which each component has area less than
$A/3$. Because there are more than three components, by Lemma
\ref{lem:small hexagons lose}, the tiling has perimeter greater
than 3.

In both cases, it follows by our perimeter bound (Proposition
\ref{prop:perimeter_bound}) that the hexagon tiling is not
perimeter minimizing.
\end{proof}

\begin{remsub}
We first attempted to prove Lemma \ref{lem:6 or more hexes is bad}
by shrinking and expanding hexagons to preserve perimeter and areas
and eventually contradict regularity. Although the perimeter of
the tiling remains the same under such variations, the problem ---
as Wacharin Wichiramala explained to us --- is that it is not always
possible to preserve the regions' areas in order to reach such a contradiction.
\end{remsub}

\begin{lemsub} \label{lem:hex_sidelength_relation}
In a hexagon with interior angles of $2\pi/3$, the sum of the lengths of
any two adjacent sides, equals the sum of the lengths
of the opposite pair of adjacent sides.
\end{lemsub}

\begin{proof}
Orient the hexagon so that the two sides $p_1, p_2$ not being
summed are horizontal, with one pair of remaining sides $q_1, q_2$
on the left and the other pair of remaining sides $r_1, r_2$ on
the right. Let $d$ be the distance between $p_1$ and $p_2$. The
sides $q_1, q_2$ make an angle of $\pi/6$ with the vertical, so
the sum of the lengths of $q_1, q_2$ is $d/\sin(\pi/3)$. Similarly, the sum of
the lengths of $r_1, r_2$ is $d/\sin(\pi/3)$. Hence, the two sums
are equal, as desired.
\end{proof}

\begin{lemsub} \label{lem:trans preserves perimeter}
Given a hexagon tiling with three components in the tiling plus exterior, there exists a second hexagon tiling with the following properties: the edges of the new tiling are parallel to the edges of the original tiling; the perimeters of the two tilings are equal; and any two parallel edges of the new tiling are congruent. In particular, there exists a second hexagon tiling with the same perimeter, which still contains three components (all of the same area) in the tiling plus exterior. \end{lemsub}

%%%%%%%%%%%%%%%%%%%%%%%%%%%%%%%%%%%%%%%%%%%%%%%%%%%%%%%%%%%%
\begin{figure}[ht]
\begin{center}
\includegraphics*[height=1.5in]{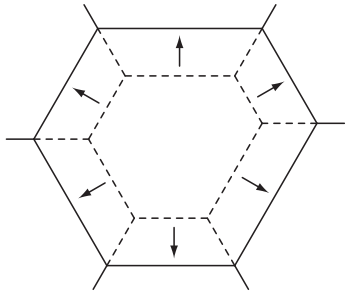}%
\end{center}
\caption{\label{fig:hexvariation} Expanding any component in a
hexagon tiling or its exterior generates a double bubble enclosing
different areas, but does not change perimeter. Hexagons can also
be shrunk in a similar fashion.}
\end{figure}
%%%%%%%%%%%%%%%%%%%%%%%%%%%%%%%%%%%%%%%%%%%%%%%%%%%%%%%%%%%%

%%%%%%%%%%%%%%%%%%%%%%%%%%%%%%%%%%%%%%%%%%%%%%%%%%%%%%%%%%%%
\begin{figure}[ht]
\begin{center}
\includegraphics*[height=2in]{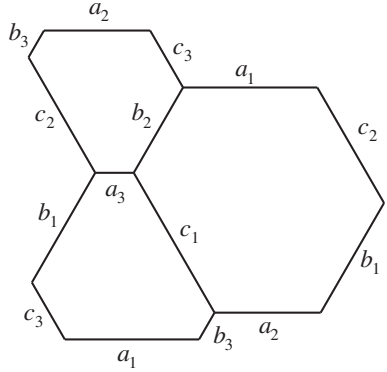}%
\end{center}
\caption{\label{fig:three-hexagon labels} A hexagon
tiling with three components in the tiling plus exterior.}
\end{figure}
%%%%%%%%%%%%%%%%%%%%%%%%%%%%%%%%%%%%%%%%%%%%%%%%%%%%%%%%%%%%

\begin{proof}
To get the desired result, we will find a variation that shrinks
and expands hexagons as in Figure \ref{fig:hexvariation} to make
the three areas equal.  Perimeter is unchanged, since $dP/dA =
\textrm{curvature} = 0$ through the whole family of double bubbles
in the variation.  We show that such a variation exists, and we
also give an alternative geometric explanation for why perimeter
remains unchanged.

Let $a_i, b_i, c_i$ ($1 \le i \le 3$) be the lengths of the edges
in the hexagon tiling as in Figure \ref{fig:three-hexagon labels};
for each $i \in \{1, 2, 3\}$, the sides of some hexagon $H_i$ are
assigned the lengths $a_{i-1}, b_{i-1}, c_{i-1}, a_{i+1}, b_{i+1},
c_{i+1}$ (where indices are taken modulo 3). In this proof, as we
modify the hexagon tiling, the values of the $a_i, b_i, c_i$ will
change accordingly.

%%%%%%%%%%%%%%%%%%%%%%%%%%%%%%%%%%%%%%%%%%%%%%%%%%%%%%%%%%%%
\begin{figure}[ht]
\begin{center}
\includegraphics*{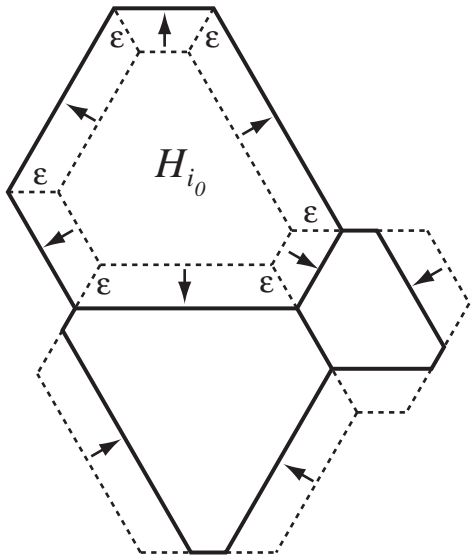}
\end{center}
\caption{\label{fig:three-hexagon expansion} In a hexagon tiling
with three components in the double bubble plus exterior, we may
expand (or shrink) one hexagon so that six edges increase by
$\epsilon$ while three edges decrease by $2\epsilon$.}
\end{figure}
%%%%%%%%%%%%%%%%%%%%%%%%%%%%%%%%%%%%%%%%%%%%%%%%%%%%%%%%%%%%

Fix some subscript $i_0$. As shown in Figure
\ref{fig:three-hexagon expansion}, for some range of real numbers
$\epsilon$ (which we will describe shortly), we may transform our
hexagon tiling into another hexagon tiling, so that the edge
lengths change as follows: $a_{i_0}, b_{i_0}, c_{i_0}$ change by
$-2\epsilon$, and the other $a_i, b_i, c_i$ change by $\epsilon$.
More specifically, the range of such $\epsilon$ are those
$\epsilon$ for which the side lengths of the resulting tiling are
positive: namely, such that $a_i+\epsilon$, $b_i+\epsilon$, and
$c_i+\epsilon$ are positive for $i \ne i_0$, and such that $a_{i_0}-2\epsilon$,
$b_{i_0}-2\epsilon$, and $c_{i_0}-2\epsilon$ are positive. Adding
up the increases and decreases in perimeter shows that under this
transformation, the perimeter of the hexagon tiling remains
unchanged.

Without loss of generality, assume that $a_1 \ge a_2 \ge a_3$
initially. By Lemma \ref{lem:hex_sidelength_relation}, we have
$(a_2 - a_3)/3 = (b_2 - b_3)/3 = (c_2 - c_3)/3$; set $\epsilon$
equal to this common value. As described in the last paragraph, we
may shrink $H_3$ so that its edge lengths ($a_i, b_i, c_i$ for $i
\ne 3$) decrease by $\epsilon$ and so that $a_3, b_3, c_3$
increase by $2\epsilon$. It is easy to check that this $\epsilon$
is in the allowed range. Furthermore, after this transformation,
we have $(a_2, b_2, c_2) = (a_3, b_3, c_3)$.

At this point, we now have $a_1 \ge a_2 = a_3$. We may now set
$\epsilon$ equal to $\frac{a_1-a_2}{3}$ and expand $H_1$ so that
its edge's lengths increase by $\epsilon$. In the resulting
hexagon tiling, we have $(a_1, b_1, c_1) = (a_2, b_2, c_2) = (a_3,
b_3, c_3)$, as desired. Furthermore, because the two
transformations we used did not change perimeter, the resulting
tiling and the initial tiling have the same perimeter. This
completes the proof.
\end{proof}

\begin{lemsub} \label{lem:cong_hexes_on_hex_torus}
Suppose that the hexagonal torus is partitioned into three
hexagons with interior angles of $2\pi/3$, where the hexagons are
congruent translations of each other. Then the hexagons are
regular hexagons of side length $\frac{1}{3}$, with sides parallel
to the three short directions of the torus. The total perimeter of
the partition is equal to three. \end{lemsub}

%%%%%%%%%%%%%%%%%%%%%%%%%%%%%%%%%%%%%%%%%%%%%%%%%%%%%%%%%%%%
\begin{figure}[ht]
\includegraphics*{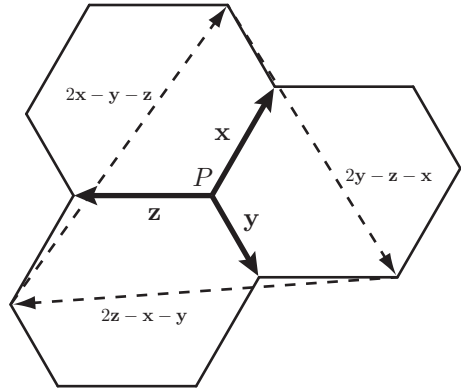}
\caption{\label{fig:three_cong_hexes_fig} In partitions of a
torus into three hexagons that are congruent
translations of each other, examining three closed geodesics
(represented by the dotted vectors) shows that if the hexagons are
not regular, then the torus is not hexagonal.}
\end{figure}
%%%%%%%%%%%%%%%%%%%%%%%%%%%%%%%%%%%%%%%%%%%%%%%%%%%%%%%%%%%%

\begin{proof}
Based on the partition, we draw a fundamental domain of the torus
as in Figure \ref{fig:three_cong_hexes_fig}. Let $P$ be the vertex
in the center of the figure where the three hexagons meet; let
$\mbf{x}$, $\mbf{y}$, $\mbf{z}$ be the vectors from $P$ outward
along the edges of the hexagon tiling, with magnitudes $x, y, z$.

The three dotted vectors all begin and end at the same point $Q$
on the torus. We tile the plane with the fundamental domain in
Figure \ref{fig:three_cong_hexes_fig}, so that each point on the
plane represents a point on the torus. We can travel between any
two points that represent $Q$ by travelling along the three dotted
vectors. Hence, any two of the dotted vectors lie along the sides
of a parallelogram fundamental domain. Now, some two of the dotted
vectors $\mathbf{v}_1$, $\mathbf{v}_2$ meet at an angle $\theta
\in [\pi/3,\pi/2]$. The fundamental domain they define has area
$\|\mathbf{v}_1\|\|\mathbf{v}_2\|\sin\theta = \sqrt{3}/2$ (the
area of the hexagonal torus). Then equality must hold in the inequalities $\|\mathbf{v}_1\|,
\|\mathbf{v}_2\| \ge 1$ and $\sin\theta \ge \sqrt{3}/2$; that is,
$\mathbf{v}_1$ and $\mathbf{v}_2$ must have magnitude 1 and meet at an angle equal to $\pi/3$. It follows that
the three dotted vectors must all have magnitude 1 and lie along
short directions of the torus.

Specifically, $2\mbf{x} - \mbf{y} - \mbf{z}$ has magnitude 1.
Hence,
\begin{eqnarray*}
1 &=& (2\mbf{x} - \mbf{y} - \mbf{z}) \cdot (2\mbf{x} - \mbf{y} - \mbf{z}) \\
&=& 4\mbf{x} \cdot \mbf{x} + \mbf{y} \cdot \mbf{y} + \mbf{z} \cdot \mbf{z} - 4 \mbf{x} \cdot \mbf{y} - 4\mbf{x} \cdot \mbf{z} + 2\mbf{y} \cdot \mbf{z} \\
&=& 4x^2 + y^2 + z^2 - 4xy\cos(2\pi/3) - 4xz\cos(2\pi/3) + 2yz\cos(2\pi/3) \\
&=& (x+y+z)^2 + 3x^2 - 3yz,
\end{eqnarray*}
or
\begin{equation}
1 - (x+y+z)^2 = 3(x^2 - yz). \label{eqn:cong_hexes_eqn}
\end{equation}
Similarly, $1 - (x + y + z)^2 = 3(y^2 - zx)$. Thus, $x^2 - yz =
y^2 - zx$, or $(x-y)(x+y+z) = 0$, implying that $x = y$. Likewise, $y = z$. By Equation
\ref{eqn:cong_hexes_eqn}, it follows that $x = y = z =
\frac{1}{3}$. Therefore, the partition consists of three regular
hexagons of side length $\frac{1}{3}$. It easily follows that the
sides of the hexagons are parallel to the three short directions
of the torus, and that the perimeter of the partition is three.
\end{proof}

\begin{corsub} \label{cor:three_hexes_are_standard} If a
hexagon tiling on the hexagonal torus contains three components in
the tiling plus exterior, then the tiling is a standard hexagon
tiling.
\end{corsub}

\begin{proof} By Lemma \ref{lem:trans preserves perimeter}, the edges
of the hexagon tiling are parallel to the edges of a hexagon
tiling satisfying the hypotheses of Lemma
\ref{lem:cong_hexes_on_hex_torus}. By Lemma
\ref{lem:cong_hexes_on_hex_torus}, the edges of the new tiling are
parallel to the three short directions of the torus. Therefore, so
are the edges of the original tiling. \end{proof}

\begin{corsub} \label{cor:standard_hex_has_p_three} Any standard hexagon tiling has perimeter three. \end{corsub}

\begin{proof} By Lemma \ref{lem:trans preserves perimeter}, a standard hexagon tiling has the same perimeter as a partition of the torus satisfying the conditions of Lemma \ref{lem:cong_hexes_on_hex_torus}. By Lemma \ref{lem:cong_hexes_on_hex_torus}, this partition has perimeter three, so the original standard hexagon tiling has perimeter three as well.
\end{proof}

\begin{propsub} \label{prop:only hex is standard hex}
If a hexagon tiling is perimeter minimizing, then it is the
standard hexagon tiling.
\end{propsub}

\begin{proof}
By Lemma \ref{lem:multiple of three}, the number of components in
a minimizing hexagon plus exterior tiling is divisible by three.
By Lemma \ref{lem:6 or more hexes is bad}, a hexagon tiling with
six or more components is not perimeter minimizing. Thus, a
perimeter-minimizing hexagon tiling must have exactly three
components.

By Lemma \ref{lem:trans preserves perimeter}, a
perimeter-minimizing three-component hexagon tiling has the same
perimeter as a hexagon tiling in which the area of each component
is one-third the area of the torus.  By Lemma \ref{lem:small
hexagons lose}, this perimeter is greater than three except on the
hexagonal torus. Hence, by Proposition \ref{prop:perimeter_bound},
any perimeter-minimizing three-component hexagon tiling must lie on the hexagonal
torus. By Corollary
\ref{cor:three_hexes_are_standard}, any such tiling is the
standard hexagon tiling.
\end{proof}

%%%%%%%%%%%%%%%%%%%%%%%%%%%%%%%%%%%%%%%%%%%%%%%%%%%%%%%%
\begin{figure}[htbp]
\begin{center}
\includegraphics*[height = 1.5in]{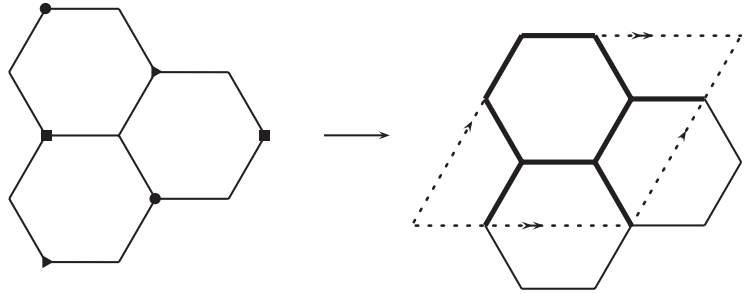}
\end{center}
\caption{\label{fig:hex_cut} This tiling of three congruent
hexagons fits on the hexagonal torus and is one of many perimeter-minimizing standard hexagon tilings.}
\end{figure}
%%%%%%%%%%%%%%%%%%%%%%%%%%%%%%%%%%%%%%%%%%%%%%%%%%%%%%%%

\begin{figure}[ht]
\begin{center}
\includegraphics*{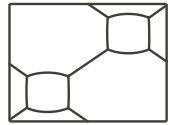}%
\caption{\label{fig:oct-square} An octagon-square tiling always
loses to some other double bubble.}
\end{center}
\end{figure}

\paragraph{\textbf{Octagon-Square Tilings}}
Recall that by Definition \ref{defn:octagon-square}, an
octagon-square tiling is a tiling in which one of the three
regions consists of two curvilinear quadrilaterals, and each of
the other two regions consists of one curvilinear octagon (see
Figure \ref{fig:oct-square}).  The proof that an octagon-square
tiling is not perimeter minimizing generalizes an argument
described to us by Gary Lawlor eliminating most octagon-square
tilings on rectangular tori. Lemma
\ref{lem:band-lens-perimeter-bound} and Corollary
\ref{cor:interface two perimeter bound} compare a purported minimizing
octagon-square tiling to the band lens in order to establish an
upper bound on the length of the interface between the two
eight-sided regions. Then, Proposition \ref{prop:taxicab average}
and Lemma \ref{lem:overlapping bands} establish a lower bound on
the length of the same interface by examining projections of the
boundary curves of the tiling onto four different directions.
Proposition \ref{prop:oct-square bad} shows that these upper and
lower bounds cannot be satisfied at the same time, implying that
no octagon-square tiling is perimeter minimizing.

%%%%%%%%%%%%%%%%%%%%%%%%%%%%%%%%%%%%%%%%%%%%%%%%%%%%%%%%%%%%
%%%%%%%%%%%%%%%%%%%%%%%%%%%%%%%%%%%%%%%%%%%%%%%%%%%%%%%%%%%%

\begin{lemsub}
\label{lem:band-lens-perimeter-bound} Let $A_1 \le A_2 \le A_0$ be
three areas whose sum is equal to the area of the torus. Consider
a lens enclosing area $A_1$ and with perimeter $P$ and diameter
$D$. For any pair of areas from $A_0, A_1, A_2$, a minimizing
double bubble enclosing those areas has perimeter less than or
equal to
\[
2 + P - D = 2 + \sqrt{A_1 \left(\frac{8\pi}{3}-2\sqrt{3}\right)},
\]
with equality only if there is a band lens that encloses area
$A_1$ in the lens and $A_2$ in the band.
\end{lemsub}

\begin{proof}
The expression for $P - D$ is easily derived, and we do not
include the proof here (see, for instance, the derivations in
Proposition \ref{prop:area_and_perim_for band_lens}).

Since all of the double bubbles that enclose two of the areas
from $\{A_0, A_1, A_2\}$ have the same perimeter, it suffices to
prove the claim for a minimizing double bubble enclosing the areas
$A_1, A_2$.

If $D \ge 1$, then $2 + P - D > 2 + 2D - D \ge 3$. Thus, by the
perimeter bound (Proposition \ref{prop:perimeter_bound}), $2 + P - D$
exceeds the perimeter of a minimizing double bubble enclosing
areas $A_1, A_2$.

Suppose instead that $D < 1$. Then the lens can be embedded on the
torus so that its axis lies along a shortest closed geodesic
$\gamma_1$ (Remark \ref{rem:lens-fits}). There exists another
shortest closed geodesic $\gamma_2$ such that $\gamma_1$ and
$\gamma_2$ divide the complement of the lens into one part with
area $A_0$ and another part with area $A_2$. The portions of these
geodesics that lie outside the lens, along with the boundary of
the lens, thus form the boundary of a double bubble enclosing
areas $A_1, A_2$. The perimeter of the double bubble is at most
the total perimeter of $\gamma_2$, the lens, and the portion of
$\gamma_1$ lying outside the lens: $1 + P + (1-D) = 2 + P - D$.
Equality holds only if $\gamma_2$ does not intersect the interior
of the lens --- i.e., only if the double bubble is a band lens.
Hence we have the desired result.
\end{proof}

\begin{corsub} \label{cor:interface two perimeter bound}
Let $A_1 \le A_2 \le A_0$ be three areas whose sum is equal to the
area of the torus. A minimizing double bubble for any pair of
areas from $A_0, A_1, A_2$ has perimeter less than
\[
2 + 2\sqrt{\pi A_1}.
\]
\end{corsub}

\begin{remsub} Rather than examining partitions of the torus into a lens and two
bands as in Lemma \ref{lem:band-lens-perimeter-bound}, one could
prove the corollary directly by examining double bubbles for which
the region of area $A_1$ is a circle and the other two regions are
bands. \end{remsub}

\begin{lemsub} \label{lem:sum of closed geodesics}
If two closed geodesics with lengths $L_1, L_2$ meet at an angle $\phi \in (0, \pi/2]$, then $L_1 + L_2 \ge 1 + 2\cos\phi$. If $L_1, L_2 \in [1, 3/2]$, then $\phi \ge \cos^{-1}(7/9)$. \end{lemsub}

\begin{proof}
The torus can be naturally viewed as the plane with various points
identified. In the plane, choose a point $P$ that represents the intersection of the two closed geodesics. Draw a segment
$\overline{PQ}$ of length $L_1$ representing the corresponding closed geodesic (on the torus). Similarly, draw a segment $\overline{PR}$ of length $L_2$
representing the other closed geodesic, such that segments
$\overline{PQ}$ and $\overline{PR}$ meet at an angle $\phi$.

Because $Q$ and $R$ represent the same point on the torus,
$\overline{QR}$ represents a closed geodesic with length at least
1. Thus, applying the Law of Cosines to triangle $PQR$,
\[
1 \le QR^2 = L_1^2 + L_2^2 - 2 L_1L_2\cos\phi = (L_1+L_2)^2 - (2+2\cos\phi)L_1L_2.
\]
From $0 \le (L_1-1)(L_2-1)$, we derive $L_1 + L_2 - 1 \le L_1L_2$, so that the above inequality becomes
\[
1 \le (L_1+L_2)^2 - (2+2\cos\phi)(L_1+L_2-1),
\]
or
\[
0 \le (L_1+L_2-1-2\cos\phi)(L_1+L_2-1).
\]
Because $L_1 + L_2$ is greater than 1, it must be greater than or equal to $1+2\cos\phi$, as desired.

Suppose now that $L_1, L_2 \in [1, 3/2]$. For a fixed $\phi$, $L_1^2 + L_2^2 - 2L_1L_2\cos\phi$ is convex in each variable $L_1, L_2$, so it is maximized for some $(L_1, L_2)$ with $L_1, L_2 \in \{1, 3/2\}$. For the four pairs $(L_1,L_2)$ satisfying this condition, $L_1^2 + L_2^2 - 2L_1L_2\cos\phi < 1$ when $\cos \phi > 7/9$. Therefore, $\phi \ge \cos^{-1}(7/9)$.
\end{proof}

\begin{remsub}
The next result, Proposition \ref{prop:taxicab average},
generalizes a result about ``taxicab distances'' in the plane. The
taxicab distance between any two points, defined with respect to
two given orthogonal directions, is the minimum length of a
piecewise linear path between them that travels only along the
given orthogonal directions (analogous to the streets a taxicab
travels on). The length of a segment between two points is at
least $\frac{1}{\sqrt{2}}$ times the taxicab distance between
them. (Setting $\phi = \pi/2$ in Proposition \ref{prop:taxicab average}
gives this estimate.) Given two lines that are \textit{not} orthogonal (to take an extreme example, lines parallel to two short directions on the hexagonal
torus), it is useful to have an analogous estimate. We derive such
an estimate relating length to the sum of two taxicab distances,
one corresponding to each of the two given lines.
\end{remsub}

\begin{propsub} \label{prop:taxicab average} Given two pairs of orthogonal lines in the plane, suppose that each line in the first pair meets each line in the second pair at an angle $\phi$ or $\pi/2-\phi$, with $\phi \in [\pi/4,\pi/2]$. For any segment $\gamma$, let $T(\gamma)$ be the sum of the lengths of the projections of $\gamma$ onto the four lines. Then
\[
T(\gamma) \le 2(\sin(\phi/2)+\cos(\phi/2)) \emph{\,\textrm{Length}}(\gamma).
\]
\end{propsub}

\begin{proof}
Without loss of generality, assume that one pair of lines are the
coordinate axes, and that the other lines lie clockwise angles $\phi$ from the positive $x$- and $y$-axes.
Without loss of generality, further assume that $\gamma$ has one
endpoint at the origin and points into the first quadrant of the
plane.

%%%%%%%%%%%%%%%%%%%%%%%%%%%%%%%%%%%%%%%%%%%%%%%%%%%%%%%%%%%%
\begin{figure}[ht]
\begin{center}
\includegraphics*{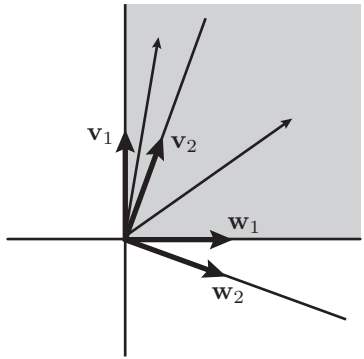}%
\caption{\label{fig:projection vector} If a vector in the first
quadrant has fixed magnitude, then the sum of its projections onto
the four directions $\mathbf{v}_i$, $\mathbf{w}_i$ drawn above attains a local maximum in two places: when it bisects the angle between $\mathbf{v}_1$ and $\mathbf{v}_2$, and when it bisects the angle between $\mathbf{w}_1$ and $\mathbf{v}_2$.}
\end{center}
\end{figure}
%%%%%%%%%%%%%%%%%%%%%%%%%%%%%%%%%%%%%%%%%%%%%%%%%%%%%%%%%%%%

Let $\mbf{u}$ be a vector parallel to $\gamma$, whose magnitude
equals the length of $\gamma$. Let $\mbf{v}_1$, $\mbf{w}_1$,
$\mbf{v}_2$, $\mbf{w}_2$ be unit vectors pointing along the four
given lines, chosen as in Figure \ref{fig:projection vector}. Then
\[
T(\gamma) = \mbf{u} \cdot (\mbf{v}_1 + \mbf{w}_1 + \mbf{v}_2) + |\mbf{u} \cdot \mbf{w}_2|.
\]

If $\mbf{u} \cdot \mbf{w}_2$ is negative, then
\[
T(\gamma) = \mbf{u} \cdot (\mbf{v}_1 + \mbf{w}_1 + \mbf{v}_2 -
\mbf{w}_2) \le \|\mbf{u}\| \|\mbf{v}_1 + \mbf{w}_1 + \mbf{v}_2 -
\mbf{w}_2\|,
\]
with equality when $\mbf{u}$ is parallel to $\mbf{v}_1 + \mbf{w}_1
+ \mbf{v}_2 - \mbf{w}_2$. When this is true, $\mbf{u}$ bisects the
angle formed by $\mbf{v}_1$ and $\mbf{v}_2$, and it is easy to
compute $\mbf{u} \cdot \mbf{v}_1 = \mbf{u} \cdot \mbf{v}_2 =
\|\mbf{u}\|\cos(\phi/2)$ and $\mbf{u} \cdot \mbf{w}_1 = \mbf{u}
\cdot \mbf{w}_2 = \|\mbf{u}\|\sin(\phi/2)$. Hence, when $\mbf{u}
\cdot \mbf{w}_2$ is negative,
\[
T(\gamma) \le \|\mbf{u}\| (2\cos(\phi/2) + 2\sin(\phi/2)) = 2(\cos(\phi/2) + \sin(\phi/2))\textrm{Length}(\gamma).
\]

Similarly, if $\mbf{u} \cdot \mbf{w}_2$ is nonnegative, then
$T(\gamma)/\|\mbf{u}\|$ is maximized when $\mbf{u}$ bisects the
angle formed by $\mbf{w}_1$ and $\mbf{v}_2$. We then find that
\begin{eqnarray*}
T(\gamma) &\le& 2(\cos(\pi/4-\phi/2) + \sin(\pi/4-\phi/2))\textrm{Length}(\gamma) \\
&\le& 2(\cos(\phi/2) + \sin(\phi/2))\textrm{Length}(\gamma),
\end{eqnarray*}
where the second inequality is true because $\phi \in [\pi/4,\pi/2]$. \end{proof}

\begin{lemsub} \label{lem:oct-square-inequality}
If $\cos^{-1}(7/9) \le \phi \le \pi/2$, then
\[
\frac{\max\{2, 1+2\cos\phi\}(1 + \sin\phi)-3}{
\sin\left(\frac{\max\{\phi,\pi/2-\phi\}}{2}\right) +
\cos\left(\frac{\max\{\phi,\pi/2-\phi\}}{2}\right) - 1} \ge 2.
\]
\end{lemsub}

\begin{proof} The denominator of the left hand side is at most $\sqrt{2} -
1$. For $\cos^{-1}(7/9) \le \phi \le \pi/4$, we have $\cos \phi
\ge \sqrt{2}/2$ and $\sin\phi > 2-\sqrt{2}$. Hence, the numerator
of the left hand side is greater than $2\sqrt{2}-2$. Therefore,
the left hand side is greater than 2 when $\phi \le \pi/4$.

For $\pi/4 \le \phi \le \pi/3$, $\max\{2, 1 + 2\cos\phi\} =
1+2\cos\phi$ and $\max\{\phi,\pi/2-\phi\} = \phi$. It is easy to
check that on the interval $\phi \in [\pi/4,\pi/3]$, the numerator
of the left hand side is decreasing while the denominator is
increasing. Because the inequality holds when $\phi = \pi/3$ (with
equality), it holds along the entire interval.

For $\pi/3 \le \phi \le \pi/2$, $\max\{2, 1 + 2\cos\phi\} = 2$ and
$\max\{\phi,\pi/2-\phi\} = \phi$. Clearing the denominator of the
desired inequality, then simplifying, yields the following
inequalities (each equivalent to the desired inequality)
\begin{eqnarray*}
2\sin\phi - 1 &\ge& 2\sin(\phi/2) + 2\cos(\phi/2) - 2 \\
4\sin(\phi/2)\cos(\phi/2) - 1 &\ge& 2\sin(\phi/2) + 2\cos(\phi/2) - 2 \\
(2\sin(\phi/2)-1)(2\cos(\phi/2)-1) &\ge& 0.
\end{eqnarray*}
The last inequality holds for $\pi/3 \le \phi \le \pi/2$.
\end{proof}

\begin{figure}[ht]
\begin{center}
\includegraphics*{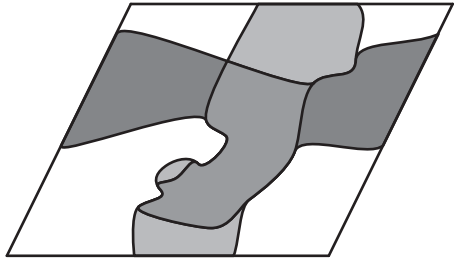}
\end{center}
\caption{\label{fig:overlapping_bands_fig} Given two bands with
overlapping boundary, we can project the boundaries onto various
directions to show that the length of the overlap is large.
(Above, the bands are filled with dark and light shades, with a
medium shade along their intersection.)}
\end{figure}

\begin{lemsub} \label{lem:overlapping bands} Let $\mathcal{B}_1$
and $\mathcal{B}_2$ be the boundaries of two bands with distinct
homologies, where $\mathcal{B}_1$ and $\mathcal{B}_2$ may overlap
each other nontrivially (see Figure
\ref{fig:overlapping_bands_fig}). If the length of $\mathcal{B}_1
\cup \mathcal{B}_2$ does not exceed 3, then the length of
$\mathcal{B}_1 \cap \mathcal{B}_2$ is greater than or equal to 2.
\end{lemsub}

%%%%%%%%%%%%%%%%%%%%%%%%%%%%%%%%%%%%%%%%%%%%%%%%
%%%%%%%%%%%%%%%%%%%%%%%%%%%%%%%%%%%%%%%%%%%%%%%%

\begin{proof} Assume that the length of $\mathcal{B}_1 \cup
\mathcal{B}_2$ is as at most 3, and let $\gamma_1, \dots,
\gamma_k$ be the components of $\mathcal{B}_1 \cap \mathcal{B}_2$.
For $i = 1, 2$, let $L_i$ be the length of the closed geodesic
$g_i$ with the same homology as the components of $\mathcal{B}_i$.
Observe that $L_i \le 3/2$, because otherwise $\mathcal{B}_i$
would have length greater than 3. Let $g_1$ and $g_2$ meet each
other at an angle $\phi \in (0, \pi/2]$.

Below, we will take each closed curve in $\mathcal{B}_1$ and project
it onto two directions --- onto the direction of $g_1$ (Figure
\ref{fig:bandprojections}(a)) and onto the direction orthogonal to
$g_2$ (Figure \ref{fig:bandprojections}(b)). We will also take similar
projections of the components of $\mathcal{B}_2$. By analyzing
these projections' effects on various pieces of $\mathcal{B}_1 \cup \mathcal{B}_2$, we will obtain our desired inequality.

For each $\gamma_i$, let $T(\gamma_i)$ be the sum of the
projection of $\gamma_i$ onto four directions: the two directions
of the closed geodesics $g_1$ and $g_2$, and the two directions
orthogonal to these closed geodesics.

%%%%%%%%%%%%%%%%%%%%%%%%%%%%%%%%%%%%%%%%%%%%%%%%%%%%%%%%%%%%
\begin{figure}[ht]
\begin{center}
\begin{tabular} {@{\extracolsep{1.5 cm}} cc}
\includegraphics*[height=1in]{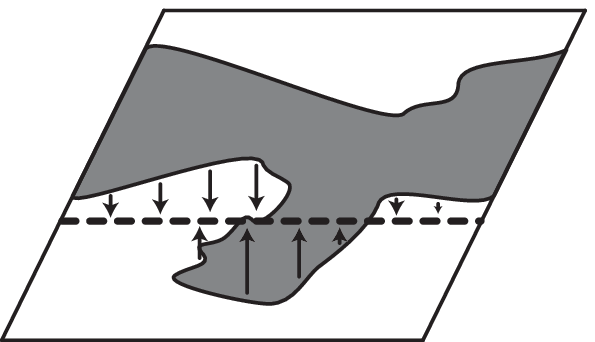}
&
\includegraphics*[height=1in]{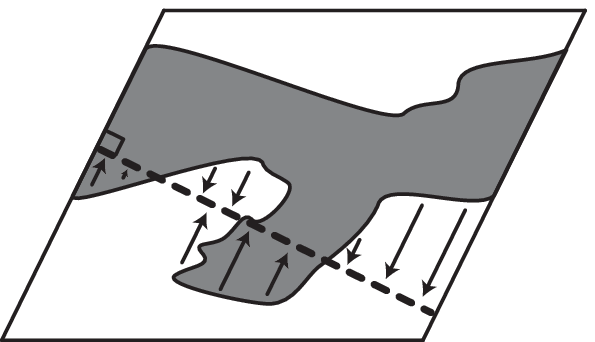}
\\
(a) & (b)
\end{tabular}
\end{center}
\caption{\label{fig:bandprojections} By projecting each component
of the boundary of the bands in two directions, we bound the
length of the intersection of the two boundaries.}
\end{figure}
%%%%%%%%%%%%%%%%%%%%%%%%%%%%%%%%%%%%%%%%%%%%%%%%%%%%%%%%%%%%

Take a component of $\mathcal{B}_1$, and project it onto the
direction of $g_1$ and onto the direction orthogonal to $g_2$. The
first projection has length at least $L_1$, and the second
projection has length at least $L_1\sin\phi$. Project the other
component of $\mathcal{B}_1$ in a similar fashion.

Likewise, project each components of $\mathcal{B}_2$ onto the direction of $g_2$ and onto the direction orthogonal to $g_1$. In all, the projections of the components of $\mathcal{B}_1$ and $\mathcal{B}_2$ have total length greater than or equal to
\[
2(L_1 + L_2)(1 + \sin \phi).
\]

Each component of $\mathcal{B}_1 \cup \mathcal{B}_2 -
\mathcal{B}_1 \cap \mathcal{B}_2$ is projected exactly twice, so
they contribute at most
\[
2 \cdot \textrm{Length}(\mathcal{B}_1 \cup \mathcal{B}_2 - \mathcal{B}_1 \cap \mathcal{B}_2)
\]
to the total length of the projections. Each $\gamma_i$, on the other hand, contributes at most $T(\gamma_i)$. Therefore,
\begin{eqnarray*}
2(L_1 + L_2)(1 + \sin \phi) &\le& 2 \cdot \textrm{Length}(\mathcal{B}_1 \cup \mathcal{B}_2 - \mathcal{B}_1 \cap \mathcal{B}_2) + \sum_{i=1}^{k} T(\gamma_i) \\
&\le& 2\left(3-\sum_{i=1}^{k} \textrm{Length}(\gamma_i)\right) + \sum_{i=1}^{k} T(\gamma_i).
\end{eqnarray*}
We estimate the two sides of this inequality as follows. First,
$L_1 + L_2 \ge 2$, and by Lemma \ref{lem:sum of closed geodesics},
$L_1 + L_2 \ge 1 + 2\cos\phi$. Hence,
\[
2\max\{2,1+2\cos\phi\}(1+\sin\phi) \le 2(L_1 + L_2)(1+\sin\phi).
\]
Second, from the two pairs of orthogonal directions we project
onto, we can choose one direction from each pair so that they meet
at an angle of $\max\{\phi,\pi/2-\phi\} \in [\pi/4,\pi/2]$.
Applying Proposition \ref{prop:taxicab average} then shows that
\[
2\left(3-\sum_{i=1}^{k} \textrm{Length}(\gamma_i)\right) + \sum_{i=1}^{k} T(\gamma_i)
\]
is at most
\[
6 +\ \left(2\sin\left(\frac{\max\{\phi,\pi/2-\phi\}}{2}\right) + 2\cos\left(\frac{\max\{\phi,\pi/2-\phi\}}{2}\right) - 2\right)\sum_{i=1}^{k}\textrm{Length}(\gamma_i).
\]

Therefore, to prove that $\sum_{i=1}^{k} \textrm{Length}(\gamma_i) \ge 2$, it suffices to prove that
\[
\frac{\max\{2, 1+2\cos\phi\}(1 + \sin\phi)-3}{
\sin\left(\frac{\max\{\phi,\pi/2-\phi\}}{2}\right) + \cos\left(\frac{\max\{\phi,\pi/2-\phi\}}{2}\right) - 1} \ge 2.
\]
By Lemma \ref{lem:sum of closed geodesics}, $\cos^{-1}(7/9) \le
\phi$; also, by definition, $\phi \le \pi/2$. By Lemma
\ref{lem:oct-square-inequality}, the above inequality holds for
such $\phi$.
\end{proof}

\begin{propsub} \label{prop:oct-square bad}
No octagon-square tiling is perimeter minimizing.
\end{propsub}

%%%%%%%%%%%%%%%%%%%%%%%%%%%%%%%%%%%%%%%%%%%%%%%%
\begin{figure}[ht]
\begin{center}
\begin{tabular} {@{\extracolsep{1.5 cm}} cc}
\includegraphics*[height=1in]{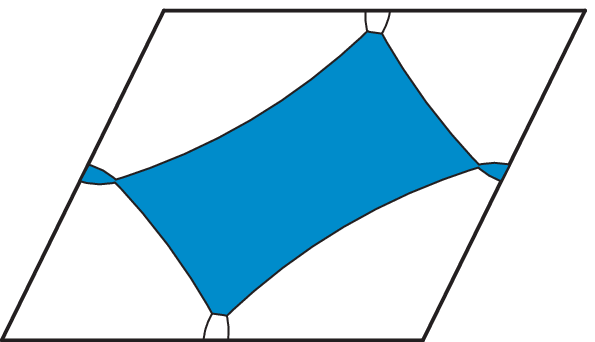}
&
\includegraphics*[height=1in]{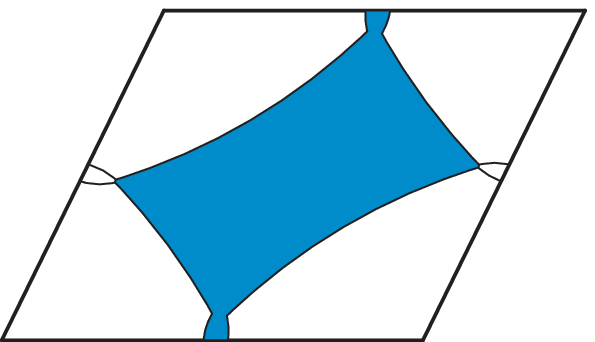}
\\
(a) & (b)
\end{tabular}
\end{center}
\caption{\label{fig:oct-square_bands} Given the above octagon-square tiling, the boundaries of the two colored bands either have a union with length greater than three, or an intersection with length greater than or equal to two. In both
cases, it can be shown that the tiling is not
perimeter minimizing.}
\end{figure}
%%%%%%%%%%%%%%%%%%%%%%%%%%%%%%%%%%%%%%%%%%%%%%%%

\begin{proof} Let $\mathcal{B}_1$ be the boundary of the band enclosing one eight-sided region and one four-sided component (Figure \ref{fig:oct-square_bands}(a)). Let $\mathcal{B}_2$ be the boundary of the band enclosing the same eight-sided region and the other four-sided component (Figure \ref{fig:oct-square_bands}(b)). If the length of $\mathcal{B}_1 \cup \mathcal{B}_2$ is greater than 3, then by the perimeter bound (Proposition \ref{prop:perimeter_bound}), the octagon square is not perimeter minimizing.

Otherwise, by Lemma \ref{lem:overlapping bands}, the length of
$\mathcal{B}_1 \cap \mathcal{B}_2$ is at least 2. Let $A_1$ be the
area of the region containing the four-sided components. By the
isoperimetric inequality for the plane, the perimeter of this
region is at least $2\sqrt{\pi A_1}$. Therefore, the perimeter of
the octagon-square tiling is at least $2 + 2\sqrt{\pi A_1}$. By
Corollary \ref{cor:interface two perimeter bound}, the tiling is
not a minimizer. \end{proof}

\begin{propsub} \label{prop:only tiling is 3 hex}
The only tiling that can be perimeter minimizing is the standard
hexagon tiling.
\end{propsub}

\begin{proof}
By Corollary \ref{cor:hex-or-oct}, the tiling must be a hexagon
tiling or an octagon-square tiling.  By Proposition
\ref{prop:oct-square bad}, no octagon-square tiling is perimeter
minimizing.  By Proposition \ref{prop:only hex is standard hex},
the only minimizing hexagon tiling is the standard hexagon tiling.
\end{proof}

\section{Perimeter-Minimizing Double Bubbles on the Flat Two-Torus}\label{S:main}

The Main Theorem characterizes the perimeter-minimizing double
bubbles on all flat two-tori, as depicted in Figure \ref{fig:introduce_minimizers}.  Corollaries \ref{cor:cylinder
corollary} and \ref{cor:free-boundary} characterize the minimizing
double bubbles on the flat infinite cylinder and the flat infinite strip
with free boundary. Remark \ref{rem:klein} makes a conjecture about the
double bubble problem on the closely related Klein bottle.

\begin{main}\label{thm:main theorem}
A perimeter-minimizing double bubble on a flat two-torus is of one
of the following five types (perhaps after relabelling the three regions), depending on the areas to be enclosed: \\[-8pt]
\begin{quote}
(i) the standard double bubble, \\[6pt]
(ii) a standard chain, \\[6pt]
(iii) the band lens, \\[6pt]
(iv) the double band, or \\[6pt]
(v) the standard hexagon tiling (only for the hexagonal torus).
\end{quote}

\end{main}

\begin{proof}
Proposition \ref{prop:topo classification} shows that a minimizer
must be a contractible double bubble, a tiling, a swath, a band
adjacent to a set of components whose union is contractible, or
the double band (after perhaps relabelling the two regions and the
exterior).

By Proposition \ref{prop:sdb_is_minimizer}, the standard double
bubble is the only possible contractible minimizer.  By
Proposition \ref{prop:chain_is_minimizer}, a standard chain is the
only possible minimizer among swaths.  By Proposition
\ref{prop:band lens_is_minimizer}, the band lens is the only
possible minimizer among bands adjacent to a set of components
whose union is contractible.  By Proposition \ref{prop:only tiling
is 3 hex}, the standard hexagon tiling (on the hexagonal torus) is
the only possible minimizing tiling.
\end{proof}

\begin{rem}
Using the phase portrait computations of Section \ref{S: numerical
comparisons}, it can be confirmed that all five types of
minimizers do occur.
\end{rem}

\begin{cor} \label{cor:cylinder corollary}
A perimeter-minimizing double bubble on the flat infinite cylinder
is of one of the following four types, depending on the areas to
be enclosed: \\[-8pt]
\begin{quote}
(i) the standard double bubble, \\[6pt]
(ii) a standard chain, \\[6pt]
(iii) the band lens, or \\[6pt]
(iv) the double band.
\end{quote}
\end{cor}

\begin{proof}
Without loss of generality, assume that the cylinder has
circumference 1. Let $\Sigma$ be a perimeter-minimizing double
bubble on the cylinder that is not of one of the four types listed
above. Because $\Sigma$ has finite perimeter, it can be placed on
some long rectangular torus (for example, one with length greater than 3).
Then by Theorem \ref{thm:main theorem}, there exists a double bubble $\Sigma'$ on the torus
enclosing the same areas as $\Sigma$ but with less perimeter. Because the torus has length greater than 3, and because $\Sigma'$ has perimeter less
than or equal to 3, it is contractible or wraps around a short
direction of the torus. Therefore, $\Sigma'$ can be placed on the
cylinder, and $\Sigma$ cannot be perimeter minimizing.
\end{proof}

\begin{cor}\label{cor:jass}
Every minimizing double bubble on the infinite cylinder has a line of reflective symmetry that is perpendicular to a shortest closed geodesic, or is congruent to a double bubble with this property.
\end{cor}

\begin{proof} Given a minimizing standard double bubble, it can be rotated so that its axis of symmetry is perpendicular to a shortest closed geodesic. Each of the other minimizers in Corollary \ref{cor:cylinder corollary} has a line of reflective symmetry with the required property. (A standard chain on the infinite cylinder has a line of reflective symmetry that is parallel to a shortest closed geodesic, but it also has two lines of reflective symmetry perpendicular to a shortest closed geodesic.) \end{proof}

\begin{figure}[htbp]
\begin{center}
\includegraphics*{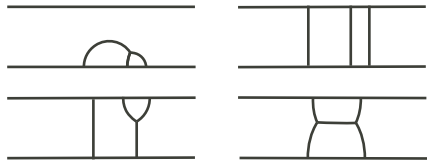}%
\end{center}
\caption{\label{fig:strip-figs} The four perimeter-minimizing double bubbles on the flat infinite strip with free boundary are each half of some perimeter-minimizing double bubble on the infinite cylinder.}
\end{figure}

\begin{cor} \label{cor:free-boundary}
Any minimizing double bubble on the flat infinite strip with free boundary is the quotient of a minimizing double bubble on an infinite cylinder (with twice the width of the strip) by a reflection across one of its lines of symmetry, using the free boundary of the strip to complete the double bubble, as in Figure \ref{fig:strip-figs}.
\end{cor}

\begin{proof}
The infinite strip is the quotient space of the infinite cylinder by some reflection.  Hence we can think of the
infinite cylinder as a double cover of the infinite strip. Lifting a
minimizer $\Sigma$ on the strip yields a double bubble $\Sigma'$
on the cylinder with twice the perimeter and enclosed areas as $\Sigma$.  If
$\Sigma'$ is a minimizing double bubble, then we are done.
Otherwise, by Corollary \ref{cor:jass}, there exists a minimizer $\Gamma'$ that encloses the same areas as $\Sigma'$ and that has at least one line of reflective symmetry perpendicular to a shortest closed geodesic. Quotienting by the corresponding reflection
yields the original infinite strip with free boundary, along with a new double bubble
$\Gamma$ enclosing the same areas as $\Sigma$ with less perimeter,
a contradiction.
\end{proof}

\begin{rem}
For finite rectangular cylinders with free boundary (i.e., the
quotient space of a rectangular torus by some reflection), a
result analogous to Corollary \ref{cor:free-boundary} holds if the
distance between the two components of the boundary is
sufficiently small with respect to the length of each component of
the boundary.
\end{rem}

%%%%%%%%%%%%%%%%%%%%%%%%%%%%%%%%%%%%%%%%%%%%%%%%%%%%%%%%
%%%%%%%%%%%%%%%%%%%%%%%%%%%%%%%%%%%%%%%%%%%%%%%%%%%%%%%%
\begin{rem} {\rm (Klein Bottles)}\label{rem:klein}
We conjecture that perimeter-minimizing double bubbles on the flat
Klein bottle are of virtually the same form as on the flat torus.
Like the flat torus, we can represent the flat Klein bottle as a
planar parallelogram with opposite sides identified and no angle
less than 60 degrees, except that we identify the top and bottom
horizontal sides with a flip. Because of this flip, on a Klein bottle with
reasonably long top and bottom, any minimizing double band is
bounded by two closed geodesics each of length 2 and with total
perimeter 4 (instead of 3). This is an important difference, and
for proofs on the torus where we required a perimeter bound of 3,
it will be necessary to modify our approach to accommodate the
higher bound for the Klein bottle. Note also that the placement of
perimeter-minimizing double bubbles on a wide Klein bottle is less
flexible than on a torus. The standard chain, band lens, and
double band must all be oriented symmetrically about one of the
two unique shortest closed geodesics of a wide Klein bottle (which
are always vertical). However, all of the conjectured minimizers
(including the standard hexagon tiling) do remain well-defined, so
long as they are placed carefully.

We also conjecture that the perimeter-minimizing double bubbles on
the infinite cylinder (which can be thought of as an infinitely long rectangular torus) carry over in the analogous space corresponding to the Klein bottle, that is, a M\"{o}bius strip of infinite length --- the surface between two horizontal lines in the plane, where the lines are identified with a flip.  Here the same comments apply, but the absence of tilings may make the full solution more accessible.
\end{rem}

\section{Formulas for Perimeter and Area}\label{S:perimeter-and-area}

This section gives formulas for the areas of the enclosed regions and
the perimeter of the four potential minimizers different from the
double band: the \sdb\ (Proposition \ref{prop:area_and_perim_for sdb}), the standard chain (Propositions \ref{prop:area_and_perim_for_sym_chains},
\ref{prop:perim_from_area_for_sym_chains}), the \bl\ (Proposition
\ref{prop:area_and_perim_for band_lens}), and the standard hexagon
tiling (Proposition \ref{prop:hex_three_parameterization}). In addition, Lemmas \ref{lem:symm chain fits} and \ref{lem:hex-tiling-wins-if-we-think-so} prove that when a standard chain with axis-length one or a standard hexagon tiling \textit{seems} to be minimizing, based on formal perimeter calculations that we can perform even if the double bubble actually does not exist on a given torus, then the double bubble actually exists and is minimizing. As described in Section \ref{S: numerical comparisons}, we use these formulas and results to produce the phase portraits in Figure \ref{fig:param-diags}.

\begin{rem}
We begin by stating a few simple geometric formulas that describe
circular arcs meeting chords. Given a circular arc subtended by a chord of length $C$ that it meets at an angle
$\theta$, the area between the arc and the chord, the length of
the arc, and the radius of the circle are given by:
\begin{eqnarray}
A(\theta, C) &=& \frac{C^2 (\theta- \sin \theta \cos \theta)}{4
\sin^2{\theta}} \label{eqn:A} \\
L(\theta, C) &=& \frac{C \theta}{\sin{\theta}} \label{eqn:L}\\
R(\theta, C) &=& \frac{C}{2 \sin \theta} \label{eqn:R}.
\end{eqnarray}\\[-12pt]
\end{rem}

\subsection{Formulas for the Standard Double Bubble} \label{SS:Param_SDB}

\par

Proposition \ref{prop:area_and_perim_for sdb} gives formulas for
perimeter and area of the standard double bubble.  These formulas
first appeared in [F, Section 2].

\begin{propsub}
\label{prop:area_and_perim_for sdb} Let $R_1$ and $R_2$ be the
regions of higher and lower pressure enclosed in a standard double
bubble. The area and perimeter formulas of a standard double
bubble in terms of the separating chord length $C$ and the angle
$\theta$ at which it meets the interior arc are
\begin{eqnarray*}
A_{R_1}(\theta, C) &=& A(2 \pi/3 - \theta, C) + A(\theta, C) \\[6pt]
A_{R_2}(\theta, C) &=& A(2 \pi/3 + \theta, C) - A(\theta, C) \\[6pt]
P(\theta, C) &=& L(2 \pi/3 + \theta, C) + L(2 \pi/3 - \theta, C)+ L(\theta, C).
\end{eqnarray*}
\end{propsub}

\begin{proof}
The arc between $R_1$ and the exterior meets the separating chord
$C$ at $\frac{2\pi}{3} - \theta$, and the arc between $R_2$ from
the exterior meets the separating chord $C$ at $\frac{2\pi}{3} +
\theta$. The given formulas follow immediately.
\end{proof}

\subsection{Formulas for the Standard Chain}\label{SS:Param_2_Comp_Sym_Chains}

To parameterize the standard chain, we distinguish between those
in which the enclosed regions have unequal pressure (Proposition
\ref{prop:area_and_perim_for_sym_chains}) and those in which the
enclosed regions have equal pressure (Proposition
\ref{prop:perim_from_area_for_sym_chains}).

%%%%%%%%%%%%%%%%%%%%%%%%%%%%%%%%%%%%%%%%%%%%%%%%%%%%%%%%%%%%%%%%%
\begin{figure}[htbp]
\includegraphics*{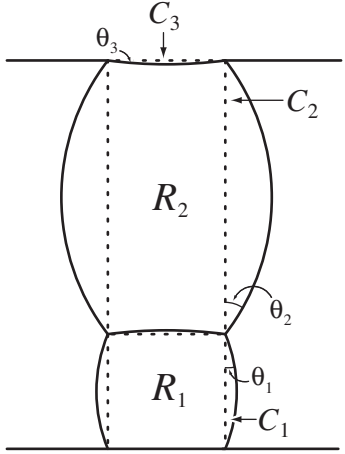}
\caption{\label{fig:fig_for_param_sym_chains} This picture
illustrates the basic parameters involved in computing perimeter
and area for a standard chain.}
\end{figure}
%%%%%%%%%%%%%%%%%%%%%%%%%%%%%%%%%%%%%%%%%%%%%%%%%%%%%%%%%%%%%%%%%

\begin{propsub}
\label{prop:area_and_perim_for_sym_chains} Consider a standard
chain with axis-length $L_0$ and with two components of unequal
pressure. Let $R_1$ and $R_2$ be the enclosed regions, and suppose
that $R_1$ has greater pressure than $R_2$. Let $C_1$ be the
length of either chord subtending an arc separating $R_1$ from the
exterior, where the chord and arc meet at (interior) angles
$\theta_1$. (See Figure \ref{fig:fig_for_param_sym_chains}.) The
area and perimeter formulas are
\begin{eqnarray*}
A_{R_1}(\theta_1, C_1) &=& 2 A(\theta_1, C_1) + 2 A(\pi/6 - \theta_1, C_3) + C_1 C_3, \\[6pt]
A_{R_2}(\theta_1, C_1) &=& 2 A(\pi/3 - \theta_1, C_2) - 2A(\pi/6 - \theta_1 , C_3) + C_2 C_3, \\[6pt]
P(\theta_1, C_1) &=& 2 L(\theta_1, C_1) + 2 L(\pi/3 - \theta_1,
C_2) + 2 L(\pi/6 - \theta_1, C_3),
\end{eqnarray*}
where $C_2$ and $C_3$ are given by
\begin{eqnarray*}
C_2 &=& L_0 - C_1, \\
C_3 &=& \frac{C_1 (L_0-C_1) \sin (\thirty -
\theta_1)}{(L_0-C_1)\sin \theta_1 - C_1 \sin (\sixty - \theta_1)}.
\end{eqnarray*}
A standard chain with these parameters exists on an infinite
cylinder of circumference $L_0$ if and only if $C_1 \in (0, \half
L_0)$ and $\theta_1 \in (\sin^{-1}(C_1), \pi/6)$.
\end{propsub}

\begin{proof}
Let $R_0$ be the low-pressure exterior region. By the results of
Section \ref{SS:swaths}, the vertices of $R_1$ and $R_2$ form a
rectangle.  For $i = 1, 2$, let either arc separating $R_i$ from
$R_0$ have curvature $\kappa_{i}$ and radius of curvature $r_{i} =
1/\kappa_{i}$, let $C_{i}$ denote the length of the chord with the
same endpoints as the arc, and let $\theta_{i}$ be the (interior)
angle the chord and the arc make. Define $\kappa_3$, $r_3$, $C_3$,
and $\theta_3$ similarly with respect to the interface between
$R_1$ and $R_2$.

Because the arcs in the boundary of the chain meet at $2\pi/3$, we
have
\[
\theta_{2} = \pi/3 - \theta_{1} \qquad \textrm{and} \qquad
\theta_{3} = \pi/6 - \theta_{1}.
\]
Of course,
\[
C_{2} = L_0 - C_{1}.
\]
Also, by the Extended Law of Sines,
\[
C_{i} = 2r_{i}\sin\theta_{i} \qquad \textrm{for\ }i = 1, 2, 3,
\]
allowing us to solve for $r_{1}$ and $r_{2}$ in terms of
$\theta_{1}$ and $C_{1}$. By the cocycle condition (as in
Proposition \ref{prop:regularity}),
\[
- \frac{1}{r_{1}} + \frac{1}{r_{2}} + \frac{1}{r_{3}} = - \kappa_{1} + \kappa_{2} + \kappa_{3} = 0,
\]
and we may solve for $r_3$ in terms of $r_1$ and $r_2$. These
equations allow us to find $C_{3}$ in terms of $\theta_{1}$ and
$C_{1}$:
\begin{eqnarray*}
C_{3} &=& 2\frac{1}{1/r_1-1/r_2}\sin\theta_3 = \frac{2r_1r_2}{r_2-r_1}\sin\theta_3 \\
&=& \frac{C_1C_2}{C_2\sin\theta_1 - C_1\sin\theta_2}\sin\theta_3 = \frac{C_{1}(L_0-C_{1})\sin(\pi/6-\theta_{1})}{(L_0-C_{1})\sin\theta_{1} - C_{1}\sin(\pi/3-\theta_{1})}.
\end{eqnarray*}

It remains to show that $C_1 \in (0, L_0/2)$, $\theta_1 \in
(\sin^{-1}(C_1/L_0), \frac{\pi}{6})$, and that a corresponds
standard chain exists (on some infinite cylinder) if these bounds
hold.

Because $r_2 > r_1$ and $\theta_2 > \theta_1$, we have $L_0 - C_1
= C_2 = 2r_2\sin\theta_2
> 2r_1\sin\theta_1 = C_1$, proving the bound on $C_1$.
The upper bound on $\theta_1$ holds because $\pi/6 - \theta_1 =
\theta_3 > 0$.

%%%%%%%%%%%%%%%%%%%%%%%%%%%%%%%%%%%%%%%%%%%%%%%%%%%%%%%%%%%%%%%%%
\begin{figure}[htbp]
\includegraphics*[bb=185 675 414 722, height=.75in]{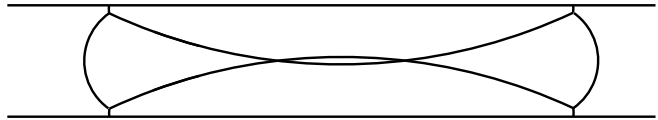}
\caption{\label{fig:self_int_chain} If the angles and lengths in the parameterization of a standard chain do not meet certain bounds, then the boundary arcs of the chain intersect each other.}
\end{figure}
%%%%%%%%%%%%%%%%%%%%%%%%%%%%%%%%%%%%%%%%%%%%%%%%%%%%%%%%%%%%%%%%%
%%%%%%%%%%%%%%%%%%%%%%%%%%%%%%%%%%%%%%%%%%%%%%%%%%%%%%%%%%%%%%%%%
\begin{figure}[htbp]
\begin{center}
\includegraphics*[height=1.1in]{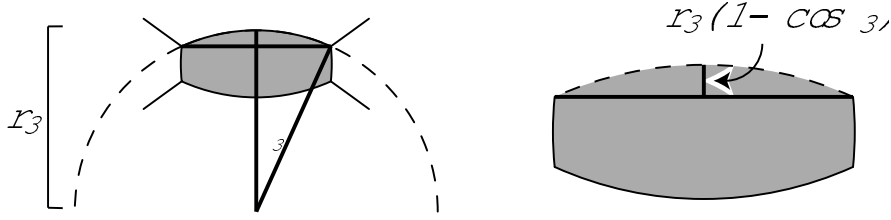}
\end{center}
\caption{\label{fig:selftangent_standard_chain}  The arcs
separating the components enclosed in a standard chain cannot
intersect each other, yielding bounds on the parameters used to
describe standard chains.}
\end{figure}
%%%%%%%%%%%%%%%%%%%%%%%%%%%%%%%%%%%%%%%%%%%%%%%%%%%%%%%%%%%%%%%%%

Given $C_1 \in (0,L_0/2)$ and $0 < \theta_1 < \pi/6$, a
corresponding standard chain exists (on some infinite cylinder) if
and only if it can be drawn so that the arcs separating $R_1$ and
$R_2$ do not intersect. Consider one of these arcs. As shown in
Figure \ref{fig:selftangent_standard_chain}, the distance between
the midpoint of the arc and the chord subtending the arc is $r_3 -
r_3\cos\theta_3$. The two arcs separating $R_1$ and $R_2$ do not
intersect each other if and only if this distance is less than
$\half C_2$ --- that is, if and only if
\[
r_3(1-\cos\theta_3) < \half C_2.
\]
(This inequality fails if and only if the arcs intersect each
other as in Figure \ref{fig:self_int_chain}.) Writing $r_3$,
$\theta_3$, and $C_2$ in terms of $C_1$ and $\theta_1$, this
inequality becomes (after much simplification)
\[
C_1 < L_0 \sin \theta_1,
\]
or $\theta_1 > \sin^{-1}(C_1/L_0)$.
\end{proof}

\begin{propsub}
\label{prop:perim_from_area_for_sym_chains} Consider a standard
chain with axis-length $L_0$ and with two regions $R_1$ and $R_2$
of equal pressure. The area and perimeter formulas are
\begin{eqnarray*}
A_{R_1}(C_3) = A_{R_2}(C_3) &=& \frac{2\pi-3\sqrt{3}}{24}L_0^2 + \frac{1}{2} L_0C_3 \\[6pt]
P(C_3) &=& \frac{2\pi}{3}L_0 + 2C_3,
\end{eqnarray*}
where $C_3$ is the length of the chord separating $R_1$ from
$R_2$.
\end{propsub}

\begin{proof}
For $i = 1, 2, 3$, define $\kappa_i$, $r_i$, $C_i$, and $\theta_i$
as in the first paragraph of the proof of Proposition
\ref{prop:area_and_perim_for_sym_chains}. (Figure
\ref{fig:fig_for_param_sym_chains} shows many of these variables,
although the standard chain in the figure does not contain two
regions of equal pressure.)

Observe that $\kappa_1 = \kappa_2$ and hence $r_1 = r_2$. Also,
$\theta_3 = 0$ and $\theta_1 = \theta_2 = \pi/6 - \theta_3 =
\pi/6$. Therefore, $C_1 = 2r_1\sin\theta_1 = 2r_2\sin\theta_2 =
C_2$, implying that $C_1 = C_2 = \half L_0$ and $r_1 = r_2 = \half
L_0/(2\sin(\pi/6)) = \half L_0$.

Although $C_1$ and $C_2$ are uniquely determined, $C_3$ is not.
Still, because all other arc lengths and angles in the chain are
fixed, we have simple formulas expressing the areas and perimeters
of such chains in terms of $C_3$. Given the fixed angles and chord
lengths found earlier, we have $A_{R_1}(C_3) = A_{R_2}(C_3) =
2A(\pi/6,\frac{L_0}{2}) + (\frac{1}{2} L_0) C_3$ and $P(C_3) =
4L(\pi/6,\frac{L_0}{2}) + 2C_3$. Using Equations
\ref{eqn:A}-\ref{eqn:L} to evaluate these expressions yields the
desired formulas.
\end{proof}

\begin{remsub} \label{rem:standard chain exists}
The plots based on this parameterization prove the existence of a
family of standard chains with fixed axis-length.  In particular,
we know that there exists a family of standard chains that wrap
around a short direction of the torus.
\end{remsub}

\begin{lemsub}\label{lem:symm chain fits}
Consider any standard chain on the unit-circumference cylinder
enclosing areas $A_1, A_2$, and any torus with area greater than $A_1 + A_2$. The perimeter of the chain is greater than or equal
to that of the minimizing double bubble(s) enclosing areas $A_1,
A_2$ on the torus. If equality holds, the standard chain fits on
the torus and is itself a minimizing double bubble on the torus.
\end{lemsub}

\begin{proof}
Consider a fundamental domain of the torus, a parallelogram with
side lengths 1 and $L$ and with interior angle $\alpha \in [\pi/3,
\pi/2]$. We orient the parallelogram so that the side with
length 1 is vertical. Also, let $R_1$ be a component of higher (or
equal) pressure enclosed in the chain, and let $R_2$ be the other
component enclosed in the chain. For $i = 1, 2, 3$ define $C_i$
and $\theta_i$ as in the first paragraph of the proof of Proposition \ref{prop:area_and_perim_for_sym_chains}.

We immerse the chain in the torus, such that the chain and the side of length 1 have the
same homology. It suffices to prove that if the chain is not embedded, then it has more perimeter than some (embedded) minimizer enclosing areas $A_1, A_2$ on the torus. Because the chain is not embedded in the torus, it does not only wrap around the torus in the direction of the side with length 1; it also
wraps around the torus in the direction of the side with length
$L$.  We consider two cases: first, where the chain wraps around
to intersect itself because the two components intersect each
other; second, where instead at least one component intersects
itself but not the other component.

%%%%%%%%%%%%%%%%%%%%%%%%%%%%%%%%%%%%%%%%%%%%%%%%%%%%%%%%%%%%%%%%%
\begin{figure}
\begin{center}
\includegraphics*[height=2in]{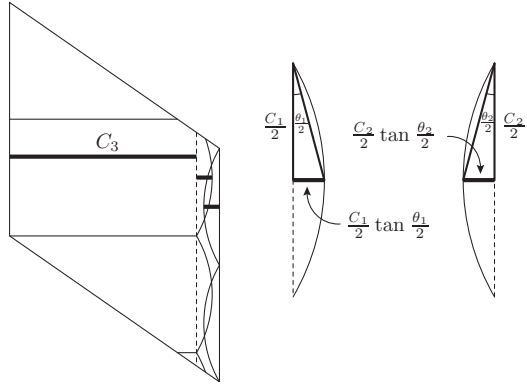}
\end{center}
\caption{\label{fig:nonembedded standard chains1} One way a
standard chain might be nonembedded is if the two components
intersect each other. In this case, by considering three segments
perpendicular to the short direction with total length at least
$\sqrt{3}/2$, it can be shown that the arcs separating the two
components in the chain are long enough (greater than $.5$) to
prove that the chain does not have minimizing perimeter.}
\end{figure}
%%%%%%%%%%%%%%%%%%%%%%%%%%%%%%%%%%%%%%%%%%%%%%%%%%%%%%%%%%%%%%%%%
First suppose that both components of the chain intersect each
other. For $i = 1, 2$, consider either boundary arc of $R_i$ that
lies along the boundary of the chain. The horizontal distance
between the midpoint of this arc and the chord subtending it is
$d_i = (C_i/2) \tan (\theta_i/2)$. Also, the vertices of the
standard chain lie on two vertical closed geodesics, and the
horizontal distance between these closed geodesics is $C_3$. Now,
because the chain wraps around the torus in a nonvertical
direction, the sum of $d_1$, $d_2$, and $C_3$ must exceed the
horizontal distance $L\sin\alpha$ between the two vertical sides
of the parallelogram. Hence,
\[
\frac{C_1}{2} \tan\left(\frac{\theta_1}{2}\right) + \frac{C_2}{2}
\tan\left(\frac{\theta_2}{2}\right) + C_3 > L\sin\alpha \ge
\sin(\pi/3) = \frac{\sqrt{3}}{2}.
\]
From our parameterization, $C_1 + C_2 = 1$ and $\theta_2 <
\theta_1 < \pi/3$. Therefore,
\[
\frac{C_1}{2} \tan\left(\frac{\theta_1}{2}\right) + \frac{C_2}{2}
\tan\left(\frac{\theta_2}{2}\right) < \frac{C_1}{2}
\tan\left(\frac{\pi}{6}\right) + \frac{C_2}{2}
\tan\left(\frac{\pi}{6}\right) =
\frac{1}{2}\tan\left(\frac{\pi}{6}\right) = \frac{\sqrt{3}}{6}.
\]
Combining these two inequalities, we find that
\[
C_3 > \frac{\sqrt{3}}{2} - \frac{\sqrt{3}}{6} = \frac{\sqrt{3}}{3}
> \frac{1}{2}.
\]
Hence, each of the arcs separating the two components enclosed in
the chain has length greater than $\frac{1}{2}$. It easily follows
that the perimeter of the chain is greater than 3. Therefore, by
the perimeter bound (Proposition \ref{prop:perimeter_bound}), the
standard chain has perimeter greater than that of the (embedded)
minimizer(s) enclosing the same areas.

%%%%%%%%%%%%%%%%%%%%%%%%%%%%%%%%%%%%%%%%%%%%%%%%%%%%%%%%%%%%%%%%%
\begin{figure}
\begin{center}
\includegraphics*[height=1in]{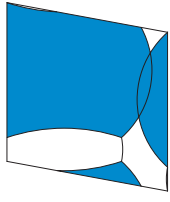}
\end{center}
\caption{\label{fig:nonembedded standard chains2} A second way a
standard chain might be nonembedded is if at least one component
intersects itself. The component then forms a band wrapping around
the torus, implying (as was the case for octagon-square tilings) that the chain has a large perimeter.}
\end{figure}
%%%%%%%%%%%%%%%%%%%%%%%%%%%%%%%%%%%%%%%%%%%%%%%%%%%%%%%%%%%%%%%%%

Next suppose that one of the components --- say, $R_i$ ---
intersects itself but not the other component. The area enclosed
by $R_i$ on the torus thus forms a band wrapping around the torus
in a nonvertical direction; let $\mathcal{B}_1$ be the boundary of
this band. Also, let $\mathcal{B}_2$ be the boundary of a second
band containing the two components of the chain (and excluding any area
that lies in the overlap of the immersed chain). By Lemma
\ref{lem:overlapping bands}, either the length of $\mathcal{B}_1
\cup \mathcal{B}_2$ is greater than 3, or the length of
$\mathcal{B}_1 \cap \mathcal{B}_2$ is greater than or equal to 2.

If the length of $\mathcal{B}_1 \cup \mathcal{B}_2$ is greater
than 3, then by the perimeter bound (Proposition
\ref{prop:perimeter_bound}), the perimeter of the standard chain
has perimeter greater than that of the (embedded) minimizer(s)
enclosing the same areas.

If instead the length of $\mathcal{B}_1 \cap \mathcal{B}_2$ is
greater than or equal to 2, then let $R_j$ be the region in the chain different from $R_i$, and let $A_j$ be the area of this region. Then $\mathcal{B}_1 \cap \mathcal{B}_2$ does not intersect the boundary of $R_j$. By the isoperimetric inequality for the plane, the
perimeter of the region is at least $2\sqrt{\pi A_j}$. Therefore,
the perimeter of the standard chain is at least $2 + 2\sqrt{\pi
A_j}$. By Corollary \ref{cor:interface two perimeter bound}, this
exceeds the perimeter of the minimizer(s) enclosing the same
areas.

\end{proof}

\subsection{Formulas for the Band Lens}\label{SS:Param_Band_Lens}
By Lemma \ref{lem:bl_short_direction}, any minimizing band lens
has the same homology as a shortest closed geodesic. Proposition
\ref{prop:area_and_perim_for band_lens} gives area and perimeter
formulas for these band lenses.

\begin{propsub}
\label{prop:area_and_perim_for band_lens} Consider a band lens
with the same homology as a shortest closed geodesic. Let $r =
(1/\mathrm{curvature})$ be the radius of curvature of the lens,
and let $d$ be the width of the band. The area and perimeter
formulas for the band lens in terms of $r$ and $d$ are
\begin{eqnarray*}
A_{\mathrm{Lens}}(r) &=& 2 r^2 \left(\frac{\pi}{3}-\frac{\sqrt 3}{4}\right), \\
A_{\mathrm{Band}}(r, d) &=& d - r^2 \left(\frac{\pi}{3}-\frac{\sqrt 3}{4}\right), \\
P(r, d) &=& \left(\frac{4 \pi}{3} - \sqrt 3\right) r + 2.
\end{eqnarray*}
Furthermore, $r < \frac{1}{\sqrt{3}}$ and $d > \frac{r}{2}$.
\end{propsub}

\begin{proof}
The chord between the vertices of the lens, with length $C =
2r\sin(\pi/3) = \sqrt{3}r$, must fit along a short direction of
the torus. Hence, $r < \frac{1}{\sqrt{3}}$. Also, because the band
and lens do not overlap, the width of the band is greater than
half the width of the lens. Hence, $d > r/2$.

Using Equations \ref{eqn:A}-\ref{eqn:L}, we can compute the area
and perimeter formulas:
\begin{eqnarray*}
A_{\mathrm{Lens}}(r) &=& 2 A\left(\frac{\pi}{3}, \sqrt 3 r\right) = 2 r^2
\left(\frac{\pi}{3}-\frac{\sqrt 3}{4}\right), \\
A_{\mathrm{Band}}(r, d) &=& d - \frac{1}{2} A_{\mathrm{Lens}}(r) = d - r^2
\left(\frac{\pi}{3}-\frac{\sqrt 3}{4}\right), \\
P(r, d) &=& 2 L \left(\frac{\pi}{3}, \sqrt 3 r\right) + (2 - \sqrt 3 r) = \left(\frac{4
\pi}{3} - \sqrt 3\right) r + 2.
\end{eqnarray*}
\end{proof}

\subsection{Formulas for the Standard Hexagon Tiling}\label{SS:Param_Standard_Hex_Tiling}

Corollary \ref{cor:standard_hex_has_p_three} showed that the
perimeter of every standard hexagon tiling is equal to three. This
section gives area formulas for the standard hexagon tiling
(Proposition \ref{prop:hex_three_parameterization}), proves that
at most one standard hexagon tiling encloses two given areas
(Lemma \ref{lem:uniqueness_by_algebra}), and describes exactly
which pairs of areas can be enclosed by a standard hexagon tiling
(Lemma \ref{lem:existence-of-some-hex-tilings}).

\begin{propsub} \label{prop:hex_three_parameterization}
In any standard hexagon tiling, there exist $a, b, c$ with sum 1
such that one region $R_1$ is a hexagon with side lengths
alternating between $a$ and $b$; another region $R_2$ is a hexagon
with side lengths alternating between $b$ and $c$; and the
exterior region $R_0$ is a hexagon with side lengths alternating
between $c$ and $a$. Conversely, given positive $a, b, c$ with sum
1, there exists a corresponding standard hexagon tiling.

The areas of the three regions are
\begin{eqnarray*}
A_{R_1} &=& \frac{\sqrt{3}}{4}(a^2 + 4ab + b^2), \\
A_{R_2} &=& \frac{\sqrt{3}}{4}(b^2 + 4bc + c^2), \\
A_{R_0} &=& \frac{\sqrt{3}}{4}(c^2 + 4ca + a^2).
\end{eqnarray*}

\end{propsub}

%%%%%%%%%%%%%%%%%%%%%%%%%%%%%%%%%%%%%%%%%%%%%%%%%%%%%%%%%%%%
\begin{figure}[ht]
\begin{center}
\begin{tabular} {@{\extracolsep{1.5 cm}} cc}
\includegraphics*[height=2in]{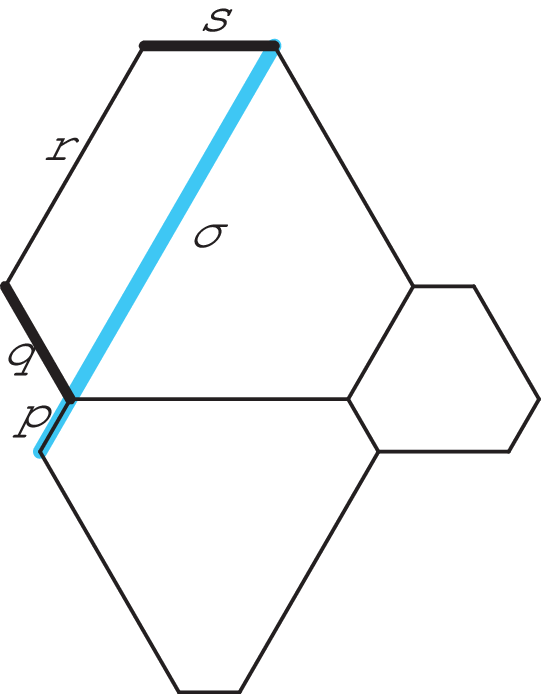}
&
\includegraphics*[height=2in]{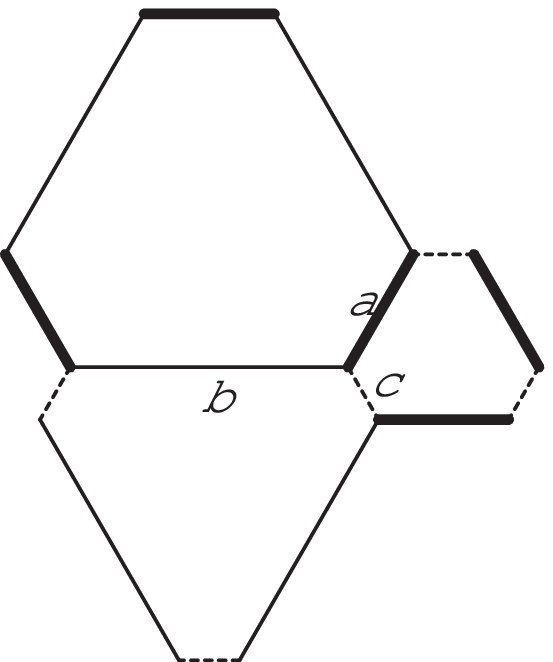}
\\
(a) & (b)
\end{tabular}
\end{center}
\caption{\label{fig:equil_triangle_on_standard_hex_tiling_fig} By
examining how far various edges of a standard hexagon tiling veer
from shortest closed geodesics, one can show that the sides of
each hexagon alternate between two different lengths.}
\end{figure}
%%%%%%%%%%%%%%%%%%%%%%%%%%%%%%%%%%%%%%%%%%%%%%%%%%%%%%%%%%%%

\begin{proof}
Consider the closed curve formed by the four segments $p, q, r, s$
depicted in Figure
\ref{fig:equil_triangle_on_standard_hex_tiling_fig}(a). Segments
$p$ and $r$ are parallel to a short direction of the torus, and
the segment $\sigma$ connecting the endpoints of the curve is
parallel to the same short direction of the torus. Therefore, the
segments $q$ and $s$ must veer the same distance away
from $\sigma$. Because they both meet $\sigma$ at $2\pi/3$, they
have equal length.

Just as $q$ and $s$ must be congruent, \textit{any} two sides of a hexagon
must be congruent if they are separated by one other side. It follows
that the sides of the hexagons have one of three lengths $a, b,
c$, as depicted in Figure
\ref{fig:equil_triangle_on_standard_hex_tiling_fig}(b). By
Corollary \ref{cor:standard_hex_has_p_three}, the tiling has
perimeter three, implying that $a + b + c = 1$.  Conversely, given
any positive $a, b, c$ such that $a + b + c = 1$, one can
construct the corresponding standard hexagon tiling as in Figure
\ref{fig:equil_triangle_on_standard_hex_tiling_fig}(b).

%%%%%%%%%%%%%%%%%%%%%%%%%%%%%%%%%%%%%%%%%%%%%%%%%%%%%%%%%%%%
\begin{figure}[ht]
\includegraphics*[height=1.5in]{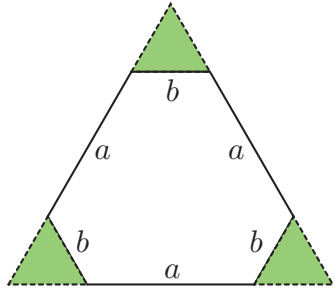}
\caption{\label{fig:extend_equil_triangles_fig} If the side
lengths of a hexagon with interior angles of $2\pi/3$ alternate
between two lengths, one can calculate the area of the hexagon by
extending three sides to form an equilateral triangle.}
\end{figure}
%%%%%%%%%%%%%%%%%%%%%%%%%%%%%%%%%%%%%%%%%%%%%%%%%%%%%%%%%%%%

To find the area of $R_1$, extend the sides of $R_1$ with length
$a$, as in Figure \ref{fig:extend_equil_triangles_fig}, in order
to form an equilateral triangle with sides of length $a + 2b$. The
new equilateral triangle consists of $R_1$ and three small
equilateral triangles with sides of length $b$. Therefore, the
area of $R_1$ equals
\[
\frac{\sqrt{3}}{{4}}(a+2b)^2 - 3 \cdot \frac{\sqrt{3}}{4}b^2 = \frac{\sqrt{3}}{4}(a^2 + 4ab + b^2).
\]
We find the analogous formulas for the areas of $R_2$ and $R_0$
similarly.
\end{proof}

\begin{lemsub} \label{lem:uniqueness_by_algebra} Given positive real numbers $A_0, A_1, A_2$, there exists at most one triple $(a, b, c)$ of nonnegative numbers with sum 1 such that
\begin{eqnarray*}
A_1 &=& \frac{\sqrt{3}}{4}(a^2 + 4ab + b^2), \\
A_2 &=& \frac{\sqrt{3}}{4}(b^2 + 4bc + c^2),\\
A_0 &=& \frac{\sqrt{3}}{4}(c^2 + 4ca + a^2).
\end{eqnarray*}
\end{lemsub}

\begin{proof} Suppose, for sake of contradiction, that there existed two such triples $(a_1$, $b_1$, $c_1)$ and $(a_2, b_2, c_2)$. Without loss of generality, assume that $a_1 < a_2$. Because $a_1^2 + 4a_1b_1 + 4b_1^2 = a_2^2 + 4a_2b_2 + b_2^2$, it follows that $b_1 > b_2$. Similarly, because $b_1^2 + 4b_1c_1 + 4c_1^2 = b_2^2 + 4b_2c_2 + 4c_2^2$, it follows that $c_1 < c_2$. Likewise, $c_1 < c_2$ implies that $a_1 > a_2$, a contradiction. \end{proof}

\begin{corsub} \label{cor:uniqueness_of_standard_hex_tiling} Given prescribed areas $A_1, A_2$, there exists at most one standard hexagon tiling enclosing those areas. \end{corsub}

\begin{proof}
The result follows from our parameterization (Proposition
\ref{prop:hex_three_parameterization}) and Lemma
\ref{lem:uniqueness_by_algebra}.
\end{proof}

\begin{lemsub} \label{lem:existence-of-some-hex-tilings} Given positive areas $A_0, A_1, A_2$ with sum $\sqrt{3}/2$ (the area of the hexagonal torus), there exists a standard hexagon tiling enclosing any two of those areas if and only if $\frac{2}{\sqrt[4]{3}}\left(\sqrt{A_j} + \sqrt{A_k}\right) > 1$ for $0 \le j < k \le 2$ --- that is, if and only if the side length of an equilateral triangle with area $A_j$ and the side length of an equilateral triangle with area $A_k$ have a sum greater than 1. (When equality holds for some pair, say $A_1, A_2$, a degenerate hexagon tiling exists --- consisting of two equilateral triangles enclosing areas $A_1, A_2$ and with an exterior hexagonal region of area $A_0$.)
\end{lemsub}

%%%%%%%%%%%%%%%%%%%%%%%%%%%%%%%%%%%%%%%%%%%%%%%%%%%%%%%%%%%%
\begin{figure}[ht]
\begin{center}
\begin{tabular} {c}
\includegraphics*[height=1.75in]{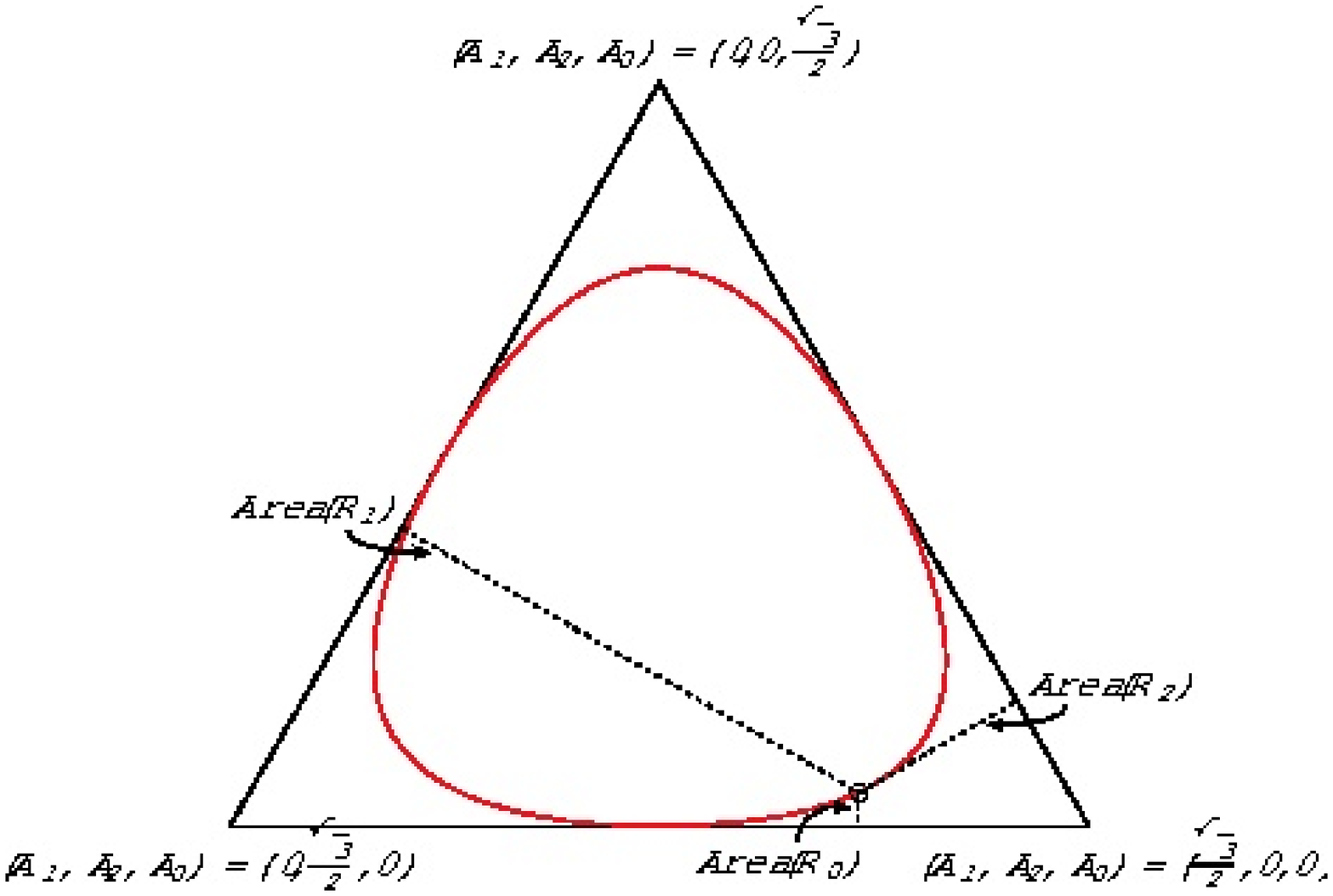} \\
(a) \\[12pt]
\includegraphics*[height=1in]{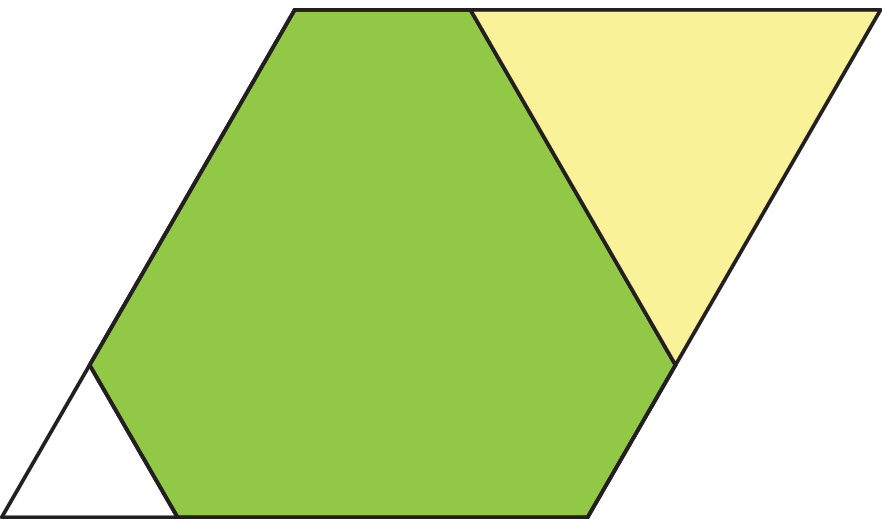}
\\
(b)
\end{tabular}
\end{center}
\caption{\label{fig:hex_tile_exist_curve_fig} As in the phase portraits in Figure \ref{fig:param-diags}, the triangle in (a) represents the possible area pairs $(A_1, A_2)$ on the hexagonal torus, where the areas of the three regions are given by
the distances to the three sides. The region bounded by the curve represents exactly those area pairs enclosed by some standard hexagon tiling. The curve itself represents area pairs enclosed by degenerate hexagon tilings consisting of a hexagon and
two triangles, such as the one depicted in (b).}
\end{figure}
%%%%%%%%%%%%%%%%%%%%%%%%%%%%%%%%%%%%%%%%%%%%%%%%%%%%%%%%%%%%

\begin{proof} Consider the following map from $\{(a, b) \mid 0 \le a, b \textrm{\ and\ } a + b \le 1\}$ to
$\mathbb{R}^2$,
\[
(a,b) \mapsto \left( \frac{\sqrt{3}}{4}(a^2 + 4ab + b^2),
\frac{\sqrt{3}}{4}(b^2 + 4b(1-a-b) + (1-a-b)^2) \right).
\]
By our parameterization (Proposition
\ref{prop:hex_three_parameterization}), this map sends a pair of
side lengths of a hexagon tiling to the areas of the region with
side lengths $a, b$ and of the region with side lengths $b,
1-a-b$. Thus, as we vary $a$ and $b$, we obtain all ordered pairs
of areas that can be enclosed by a standard hexagon tiling.

By Lemma \ref{lem:uniqueness_by_algebra}, the above map is
injective. Also, it maps the boundary of its domain continuously
to a simple closed curve $\mathcal{C}$ in $\mathbb{R}^2$, as
depicted in Figure \ref{fig:hex_tile_exist_curve_fig}(a). (The
figure is drawn with barycentric coordinates, not Cartesian
coordinates.) Because the map is injective and continuous, it
follows that its range is precisely the closed region bounded by
$\mathcal{C}$. This boundary $\mathcal{C}$ consists of three
pieces, corresponding to $a = 0, b = 0, c = 0$. When $a = 0$, we
have $1 = b + c = \frac{2}{\sqrt[4]{3}}(\sqrt{A_1} + \sqrt{A_2})$;
similarly, along the other two pieces of the boundary, we have $1
= \frac{2}{\sqrt[4]{3}}(\sqrt{A_j} + \sqrt{A_k})$ for $(j, k) =
(2, 0)$ and $(0, 1)$. Hence, the interior of the region bounded by
these three pieces consists of those area triples $(A_0, A_1,
A_2)$ such that $\frac{2}{\sqrt[4]{3}}(\sqrt{A_j} + \sqrt{A_k}) >
1$ for $0 \le j < k \le 2$. Therefore, for each of these triples
and no others, some standard hexagon tiling encloses some two of the
three areas.
\end{proof}

\begin{lemsub} \label{lem:hex-tiling-wins-if-we-think-so} If a minimizing double bubble enclosing areas $A_1,
A_2$ on the hexagonal torus has perimeter equal to three, then there exists a
(perimeter-minimizing) standard hexagon tiling enclosing those
areas.
\end{lemsub}

\begin{proof}
Let the exterior region of the given double
bubble have area $A_0$, and without loss of generality assume that
$A_1 \le A_2 \le A_0$.  By Lemma \ref{lem:band-lens-perimeter-bound},
\[
2 + \sqrt{A_1 \left(\frac{8\pi}{3}-2\sqrt{3}\right)}
\]
is greater than or equal to the given minimal perimeter, three.
The (non-sharp) bound $A_1 \ge \frac{4}{25}$ follows easily.
Hence, $A_0 \ge A_2 \ge A_1 \ge \frac{4}{25}$. For $0 \le j < k
\le 2$, we thus have
\[
\sqrt{A_j} + \sqrt{A_k} \ge \frac{4}{5} \ge \frac{\sqrt[4]{3}}{2}.
\]
By Lemma \ref{lem:existence-of-some-hex-tilings}, there exists a
standard hexagon tiling enclosing areas $A_1, A_2$. By Corollary
\ref{cor:standard_hex_has_p_three}, the tiling has perimeter
three, so it is perimeter minimizing. \end{proof}

\section{Phase Portraits} \label{S: numerical comparisons}

The phase portraits of Figure \ref{fig:param-diags} show which of
the minimizers of Theorem \ref{thm:main theorem} win for which
pairs of prescribed areas on four different tori --- where we
represent each torus as a parallelogram with side lengths $1, L$
and interior angle $\theta$. To create these phase portraits, we
plot perimeter as a function of area pairs for each candidate in
Theorem \ref{thm:main theorem} using the parameterizations from
Section \ref{S:perimeter-and-area}. For each candidate in
\ref{thm:main theorem} besides the standard chain, we know that at
most one double bubble exists enclosing any given pair of areas,
so we can solve for its perimeter numerically. Because we do not
know that at most one standard chain (of a fixed homology) exists
enclosing two given areas, we instead create the perimeter-areas
plot for the standard chain parametrically. Also, using Lemma
\ref{lem:wrap only once}, it is possible to show that any
minimizing standard chain has axis-length $1$, $L$, or $\sqrt{1 +
L^2 \pm 2\cos\theta}$, so we do not plot perimeter as a function
of area pairs for other types of standard chains.

\begin{figure}[htbp]
\begin{center}
\includegraphics*[height=3.25in]{\Path/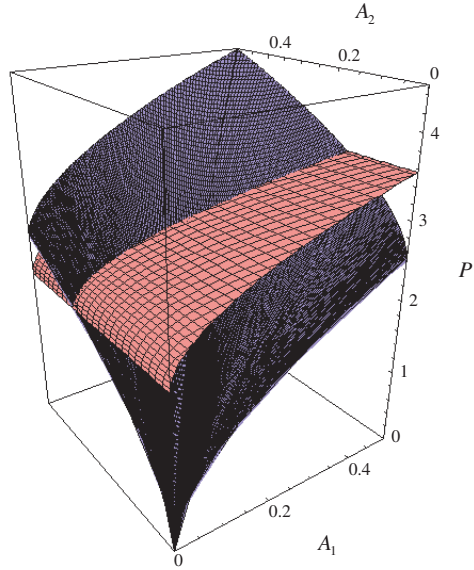}%
\end{center}
\caption{\label{fig:finalplot} This \emph{Mathematica} plot of the
intersection of the perimeter plots for the standard double bubble and band lens illustrates the process used to create the phase portraits in Figure \ref{fig:param-diags}. The plots suggest that perimeters do not fluctuate wildly along certain sets of unusual area pairs, increasing our confidence that the phase portraits are accurate.}
\end{figure}

After plotting all these perimeter-areas plots simultaneously, we
then view the combined plot from below so that at each area
pair, the plot(s) of the minimal calculated perimeter is
visible. Figure \ref{fig:finalplot} illustrates a
simplified version of this process, comparing the perimeters of
only two candidates, the standard double bubble and the band lens.

This process works if whenever the perimeter calculated for a given double bubble is minimizing for two given areas, the double bubble actually exists on the torus for those areas. Indeed, this is the case. When the standard double bubble does not fit on the torus, it has diameter greater than 1, so by Propositions \ref{prop:perimeter is 3
times diameter} and \ref{prop:perimeter_bound} it cannot have
minimizing perimeter.  Lemma \ref{lem:band-lens-perimeter-bound}
shows that when the calculated perimeter for the band lens is
minimizing, that band lens exists. By Lemma \ref{lem:symm chain fits},
when a standard chain with axis-length one has minimizing
perimeter, it can be embedded on the torus; and according to the plots, no standard chain with axis-length greater than one has minimizing perimeter. As the simplest case,
we know that the double band exists for all area pairs.  Finally, by Lemma \ref{lem:hex-tiling-wins-if-we-think-so}, when the minimal perimeter on the hexagonal torus is 3 for two given areas, there exists a standard hexagon tiling enclosing those two areas. \\

%%%%%%%%%%%%%%%%%%%%%%%%%%%%%%%%%%%%%%%%%%%%%%%%%%%%%%%%
%%%%%%%%%%%%%%%%%%%%%%%%%%%%%%%%%%%%%%%%%%%%%%%%%%%%%%%%

The following conjecture describes a simple connection \textit{among}
the various phase portraits. As usual, we assume that the shortest
geodesics on the torus and the cylinder have unit length.

\begin{figure}[htbp]
\begin{center}
\includegraphics*[height=2.25in]{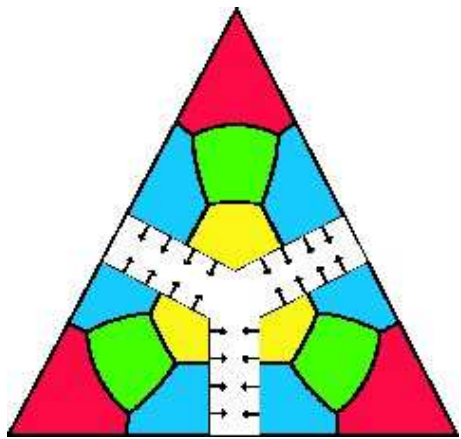}%
\end{center}
\caption{\label{fig:cutout} Removing pieces of this phase portrait
for a rectangular torus with length $1.2$ yields the phase
portrait plot for some torus with less area.}
\end{figure}

\begin{conj}
Given areas $A_0 \ge A_2 \ge A_1$ and a torus with
area $A_0 + A_1 + A_2$, the non-tiling minimizing double bubble(s)
enclosing $A_1, A_2$ on this torus are the same as the minimizing double
bubble(s)  enclosing $A_1, A_2$ on the infinite cylinder. Consequently the
phase portrait of any flat two-torus can be obtained from the phase
portrait of any flat two-torus of larger area (or from the phase portrait
of the infinite cylinder) as in Figure \ref{fig:cutout}. In particular,
tori of equal areas have identical phase portraits.
\end{conj}

To prove Conjecture 8.1 it would suffice to show that every minimizing
standard chain has axis-length one. The proof would also use the fact
(which we can show using Theorem \ref{thm:main theorem}) that in a
minimizer, a region of greater pressure has no more area than a region of
less pressure, as well as Proposition
\ref{prop:perimeter is 3 times diameter}, Lemma
\ref{lem:bl_short_direction}, and Lemma \ref{lem:symm chain fits}.

%%%%%%%%%%%%%%%%%%%%%%%%%%%%%%%%%%%%%%%%%%%%%%%%%%%%%%%%
%%%%%%%%%%%%%%%%%%%%%%%%%%%%%%%%%%%%%%%%%%%%%%%%%%%%%%%%

\end{document}